\documentclass[11pt,leqno,oneside,letterpaper, openany]{amsart}
\usepackage[letterpaper,left=1.8truein, top=1.1truein]{geometry}
\usepackage{amsmath,amsfonts,amssymb,amscd,amsthm,amsbsy,epsf,graphics} 
\usepackage{stmaryrd}
\usepackage{setspace}
\usepackage[title]{appendix}
\usepackage{comment} 
\usepackage[normalem]{ulem}
\usepackage{mathrsfs}
\usepackage{enumitem} 

\usepackage[dvipsnames]{xcolor}
\usepackage{color}

\definecolor{org}{RGB}{225,119,0}

\usepackage[colorlinks=true,linkcolor=blue,citecolor=blue, urlcolor=black]{hyperref} 

\usepackage{bm}
\usepackage{caption,subcaption,graphicx,epsfig}
\usepackage{tikz-cd,calc}
\usetikzlibrary{knots}
\usetikzlibrary{decorations.markings}

\usepackage{lipsum}

\makeatletter
\def\l@subsection{\@tocline{2}{0pt}{4pc}{5pc}{}}
\makeatother
 
\let\oldtocsection=\tocsection
 
\let\oldtocsubsection=\tocsubsection
 
\let\oldtocsubsubsection=\tocsubsubsection
 
\renewcommand{\tocsection}[2]{\hspace{0em}\oldtocsection{#1}{#2}}
\renewcommand{\tocsubsection}[2]{\hspace{0em}\oldtocsubsection{#1}{#2}}
\renewcommand{\tocsubsubsection}[2]{\hspace{2em}\oldtocsubsubsection{#1}{#2}}

\DeclareGraphicsExtensions{.pdf}

\newtheorem{theorem}{Theorem}[section] 
\newtheorem*{theorem*}{Theorem}
\newtheorem{corollary}[theorem]{Corollary}
\newtheorem*{corollary*}{Corollary}
\newtheorem{lemma}[theorem]{Lemma}
\newtheorem{proposition}[theorem]{Proposition}
\newtheorem*{proposition*}{Proposition}
\newtheorem{conjecture}[theorem]{Conjecture}

\newtheorem*{question*}{Question}

\newtheorem{innercustomthm}{Theorem}
\newenvironment{customthm}[1]
  {\renewcommand\theinnercustomthm{#1}\innercustomthm}
  {\endinnercustomthm}

\theoremstyle{definition}
\newtheorem{definition}[theorem]{Definition} 
\newtheorem{example}[theorem]{Example}

\newtheorem*{problem*}{Problem}

\newtheorem{remark}[theorem]{Remark} 
\newtheorem*{remark*}{Remark}

\newtheorem*{acknowledgement*}{Acknowledgements}

\newtheoremstyle{cases}
  {12pt plus 6 pt}
  {2pt}
  {\bfseries}   
  {}
  {\bfseries}
  {.}
  {.5em}
  {}
\theoremstyle{cases}

\numberwithin{subcase}{case} 
\numberwithin{subsubcase}{subcase}
\numberwithin{equation}{subsection}

\graphicspath{{figures/}}

\begin{document}

\title{Slope detection and toroidal $3$-manifolds}

\author[Steven Boyer]{Steven Boyer}
\thanks{Steven Boyer was partially supported by NSERC grant RGPIN-2018-06549}
\address{D\'epartement de Math\'ematiques, Universit\'e du Qu\'ebec \`a Montr\'eal, 201 President Kennedy Avenue, Montr\'eal, Qc., Canada H2X 3Y7.}
\email{boyer.steven@uqam.ca}
\urladdr{http://www.cirget.uqam.ca/boyer/boyer.html}

\author[Cameron McA. Gordon]{Cameron McA. Gordon} 
\address{Department of Mathematics, University of Texas at Austin, 1 University Station, Austin, TX 78712, USA.}
\email{gordon@math.utexas.edu}

\author[Ying Hu]{Ying Hu}
\address{Department of Mathematical and Statistical Sciences, University of Nebraska Omaha, 6001 Dodge Street, Omaha, NE 68182-0243, USA.}
\email{yinghu@unomaha.com}
\urladdr{https://yinghu-math.github.io}

\thanks{2010 Mathematics Subject Classification. Primary 57M25, 57M50, 57M99}

\thanks{Key words: $L$-space conjecture, toroidal manifolds, slope detection, branched surfaces, universal circle actions, homology spheres, satellite knots, left-orderable, cyclic branched covers}

\begin{abstract}
The $L$-space conjecture asserts the equivalence, for prime 3-manifolds, of three properties: not being an $L$-space, having a left-orderable fundamental group, and admitting a co-oriented taut foliation. We investigate these properties for toroidal $3$-manifolds using various notions of slope detection. Our main technical result gives sufficient conditions for certain slopes on the boundaries of rational homology solid tori to be detected by left-orders, foliations, and Heegaard Floer homology, using Thurston's universal circle actions, Li's result on laminar branched surfaces, and Rasmussen-Rasmussen's result on L-space intervals, respectively. This leads to a proof that toroidal integer homology spheres have left-orderable fundamental groups, as predicted by the $L$-space conjecture. It also allows us to show that the cyclic branched covers of prime satellite knots are not $L$-spaces and have left-orderable fundamental groups, as conjectured by Gordon and Lidman. Similarly we show that a cyclic branched cover of a satellite knot admits a co-oriented taut foliation when it has a fibred companion. A partial extension of these results to toroidal links leads to a proof that prime quasi-alternating links are either hyperbolic or $(2, m)$-torus links, which generalises Menasco's classical theorem that non-split alternating links are either hyperbolic or $(2, m)$-torus links.
\end{abstract}

\maketitle 

 \setcounter{tocdepth}{1}
{\footnotesize

\parskip=.2em
\tableofcontents
} 

\newpage
\section{Introduction}
\label{sec: intro}
The $L$-space conjecture contends that a closed, connected, orientable, irreducible $3$-manifold is not an $L$-space if and only if it admits a co-oriented taut foliation and if and only if it has a left-orderable fundamental group; see \cite[Conjecture 1]{BGW} and  \cite[Conjecture 5]{Juh}. We will sometimes abbreviate these properties by saying that the manifold is $NLS$, $CTF$, or $LO$, respectively. The only implication known is that an $L$-space admits no co-oriented taut foliations, i.e. $CTF$ implies $NLS$ (\cite{OS1, KR, Bow}). Though the conjecture has been verified in a number of cases, including all graph manifolds (\cite{BRW, LS, BGW, BC1, HRRW}), it remains widely open. Here we investigate the three properties in question in the case of toroidal $3$-manifolds. 

One feature of this article is the unified approach to the study of the three properties in the $L$-space conjecture based on various notions of slope detection, which reduces the verification that a toroidal manifold has one of these properties to the question of whether certain slopes are detected in the appropriate sense. Some of the main contributions of this paper concern slope detection, but before describing them and our motivation for proving them, we list three applications of the method; for a fuller account of our applications, see \S \ref{sec: results}.

Eftekhary \cite{Eft} (see also \cite{HRW1}) has shown that irreducible toroidal integer homology 3-spheres are $NLS$, so the $L$-space conjecture predicts that they are also $LO$ and $CTF$. We show

\begin{customthm}{\ref{thm: tor zhs implies lo}}
An irreducible toroidal integer homology 3-sphere has a left-orderable fundamental group. It admits a co-oriented taut foliation if it contains an incompressible torus bounding a 3-manifold that fibres over the circle.
\end{customthm}

A class of 3-manifolds for which the three properties $NLS$, $LO$ and $CTF$ have been widely investigated are the $n$-fold cyclic branched covers $\Sigma_n(K)$ of knots $K$. For torus knots the situation is completely understood \cite{GLid1}, and there are several results for various families of hyperbolic knots; see the references in \S \ref{sec: results}. For the remaining case of prime satellite knots $K$, Gordon and Lidman conjectured (\cite{GLid1, GLid2}) that $\Sigma_n(K)$ is $LO$ and $CTF$ (and hence $NLS$) for all $n \ge 2$. We prove this for $LO$ and $NLS$, and obtain a partial result for $CTF$. In fact we prove the corresponding statements for prime satellite links.

\begin{customthm}{\ref{thm: GL thm intro}}
 If   $L$ is a prime satellite  link then for all $n \ge 2$, $\Sigma_n(L)$ is not an $L$-space and has a left-orderable fundamental group. If in addition   $L$ has a fibred companion knot, then $\Sigma_n(L)$ admits a co-oriented taut foliation for all $n \geq 2$.
\end{customthm}

A third application is to quasi-alternating links: we show that Menasco's classical theorem, that a prime non-split alternating link is either hyperbolic or a $(2,m)$-torus link, extends to quasi-alternating links. (For knots this exension is called a ``folklore conjecture'' in \cite{Dey}.)

\begin{customthm}{\ref{thm: satellite not qa}}
A prime $\mathbb Z/2$-Khovanov thin link is either hyperbolic or a $(2, m)$-torus link. In particular, this holds for prime quasi-alternating links.  
\end{customthm}

These results follow from the key technical contributions of the paper, which are to provide sufficient conditions for the detectability of certain slopes and to prove a gluing theorem for co-oriented taut foliations. One of the main motivations for this work is to lay the groundwork for a systematic approach to proving the $L$-space conjecture for toroidal $3$-manifolds. We discuss the details next. 

Slope detection was first introduced by Boyer and Clay in the case of Seifert fibred manifolds to study the $L$-space conjecture for toroidal graph manifolds (\cite{BC1}). Roughly speaking, certain slopes on the boundary of a knot manifold\footnote{A knot manifold is a compact, connected, orientable, irreducible $3$-manifold with boundary an incompressible torus.} are singled out (i.e. ``detected") using Heegaard Floer homology ($NLS$-detection), left-orders ($LO$-detection), or foliations ($CTF$-detection). The importance of slope detection can be seen in the fact that a closed $3$-manifold $W$ will have property $* = NLS, LO$ or $CTF$ if it can be expressed as a union of two knot manifolds $W = M_1 \cup M_2$ in such a way that the gluing map matches a $*$-detected slope on $\partial M_1$ with a $*$-detected slope on $\partial M_2$. See Theorem \ref{thm: * gluing} below. 

More precisely, a rational slope on a torus $T$ is an isotopy class of essential simple closed curves on $T$ (cf. \S \ref{sec: conventions}). We say that a rational slope $[\alpha]$ on the boundary of a knot manifold $M$ is 

\begin{itemize}[leftmargin=*]
\setlength\itemsep{0.5em}
\item {\it $NLS$-detected} if it lies in the closure of the set of rational slopes whose associated Dehn fillings are not $L$-spaces (cf. \S \ref{subsec: NLS detn}). 

\item {\it $CTF$-detected} (or {\it foliation-detected)} if there is a co-oriented taut foliation on $M$ which intersects $\partial M$ transversely in a foliation without Reeb annuli which contains a closed leaf of slope $[\alpha]$ (cf. \S \ref{subsec: slopes of folns}). 
 
\item {\it $LO$-detected} (or {\it order-detected)} if there is a left-order $\mathfrak{o} \in LO(M)$ such that for each $\gamma \in \pi_1(M)$, the slope of $(\gamma \cdot \mathfrak{o})|_{\pi_1(\partial M)}$  is $[\alpha]$ (cf. \S \ref{subsec: gens on los} and \S \ref{subsec: lo z2 and o-d}). 
\end{itemize}
Though irrational slopes (cf. \S \ref{sec: conventions}) arise naturally in the setting of $LO$-detection and $CTF$-detection, and will be discussed in \S \ref{sec: foln gluing knot mflds} and \S \ref{sec: ord gps} respectively, our results and their proofs depend only on rational slopes, and that is where our focus will lie. 

There are topological characterisations of $*$-detected rational slopes in terms of attaching so-called $*$-generalised solid tori. A simple example of a non-hyperbolic $*$-generalised solid torus is the twisted $I$-bundle over the Klein bottle. See Corollary \ref{cor: gluing defn of nls detection}, Corollary \ref{cor: gluing defn of ctf detection} and Corollary \ref{cor: gluing defn of lo detection}.  

The $L$-space conjecture suggests the equality of the sets of $LO$-detected, $CTF$-detected, and $NLS$-detected rational slopes on $\partial M$ (cf. Remark \ref{rem: detd slopes the same}), though to date, the only known relation is that the set of rational slopes which are $CTF$-detected is contained in the set of $NLS$-detected   rational slopes; see Corollary \ref{cor: ctf detd implies nls detd}.  Moreover, owing to the work of Jake and Sarah Rasmussen \cite{RR} and of Jonathan Hanselman, Jake Rasmussen and Liam Watson \cite{HRW1}, the set of $NLS$-detected   rational slopes consists of either all  rational slopes, a closed subinterval of  rational slopes, or a singleton containing only the longitudinal slope (\S \ref{sec: $NLS$ detn}).

The following gluing theorem is a combination of work of Hanselman, Rasmussen and Watson \cite[Theorem 13]{HRW1}, Boyer and Clay \cite[Theorem 1.3]{BC2}, and  Theorem \ref{thm: fln gluing} of this paper. 

\begin{theorem}  
\label{thm: * gluing}
Suppose that $W = M_1 \cup_f M_2$ where $M_1, M_2$ are knot manifolds and $f: \partial M_1 \xrightarrow{\; \cong \; } \partial M_2$. 
If $f$ identifies a rational $\ast$-detected slope on the boundary of $M_1$ with a rational $\ast$-detected slope on the boundary of $M_2$, then $W$ has property $\ast$. 
\end{theorem}

The same conclusion also holds when the closed $3$-manifold is obtained by gluing several manifolds together along multiple tori. See Theorem \ref{thm: general * gluing}. The converse of Theorem \ref{thm: * gluing} is known to be true when $* = NLS$ by \cite[Theorem 13]{HRW1}. We conjecture the following: 

\begin{conjecture}
\label{conj: gluing}
If the union $W = M_1 \cup M_2$ of two knot manifolds is $* \in \{LO, CTF\}$, then the gluing map identifies a rational $\ast$-detected slope on the boundary of $M_1$ with a rational $\ast$-detected slope on the boundary of $M_2$. 
\end{conjecture}

The major difficulty in applying the $*$-gluing theorem is to show that certain  rational slopes are $*$-detected. Our main detection results are consolidated in Theorem \ref{thm: meridional detn} below. 

Recall that the genus $g(M)$ of an irreducible rational homology solid torus $M \not \cong S^1 \times D^2$ whose longitude  $\lambda_M$ is integrally null-homologous is the minimal genus of an oriented surface with connected boundary representing a generator of $H_2(M, \partial M) \cong \mathbb Z$. Let $T_1(M)$ denote the torsion subgroup of $H_1(M)$. 
 
\begin{theorem} 
\label{thm: meridional detn} 
Let $M \not \cong S^1 \times D^2$ be an irreducible rational homology solid torus whose longitude  $\lambda_M$ is integrally null-homologous. 
\begin{enumerate}[leftmargin=*] 
\setlength\itemsep{0.5em}
\item[{\rm (1)}]  Each  rational slope of distance $1$ from $\lambda_M$ is NLS-detected. Further, if $M$ is an integer homology solid torus, then each rational slope whose distance   from $\lambda_M$ divides $2g(M) - 1$ is NLS-detected. 
\item[{\rm (2)}] If $T_1(M)$ is a $\mathbb Z/2$ vector space, then each rational slope whose distance   from $\lambda_M$ divides $2g(M) - 1$ is $LO$-detected. In particular this holds for  rational slopes of distance $1$ from $\lambda_M$. 
\item[{\rm (3)}] If $M$ fibres over the circle, then each rational slope of distance $1$   from 
$\lambda_M$ is $CTF$-detected.
\end{enumerate}
\end{theorem}

Theorem \ref{thm: meridional detn}(2) has the following useful corollary. 

\begin{corollary}
\label{cor: ord detd fibration}
If $K$ is a non-trivial knot in the $3$-sphere, each rational slope whose distance   from $\lambda_M$ divides $2g(M) - 1$ is $LO$-detected. In particular, the meridional slope of $K$ is $LO$-detected, as are the slopes $2g(K) - 1$ and $-(2g(K) - 1)$. 
\end{corollary}

Parts (1), (2) and (3) of Theorem \ref{thm: meridional detn} are proven in Sections \ref{sec: $NLS$ detn}, \ref{sec: univ circles and knot mflds} and \ref{sec: ctf on hyperbolic fibred knot} respectively. In particuar, Theorem \ref{thm: meridional detn}(1) is a corollary of a more general $NLS$-detection result stated in Proposition \ref{prop: mu + k lambda NLS-detd}.

The key step in the proof of Theorem \ref{thm: meridional detn}(2), is to show that the universal circle action of $\pi_1(M)$ associated to a finite depth foliation on a knot manifold $M$ satisfies a certain dynamical property when restricted to the peripheral subgroup $\pi_1(\partial M)$, which is stated in Theorem \ref{thm: universal circle fixed point} below and is interesting in its own right. In this theorem,  
$\mbox{{\rm Homeo}}_{\mathbb Z}(\mathbb R) = \{f \in \mbox{{\rm Homeo}}_{+}(\mathbb R)  \; | \; f(x + 1) = f(x) + 1\}$
denotes the universal covering group of $\mbox{{\rm Homeo}}_+(S^1)$, where $f \in \mbox{{\rm Homeo}}_{\mathbb Z}(\mathbb R)$ maps to $f$ (mod $\mathbb Z$) in $\mbox{Homeo}_+(\mathbb R / \mathbb Z) = \mbox{Homeo}_+(S^1)$.  

\begin{theorem} 
\label{thm: universal circle fixed point}
Let $M \not \cong S^1 \times D^2$ be an irreducible rational homology solid torus whose longitude  $\lambda_M$ is integrally null-homologous. 
\begin{enumerate}[leftmargin=*] 
\setlength\itemsep{0.5em}
\item[{\rm (1)}]  There exists an associated universal circle action $\rho_M: \pi_1(M) \rightarrow \mbox{{\rm Homeo}}_+(S^1)$ whose restriction to $\pi_1(\partial M)$ has a global fixed point;
\item[{\rm (2)}] Suppose that $\widetilde \rho_M: \pi_1(M) \rightarrow \mbox{{\rm Homeo}}_{\mathbb Z}(\mathbb R)$ is a lift of $\rho_M$. Then the absolute value of the translation number of $\widetilde{\rho}_M(\lambda_M)$ is $2g(M) - 1$, where $g(M)$ is the genus of $M$. 
\end{enumerate}
\end{theorem}
The main problem left unresolved by our paper is in establishing detection results for $CTF$ comparable to those we prove for $LO$ and $NLS$.

\begin{conjecture}
\label{conj: ctf detection}
Let $M \not \cong S^1 \times D^2$ be an irreducible rational homology solid torus whose longitude  $\lambda_M$ is integrally null-homologous. Then each rational slope of distance $1$ from 
$\lambda_M$ is $CTF$-detected. 
\end{conjecture}

Assuming that Conjecture \ref{conj: ctf detection} holds, we can then remove the the fibring restrictions from Theorems \ref{thm: tor zhs implies lo} and \ref{thm: GL thm intro}. 

To date, the most well-developed method for constructing taut foliations is to use branched surfaces, due to the work of Tao Li in \cite{Li1, Li2}. However, building a branched surface in a general knot manifold which carries all meridional slopes is challenging, even for knot exteriors in $S^3$ (cf.~\cite{LR}). Theorem \ref{thm: meridional detn}(3) shows that Conjecture \ref{conj: ctf detection} holds when the knot manifold $M$ is fibred. This is done in \S \ref{sec: ctf on hyperbolic fibred knot} by adapting Tao Li's arguments in \cite{Li1} for constructing taut foliations on closed manifolds.

We end this introduction with a brief discussion of the subtle relationship between slope detection and Dehn filling. 

As mentioned above, there are topological characterisations of $*$-detected rational slopes in terms of attaching twisted $I$-bundles over the Klein bottle, though at first glance it may seem more natural to try and characterise them through the attachment of solid tori, i.e. by Dehn filling. Indeed, for $* = NLS$ and $LO$, it is true that if the corresponding Dehn filled manifold has property $*$, then the filling slope is $*$-detected. (For $* = CTF$, this is true with the a priori stronger notion of strong $CTF$-detection; see \S \ref{sec: foln gluing knot mflds}.) However, the following example shows that it is necessary to consider a broader set of rational slopes when studying when a manifold obtained by gluing together two knot manifolds along their boundaries has property $*$. 

\begin{example}     
\label{example: trefoil glue}
For each $n\geq 2$, let  $W_n$ be the toroidal graph manifold obtained by gluing two copies of the exterior of the  right-handed trefoil knot through the map $f_n$ which sends slope $r$ to $-r+n$ for any $r\in \mathbb{Q}$. It is well-known that the set of rational slopes for which the corresponding Dehn fillings of the trefoil exterior has a left-orderable fundamental group, resp. admits a co-oriented taut foliation, is precisely $(-\infty, 1)$. Under the gluing map $f_n$, the interval $(-\infty, 1)$ is sent to the interval $(n-1,\infty)$, so no rational slope in $(-\infty, 1)$ is matched to a  rational slope in $(-\infty, 1)$. On the other hand, $W_n$ has a left-orderable fundamental group and admits a co-oriented taut foliation for each $n\geq 2$ \cite{BC1}.  In fact, the meridional slope of the trefoil exterior is $\ast$-detected for each value of $\ast$, and since two copies of the meridional slope are identified under the gluing maps $f_n$, it follows from the gluing property of $\ast$-detected  rational slopes, Theorem \ref{thm: * gluing}, that $W_n$ has property $\ast$. 
\end{example}

Though Example \ref{example: trefoil glue} shows that the set of rational slopes whose associated Dehn fillings have property $*$ can be properly contained in the set of $*$-detected rational slopes, it is expected that the former is the set of rational slopes in the closure of the latter. This is known to hold when $* = NLS$, and so the conjectured independence of the set of $*$-detected rational slopes from $*$ (see Remark \ref{rem: detd slopes the same})  combines with the $L$-space conjecture to suggest that it is also true for $* = LO$ and $CTF$.

\subsection{Organisation of the paper} 
We state our applications of the slope detection method in Section \ref{sec: results}. In Section \ref{sec: conventions}, we list notation, terminology, and certain conventions on slopes used in the paper.
Sections \ref{sec: $NLS$ detn}, \ref{sec: foln gluing knot mflds}, and \ref{sec: ord gps} discuss $\ast$-detection of  rational slopes, where $\ast$ is $NLS, CTF$, and $LO$ respectively, and the knot manifold $\ast$-gluing theorems (Theorems \ref{thm: HRW1}, \ref{thm: fln gluing}, and \ref{thm: LO-gluing}). Consequently we deduce Corollaries \ref{cor: gluing defn of nls detection},  \ref{cor: gluing defn of ctf detection}, and \ref{cor: gluing defn of lo detection}, which characterise $\ast$-detection of rational slopes in terms of attaching $*$-generalised solid tori. Then, in \S \ref{sec: detn and gluing}, we state and prove the general $\ast$-gluing theorem (Theorem \ref{thm: general * gluing}) for gluing $3$-manifolds along multiple torus boundary components. 
Universal circle representations are discussed in \S \ref{sec: univ circles and knot mflds} and applied to prove Theorem \ref{thm: universal circle fixed point} and consequently, Theorem \ref{thm: meridional detn}(2). Proposition \ref{prop: (d) implies not ord-det} and Corollary \ref{cor: nie 1bb} are also proved, which illustrate the sharpness of our $LO$-detection results.
Section \ref{sec: ctf on hyperbolic fibred knot} is devoted to building  foliations that detect  rational slopes of distance $1$ from the longitudinal slope of a fibred knot manifold (Theorem \ref{thm: meridional detn}(3)). The final section of the paper, \S \ref{sec: cbc and (*)}, applies the material developed in the earlier sections to prove the results described in \S \ref{sec: results}.

\begin{acknowledgement*}
We thank Tye Lidman for helpful remarks concerning Heegaard Floer homology and Rachel Roberts for pointing out the relevance of Tao Li's work to foliation-detection.  
\end{acknowledgement*}

\section{Applications of slope detection results to toroidal manifolds}
In this section, we list the applications of our slope detection results to families of toroidal $3$-manifolds, including those with small order first homology and cyclic branched covers of toroidal links.  Their proofs can be found in \S \ref{sec: cbc and (*)}.

\label{sec: results} 
\subsection{Toroidal $3$-manifolds with small order first homology} 
\label{subsec: t with sofh}

Hanselman, Rasmussen and Watson have shown that if a closed, connected, orientable, irreducible $3$-manifold $W$ is a toroidal $L$-space, then either $|H_1(W)| = 5$ or $|H_1(W)| \geq 7$ (\cite[Theorem 57]{HRW1}). As they note (cf. \cite[Corollary 9]{HRW1}), this implies that prime $L$-spaces whose first homology has order $1, 2, 3, 4$, or $6$ are geometric. The proof is an application of the $NLS$-gluing theorem (\cite[Theorem 13]{HRW1}) combined with results on sets of $NLS$-detected  rational slopes on the boundaries of knot manifolds (\cite{RR}, \cite{HRW1}). Using the similar properties for $LO$ and $CTF$ detected slopes (see Theorem \ref{thm: * gluing} and Theorem \ref{thm: meridional detn}), we prove the following. 

\begin{theorem}
\label{thm: small order lo}
Suppose that $W$ is a closed, connected, orientable, irreducible, toroidal $3$-manifold such that $|H_1(W)| \in \{1, 2, 3, 4, 6\}$. Then,
\begin{enumerate}[leftmargin=*] 
\setlength\itemsep{0.5em}
\item[{\rm (1)}]  $\pi_1(W)$ is left-orderable if $|H_1(W)| \leq 4$; 
\item[{\rm (2)}]  $W$ admits a co-oriented taut foliation if either
\begin{itemize}
\item[\rm (a)] it contains an essential torus which splits it into two fibred knot manifolds, or 
\item[{\rm (b)}] it is an integer homology $3$-sphere which contains an essential torus bounding a fibred knot manifold to one side. 
\end{itemize}
\end{enumerate}
\end{theorem}
Since closed, connected, orientable, irreducible, atoroidal $3$-manifolds are geometric, we deduce: 

\begin{corollary}
\label{cor: small order nl is geometrico}
Suppose that $W$ is a closed, connected, orientable, irreducible $3$-manifold such that $|H_1(W)| \leq 4$. If $\pi_1(W)$ is non-left-orderable, then $W$ is a geometric manifold. 

\end{corollary}

\begin{remark}
\begin{enumerate}[leftmargin=*] 
\setlength\itemsep{0.5em}
\item[{\rm (1)}]  The requirement that $|H_1(W)| \ne 6$ in part (1) of Theorem \ref{thm: small order lo} and the fibredness condition in part (2) are necessitated by the fact that the current form of our $LO$- and $CTF$-detection results are weaker than their $NLS$ counterparts. For instance, in the $LO$ case, though our method works in most cases when $|H_1(W)| = 6$, it breaks down in the particular situation described in Remark \ref{rem: why not 6}, because our order-detection results require $H^2(W)$ to be a $\mathbb{Z}/2$-vector space. 
\item[{\rm (2)}] 
We note that the theorem is sharp in that there are precisely four irreducible, toroidal $L$-spaces with $|H_1(W)| = 5$ (\cite[Theorem 57]{HRW1}), each obtained by gluing trefoil exteriors, which are fibred. As $L$-spaces they admit no co-oriented taut foliations, and as graph manifolds with no co-oriented taut foliations they have non-left-orderable fundamental groups (\cite{BC1}). 
\end{enumerate}
\end{remark}

Ozsv\'ath and Szab\'o have conjectured that an irreducible integer homology $3$-sphere with infinite fundamental group is not an $L$-space, and this has been verified in the toroidal case by Eftekhary \cite{Eft} (see also \cite[Corollary 10]{HRW1}). Theorem \ref{thm: small order lo} shows that the following analogues hold, as predicted by the $L$-space conjecture. 

\begin{theorem}
\label{thm: tor zhs implies lo}
An irreducible, toroidal integer homology $3$-sphere $W$ has a left-orderable fundamental group. It admits a co-oriented taut foliation if it contains a fibred knot manifold whose boundary is incompressible in $W$. 

\end{theorem}

Theorem \ref{thm: tor zhs implies lo} shows that the fundamental group of a toroidal integer homology $3$-sphere $W$ containing an essential torus $T$ admits non-trivial representations with values in $\mbox{Homeo}_+(S^1)$, and it is natural to ask if $\mbox{Homeo}_+(S^1)$ can be replaced by $PSL(2, \mathbb R)$, especially as Lidman, Pinz\'on-Caicedo and Zentner have shown that it admits irreducible representations with values in $SU(2)$ (\cite{LPZ}). Gao has constructed hyperbolic integer homology $3$-spheres which admit no non-trivial $PSL(2, \mathbb R)$ representations (\cite{Gao}), later shown to have left-orderable fundamental groups (\cite[\S 8.8]{Du}). Her work can be used to produce an infinite family of toroidal integer homology $3$-spheres with the same property. 

\begin{proposition}
\label{prop: gao toroidal}
There is an infinite family of toroidal integer homology spheres whose fundamental groups admit no non-trivial $PSL(2, \mathbb R)$-representations. 
\end{proposition}

\subsection{The Gordon-Lidman conjecture} 
\label{subsec: gl conjecture and applications}

Gordon and Lidman initiated a study of the $L$-space conjecture for the cyclic branched covers $\Sigma_n(K)$ of knots $K$ in the $3$-sphere in \cite{GLid1}. When $K$ is the $(p, q)$-torus knot they showed that $\Sigma_n(K)$ is an $L$-space and has a non-left-orderable fundamental group precisely when $(n, |p|, |q|)$ is a Platonic triple; otherwise it has a left-orderable fundamental group and admits a co-oriented taut foliation. Partial results on the cyclic branched covers of families of hyperbolic knots have been obtained in a number of other papers, including  \cite{Hu1}, \cite{Tn}, \cite{Go}, \cite{HM}, \cite{BH}, \cite{BBG1}, \cite{BBG2}, \cite{HLL}, \cite{Tr}, and \cite{BGH}. Here we investigate cyclic branched covers 
of satellite knots, where the following is expected to occur.  

\begin{conjecture} {\rm (Gordon-Lidman)}
\label{conj: GL}
Each cyclic branched cover of $S^3$ over a prime satellite knot has a left-orderable fundamental group and admits a co-oriented taut foliation. 
\end{conjecture} 

Gordon and Lidman verified their conjecture in various cases, including cable knots (\cite[Theorem 1.3]{GLid1}, \cite[Theorem 1]{GLid2}). Theorem \ref{thm: results GL conjecture} below confirms the conjecture for $LO$ and $NLS$, and shows that it also holds for $CTF$ when for instance the companion knot is fibred. In fact the theorem proves the corresponding extension to satellite links, where we define a satellite link as follows. Let $K$ be a non-trivial knot in $S^3$, so $S^3 = V \cup X_K$ where $V$ is a solid torus and $X_K$ is the exterior of $K$. Let $P$ be a link in $V$ such that $V \setminus P$ is irreducible and $P$ is not a core of $V$. Then the image of $P$ in $S^3$ is a {\it satellite} link which we denote by $P(K)$. The link $P$ in $V$ is the {\it pattern} of $P(K)$, and $K$ is the {\it companion} of $P(K)$. A link $L$ is {\it toroidal} if its exterior $X_L$ contains an essential torus. A satellite link $L = P(K)$ is toroidal: the torus $\partial V = \partial X_K$ is essential in $X_L$. Clearly the converse is true if $L$ is a knot, but not in general.

\begin{theorem}
\label{thm: results GL conjecture}
 Suppose that $P(K)$ is a prime satellite  link with pattern $P \subset S^1 \times D^2$ and companion $K$. Let $\Sigma_n(P(K))$ denote the $n$-fold cyclic cover of $S^3$ branched over $P(K)$. Then
 \begin{enumerate}[leftmargin=*] 
\setlength\itemsep{0.5em}
\item[{\rm (1)}]  $\Sigma_n(P(K))$ is not an $L$-space for each $n\geq 2$; 
\item[{\rm (2)}]  $\pi_1(\Sigma_n(P(K)))$ is left-orderable for each $n\geq 2$; 
\item[{\rm (3)}] $\Sigma_n(P(K))$ admits a co-orientable taut foliation for each $n\geq 2$ if either of the following conditions holds:
\begin{enumerate}[leftmargin=*] 
\setlength\itemsep{0.5em}
\item[{\rm (a)}]   $K$ is fibred; 
\item[{\rm (b)}]  $P$ is braided and either $K$ is persistently foliar or $n$ has a prime factor which does not divide the winding number of $P$ in $S^1 \times D^2$.
\end{enumerate}
\end{enumerate}
\end{theorem} 
We refer the reader to \S \ref{subsec: satellite 2} for the definition of persistently foliar.  As a consequence of Theorem \ref{thm: results GL conjecture} we show

\begin{theorem}
\label{thm: GL thm intro}
 If   $L$ is a prime satellite  link then for all $n \ge 2$, $\Sigma_n(L)$ is not an $L$-space and has a left-orderable fundamental group.  If in addition   $L$ has a fibred companion knot, then $\Sigma_n(L)$ admits a co-oriented taut foliation for all $n \geq 2$.
\end{theorem}

\begin{corollary}
Each cyclic branched cover of a prime satellite  link with a fibred companion has a left-orderable fundamental group and admits a co-oriented taut foliation. In particular this holds for 
prime fibred satellite links. 
\end{corollary}

There is an extension of Theorem \ref {thm: results GL conjecture} to links in homology 3-spheres. To state this, define a {\it companion manifold} of a link $L$ with exterior $X$ in a closed $3$-manifold to be a knot manifold $M \subset \mbox{int}(X)$ for which $\partial M$ is essential in $X$.

\begin{theorem} 
\label{thm: detn and bccs} 
Suppose that   $L$ is a prime  link in an integer homology $3$-sphere $W$ whose exterior is irreducible and contains a companion manifold. If $\Sigma_n(L)$ denotes the $n$-fold cyclic cover of $W$ branched over   $L$, then 
\begin{enumerate}[leftmargin=*] 
\setlength\itemsep{0.5em}
\item[{\rm (1)}]   $\Sigma_{n}(L)$ is not an $L$-space for each $n \geq 2$;  
\item[{\rm (2)}]  $\Sigma_n(L)$ has a left-orderable fundamental group for each $n \geq 2$;   
\item[{\rm (3)}]  $\Sigma_n(L)$ admits a co-oriented taut foliation for each $n \geq 2$ if the exterior of $L$ contains a fibred companion manifold.
\end{enumerate}
\end{theorem}

Theorem \ref{thm: detn and bccs} yields Theorem \ref{thm: results GL conjecture} except for the existence of a co-orientable taut foliation on $\Sigma_n(P(K))$ in the situations described in part (3)(b) of its statement, which is dealt with in \S \ref{subsec: satellite 2}. 

We expect that Conjecture \ref{conj: GL} holds for oriented toroidal links.  
\begin{conjecture} 
    \label{conj: n-fold toroidal link}
Each cyclic branched cover of $S^3$ over a prime oriented toroidal link has a left-orderable fundamental group and admits a co-oriented taut foliation. 
\end{conjecture} 

  Our next result shows that the left-orderable part of Conjecture \ref{conj: n-fold toroidal link} holds generically.

\begin{theorem}
\label{thm: theorem y}
Suppose that $L$ is a prime link in an irreducible integer homology $3$-sphere $W$ whose exterior is irreducible and contains an essential torus $T$ which is not a Heegaard surface of $W$. Then for any orientation on $L$, $\Sigma_n(L)$ has a  left-orderable fundamental group for any $n\geq 2$.  
\end{theorem}

Theorem \ref{thm: theorem y} shows that the only case of the left-orderable part of Conjecture \ref{conj: n-fold toroidal link} which remains to be verified is when $W = S^3$ and every essential torus in $X_L$ splits $S^3$ into two solid tori $X_1$ and $X_2$ such that $L\cap X_i = L_i$ is non-empty for $i=1, 2$. 
Note that though the $n$-fold cyclic branched cover of a link $L$ depends on the orientation of the link when $n>2$, the property of being toroidal does not.

For the $NLS$ part of Conjecture \ref{conj: n-fold toroidal link} we have the following partial result.
      
\begin{theorem}
\label{thm: toroidal links nls}
If $L$ is a prime toroidal link in an integer homology $3$-sphere whose exterior is irreducible, then  $\Sigma_2(L)$ is not an $L$-space. More generally, for any orientation on $L$, 
$\Sigma_n(L)$ is not an $L$-space for $n = 2^k$, $k \geq 1$, and $n = 3\cdot2^k$, $k \geq 0$. 
\end{theorem}

Theorem \ref{thm: toroidal links nls} has an interesting consequence for links in the $3$-sphere. 
Ozsv\'ath and Szab\'o have shown that the $2$-fold cyclic branched cover of a $\mathbb{Z}/2$-Khovanov thin link $L$ is an $L$-space. More precisely, 
they established the inequality 
$$|\det(L)| \leq \mbox{dim}_{\mathbb Z/2} \widehat{HF}(\Sigma_2(L); \mathbb Z/2) \leq \mbox{dim}_{\mathbb Z/2} \widetilde{\mbox{Kh}}(L; \mathbb Z/2)$$
 (\cite[Corollary 1.2]{OS2}), where $|\det(L)|$ is the order of $H_1(\Sigma_2(L))$ and $\widetilde{\mbox{Kh}}(L; \mathbb Z/2)$ denotes the reduced Khovanov homology of $L$ with $\mathbb Z/2$ coefficients. Since $\mbox{dim}_{\mathbb Z/2} \widetilde{\mbox{Kh}}(L; \mathbb Z/2) = |\det(L)| = |H_1(\Sigma_2(L))|$ for a $\mathbb Z/2$-Khovanov thin link, it follows that $\Sigma_2(L)$ is an $L$-space. Thus, as a corollary of Theorem \ref{thm: toroidal links nls} we deduce,

\begin{corollary}
\label{cor: satellite not H-thin}
Prime toroidal links in the $3$-sphere are not $\mathbb Z/2$-Khovanov thin. 
\end{corollary}

Menasco has shown that prime toroidal links are never alternating (\cite[Corollary 2]{Me}), and it is a folklore conjecture that they are never quasi-alternating. Since quasi-alternating links are $\mathbb Z/2$-Khovanov thin (\cite{MO}), Corollary \ref{cor: satellite not H-thin} allows us to verify this. 

\begin{theorem}
\label{thm: satellite not qa} 
A prime $\mathbb Z/2$-Khovanov thin link is either hyperbolic or a $(2, m)$-torus link. In particular, this holds for prime quasi-alternating links.  
\end{theorem}

\section{Slopes and multislopes}
\label{sec: conventions}

From now on we assume that manifolds are connected and orientable, unless otherwise stated. 

Consider a compact $3$-manifold $M$ with torus boundary. Duality implies that the kernel of the homomorphism $H_1(\partial M; \mathbb Q) \to H_1(M; \mathbb Q)$ has dimension $1$. Equivalently, if $T_1(M)$ denotes the torsion subgroup of $H_1(M)$, there is a primitive element $\lambda_M \in H_1(\partial M)$, unique up to sign, whose image in $H_1(M)$ is contained in $T_1(M)$. This class is called the (rational) {\it longitude} of $M$. We will refer to any class $\mu$ which forms a basis of $H_1(\partial M)$ with $\lambda_M$ as a {\it meridional class}.   

A {\it slope} on a torus $T$ is an element of the projective space of $H_1(T; \mathbb R)$ and the slope of a non-zero class $\gamma \in H_1(T; \mathbb R)$ will be denoted by $[\gamma]$. Rational slopes, which are those corresponding to lines defined by non-zero elements of $H_1(T) \subset H_1(T; \mathbb R)$, can be identified as either $\pm$-pairs of primitive elements of $H_1(T)$ or isotopy classes of essential simple closed curves on $T$.  The {\it distance} between two rational slopes on $T$ is the absolute value of the algebraic intersection number between primitive elements of $H_1(T)$ representing the slopes. Thus the distance is $0$ if and only if the slopes coincide and is $1$ if and only primitive representatives of the slopes form a basis of $H_1(T)$. 

If $M$ is a compact $3$-manifold with torus boundary, we use $\mathcal{S}(M)$ to denote the set of slopes on its boundary and $\mathcal{S}_{rat}(M)$ to denote the set of rational slopes. 

The {\it longitudinal slope} of $M$ is the  rational slope represented by a longitudinal class $\lambda_M \in H_1(\partial M)$ and a {\it meridional slope} is a  rational slope represented by a meridional class. From time to time we will somewhat ambiguously denote the slope of a primitive class $\alpha \in H_1(\partial M)$ as simply ``$\alpha$". In particular, this is the convention when we are performing Dehn filling on $M$ along $[\alpha]$, which is denoted by $M(\alpha)$.

More generally, if $M$ is a compact $3$-manifold whose boundary is a union of tori $T_1, T_2, \ldots , T_r$, the set  of {\it multislopes} on $\partial M$ is given by 
$$\mathcal{S}(M) = \mathcal{S}(T_1) \times \mathcal{S}(T_2) \times \cdots \times \mathcal{S}(T_r) \cong (S^1)^r$$ 
and  the set  of rational multislopes by
$$\mathcal{S}_{rat}(M) = \mathcal{S}_{rat}(T_1) \times \mathcal{S}_{rat}(T_2) \times \cdots \times \mathcal{S}_{rat}(T_r)$$

\section{NLS-detection and gluing knot manifolds} 
\label{sec: $NLS$ detn}

We work with Heegaard Floer homology with $\mathbb Z/2$-coefficients throughout the paper. In particular we take ``$L$-space" to mean ``$\mathbb Z/2$-$L$-space". 

\subsection{$NLS$-detection}
\label{subsec: NLS detn}
Consider a rational homology solid torus $M$ and the set
$$\mathcal{L}(M) = \{ \mbox{rational slopes } \alpha \; | \; M(\alpha) \mbox{ is an $L$-space}\},$$ 
It is shown in the paragraph following the proof of Corollary 56 of \cite{HRW1} that $\mathcal{L}(M)$
is either empty, a closed subinterval of $\mathcal{S}_{rat}(M)$, or $\mathcal{S}_{rat}(M) \setminus \{\lambda_M\}$. (This had previously been proven for the set of $\mathbb Z$-$L$-space filling slopes in \cite[Theorem 1.6]{RR}.) 

In the case that $\mathcal{L}(M)$ contains at least two slopes, the exterior of an $L$-space knot for instance, $M$ is called {\it Floer simple}. More precisely, 
$M$ is $\mathbb Z/2$-Floer simple, and though ostensibly this could differ from the condition that $M$ be $\mathbb Z$-Floer simple, they turn out to be equivalent conditions. 
One direction of the verification is immediate; the universal coefficient theorem implies that $\mathbb Z$-$L$-spaces are $\mathbb Z/2$-$L$-spaces and therefore $M$ is $\mathbb Z/2$-Floer simple if it is $\mathbb Z$-Floer simple. Conversely, as noted in the remarks on coefficients on page 610 of \cite{HRRW}, it follows from the proof of \cite[Proposition 3.7]{RR} that a Dehn filling of a $\mathbb Z/2$-Floer simple manifold is a $\mathbb Z/2$-$L$-space if and only if it is a $\mathbb Z$-$L$-space. Thus $M$ is $\mathbb Z/2$-Floer simple if and only if it is $\mathbb Z$-Floer simple.

When $\mathcal{L}(M) = \mathcal{S}_{rat}(M) \setminus \{\lambda_M\}$, $M$ is called an {\it HF-generalised solid torus}\footnote{Our definition of $HF$-generalised solid tori differs from that found in \cite[Definition 7.2]{RR} or \cite[Definition 2.2]{Gi}, but is equivalent to it by \cite[Proposition 7.1]{RR} or \cite[Theorem 2.3 and Corollary 2.7]{Gi}.}. An example is given by the twisted $I$-bundle over the Klein bottle, which we denote by $N$, since its non-longitudinal fillings either have finite fundamental groups or are $P^3 \# P^3$. A restricted version of the following result was proven in \cite {RR}. The full version is due to Gillespie \cite{Gi}. 

\begin{proposition} {\rm (Rasmussen-Rasmussen, Gillespie)}
\label{prop: hfst}
If $M$ is an HF-generalised solid torus, then any Thurston norm minimizing 
surface $S$ in $M$ representing a generator $\eta$ of $H_2(M, \partial M)$ has genus zero. Further, $|\partial S|$ is the order of the longitude  of $M$ in $H_1(M)$. Thus $M$ is boundary-compressible if the longitude  of $M$ is null-homologous in $M$. 
\end{proposition}

\begin{proof} 
Let $d$ be the order of the rational longitude in $H_1(M)$. 
It follows from \cite[Theorem 2.3]{Gi} and \cite[Theorem 1.6, Corollary 2.3]{RR} that 
$$2g(S) + |\partial S| - 2 \leq \max \{0, 2g(S) + |\partial S| - 2\} = \|\eta\| \leq  d - 2,$$
so as $d \leq |\partial S|$ we have
$$0 \leq 2g(S) \leq d -  |\partial S| \leq 0$$
Thus $g(S) = 0$ and $|\partial S| = d$.
\end{proof}

The complement of $\mathcal{L}(M)$ in $\mathcal{S}_{rat}(M)$ is the set of 
{\it strongly NLS-detected rational slopes}, whose closure in $\mathcal{S}(M)$ 
$$\mathcal{D}_{NLS}(M) = \overline{\mathcal{S}_{rat}(M) \setminus \mathcal{L}(M)}$$ 
is called the set of {\it $NLS$-detected slopes}.

We record the following simple consequence of the definitions for later use. 

\begin{proposition} 
\label{prop: hfgst implies small detd set}
$[\lambda_M] \in \mathcal{D}_{NLS}(M)$ for each rational homology solid torus $M$ and if $M$ is an $HF$-generalised solid torus, then $\mathcal{D}_{NLS}(M) = \{[\lambda_M]\}$. 
\end{proposition}
 
It is well-known that if $M$ is the exterior of a non-trivial knot in the $3$-sphere, then there is a sequence of slopes in $\mathcal{S}_{rat}(M) \setminus \mathcal{L}(M)$ which converges to the meridional slope (cf. \cite[Theorem 1]{BS}), which is therefore $NLS$-detected. More generally, we have the following result, which includes Theorem \ref{thm: meridional detn}(1).

\begin{proposition}
\label{prop: mu + k lambda NLS-detd}
Let $M \not \cong S^1 \times D^2$ be an irreducible rational homology solid torus whose longitude  $\lambda_M$ is integrally null-homologous.  

\begin{enumerate}[leftmargin=*] 
\setlength\itemsep{0.5em}
\item[{\rm (1)}] Any rational slope of distance $1$ from $\lambda_M$ is contained in $\mathcal{D}_{NLS}(M)$. Moreover, all but at most two such slopes are contained in $\mathcal{S}_{rat}(M) \setminus \mathcal{L}(M)$, and if two, they are the endpoints of $\mathcal{L}(M)$ and are of distance $1$ from each other. 
\item[{\rm (2)}] If $M$ is an integer homology solid torus and $d \geq 1$ is a divisor of $2g(M) - 1$, then each rational slope of distance $d$ from $\lambda_M$ is $NLS$-detected.
\end{enumerate}
\end{proposition}

\begin{proof} 
Let $\lambda \in H_1(\partial M)$ represent the longitudinal slope of $M$ and $\mu  \in H_1(\partial M)$ represent a meridional slope. Then the  rational slopes of distance $1$ from $\lambda$ are 
represented by the classes $\{\mu + k \lambda \; | \; k \in \mathbb Z\}$.

Proposition \ref{prop: hfst} shows that $M$ is not an $HF$-generalised solid torus. If it isn't Floer simple, the set $\mathcal{L}(M)$ of 
$L$-space filling slopes of $M$ consists of at most one element by definition \cite{RR} and therefore $\mathcal{D}_{NLS}(M) = \overline{\mathcal{S}_{rat}(M) \setminus \mathcal{L}(M)} = \mathcal{S}(M)$, which 
implies the desired conclusions. 

Next suppose that $M$ is Floer simple, though not a HF-generalised solid torus. Since $\lambda$ is null-homologous in $H_1(M)$, there is a canonical splitting $H_1(M) = \langle \bar \mu\rangle \oplus T$ where $T$ is the torsion subgroup of  $H_1(M)$ and $\langle \bar \mu\rangle \cong \mathbb Z$ is generated by the image $\bar \mu$ of $\mu$ under the inclusion-induced homomorphism $\iota: H_1(\partial M) 
\to H_1(M)$. Further, the image of 
$\iota$ is precisely $\langle \bar \mu\rangle$. Let $\phi: H_1(M)  \to \mathbb Z$ be the homomorphism which sends $\bar \mu$ to $1$. 

The Turaev torsion of $M$ is an element $\tau(M)$ of the ring of possibly infinite formal sums $\sum_{\alpha \in H_1(M)} a_\alpha [\alpha]$ where $a_\alpha \in \mathbb Z$ for all $\alpha$ and is zero for $\phi(\alpha) \ll 0$. It turns out that  
$$\tau(M) = \sum_{\phi(\alpha) \geq 0} a_\alpha [\alpha],$$
where $a_0 \ne 0$ and $a_\alpha$ is $1$ for all but finitely many $\alpha$ with $\phi(\alpha) \geq 0$. The {\it support} of $\tau(M)$ is the subset $S[\tau(M)] = \{\alpha \in H_1(M) \; | \; a_\alpha \ne 0\}$ of $H_1(M)$. 

Set 
$$\mathcal{D}_+^\tau(M) = \{\beta - \alpha \; | \; \beta \not \in S[\tau(M)], \alpha \in S[\tau(M)], \phi(\beta) > \phi(\alpha)\} \cap \langle \bar \mu\rangle \subset H_1(M)$$
Theorem 1.6 of \cite{RR} shows that the endpoints of $\mathcal{D}_{NLS}(M)$ are contained in $\iota^{-1}(\mathcal{D}^\tau_+ (M))$. (See also Theorem 68 of \cite{HRW2}.) Further, the rays determined by these endpoints cobound a sector $\mathfrak{S}(M)$ in $H_1(\partial M; \mathbb R)$ which contains $\mathcal{L}(M)$ and whose interior is disjoint from $\iota^{-1}(\mathcal{D}^\tau_+ (M))$. 

Our hypotheses imply that $\mathcal{D}^\tau_+(M)$ is non-empty, so there are $\alpha, \beta \in H_1(M)$ such that $\alpha \in S[\tau(M)],  \beta \not \in S[\tau(M)]$, and $\beta = \alpha + r \bar \mu\;$ for some $r > 0$. Fix an integer $0 \leq s < r$ such that $\alpha + s \bar \mu \in S[\tau(M)]$ and $\alpha + (s+1) \bar \mu \not \in S[\tau(M)]$, and observe that $\bar \mu \in \mathcal{D}^\tau_+ (M)$ since $(\alpha + (s+1) \bar \mu) 
- (\alpha + s \bar \mu) = \bar \mu$ and $\phi(\bar \mu) = 1$. Hence
$$\{ \mu + k \lambda \; | \; k \in \mathbb Z\} = \iota^{-1}(\bar \mu) \subset \iota^{-1}(\mathcal{D}^\tau_+ (M))$$
Part (1) of the proposition holds if no $M(\mu + k \lambda)$ is an $L$-space. On the other hand, if some $M(\mu + k \lambda)$ is an $L$-space, $\mu  + k \lambda$ is contained in the sector $\mathfrak{S}(M)$, and as the interior of this sector is disjoint from $\iota^{-1}(\mathcal{D}^\tau_+ (M))$, the only possibility is for $\mu  + k \lambda$ to lie on one of the sector's bounding rays. Then $\mu + k \lambda$ must be an endpoint of the closed interval $\mathcal{D}_{NLS}(M)$ and so is $NLS$-detected. 
Thus there are at most two 
such slopes. If there are two, they correspond to the bounding rays of $\mathfrak{S}(M)$, and as $\mbox{int}(\mathfrak{S}(M))\cap  \{\mu + k \lambda \; | \; k \in \mathbb Z\} \subseteq \mbox{int}(\mathfrak{S}(M)) \cap \iota^{-1}(\mathcal{D}^\tau_+ (M)) = \emptyset$, there is an $\epsilon \in \{\pm 1\}$ such that the two slopes correspond to $\mu + k \lambda$ and $\mu + (k + \epsilon) \lambda$. Hence the distance between them is $1$. 
This completes the proof of (1). 

To prove (2), suppose that $M$ is an integer homology solid torus of genus $g$ and Alexander polynomial $\Delta_M(t)$. As above we can assume that 
$M$ is Floer simple, though not an HF-generalised solid torus. Then by \cite[Corollary 2.3]{RR}, $\mbox{deg} \; \Delta_M(t) = \|M\| + 1 = 2g$, where $\| M \|$ denotes the Thurston norm of a generator of $H_2(M, \partial M)$.

Under our assumptions, $\mathbb Z \cong H_1(M) = \langle \bar \mu \rangle$, and writing $H_1(M)$ multiplicatively with generator $t$ we can express $\tau(M)$ as a formal sum $\sum_{r \geq 0} a_r t^r$ where $a_r \in \{0, 1\}$ for each $r$ and $a_0 = 1$. Further, $a_r = 1$ for $r \geq 2g(M)$ (cf. \cite[Proposition 2.2]{RR}). Then by \cite[\S 5.2]{Tv} we have, 

\begin{eqnarray}
\Delta_M(t) = (1-t)\tau(M) & = & (1-t)(1 + \sum_{r = 1}^{2g-1} a_r t^r) + (1-t)(\sum_{r \geq 2g} t^r) \nonumber \\ 
& = & 1+ \sum_{r = 1}^{2g-1}(a_r- a_{r-1}) t^r  \,+\, (1 - a_{2g-1}) t^{2g} \nonumber 
\end{eqnarray} 
Since $\Delta_M(t^{-1}) \sim \Delta_M(t)$, we have $1 = a_0 = a_{2g} = 1 - a_{2g-1}$ and therefore $a_{2g-1} = 0$. Hence $(2g-1)\bar \mu \not \in S[\tau(M)]$ while $0  \in S[\tau(M)]$. Now we proceed as in the proof of the first assertion of (1): Consideration of the sequence $0, d \bar \mu, 2d \bar \mu, \ldots , (r-1)d \bar \mu, rd \bar \mu = (2g-1)\bar \mu$ leads to the conclusion that $d \bar \mu \in \mathcal{D}^\tau_+ (M)$ and therefore $\{ d\mu + k \lambda \; | \; k \in \mathbb Z\} = \iota^{-1}(\bar \mu) \subset \iota^{-1}(\mathcal{D}^\tau_+ (M))$. As in (1), this implies that the slope associated to each primitive class of the form 
$d\mu + k \lambda$ is $NLS$-detected. 
\end{proof}

It follows from part (1) of Proposition \ref{prop: mu + k lambda NLS-detd} and the results of \cite{RR} and \cite{HRW1} that there is a basis
$\{ \mu, \lambda_M\}$ of $H_1(\partial M)$ such that each  rational slope in at least one of the intervals $[-\infty, 1]$ or $[-1, \infty]$ is $NLS$-detected. We expect that analogous statements hold for the sets of $CTF$-detected  rational slopes and $LO$-detected  rational slopes.

\subsection{The $NLS$-gluing theorem for knot manifolds}
The following gluing theorem, which is Theorem 13 of \cite{HRW1}, contains the $\ast = NLS$ case of Theorem \ref{thm: * gluing}. 

\begin{theorem}
{\rm (Hanselman-Rasmussen-Watson)} 
\label{thm: HRW1}
If $M_1, M_2$ are two knot manifolds and $f: \partial M_1 \to \partial M_2$ is a homeomorphism, then $W = M_1 \cup_f M_2$ is not an $L$-space if and only if $f$ identifies a rational slope in $\mathcal{D}_{NLS}(M_1)$ with one in $\mathcal{D}_{NLS}(M_2)$. 
\end{theorem}
The rationality requirement in the theorem can be removed; the restrictions on $\mathcal{L}(M_i)$ described in the first paragraph of \S \ref{subsec: NLS detn} imply that if $f: \partial M_1 \to \partial M_2$ identifies a pair of irrational slopes, it also identifies a pair of rational slopes. 

Theorem \ref{thm: HRW1} leads to the alternate characterisation of $NLS$-detection of 
rational slopes stated in the introduction. 

\begin{corollary} 
{\rm (Hanselman-Rasmussen-Watson)} 
\label{cor: gluing defn of nls detection}
Suppose that $M$ is a knot manifold and $W = M \cup_f U$ is a manifold obtained by gluing an irreducible $HF$-generalised solid torus $U \not \cong S^1 \times D^2$ to $M$ using a homeomorphism $f: \partial U \to \partial M$ which identifies the longitudinal slope $[\lambda_U]$ of $U$ with a rational slope $[\alpha] \in \mathcal{S}_{rat}(M)$. Then $[\alpha]$ is $NLS$-detected in $M$ if and only if $W$ is not an $L$-space. 
\end{corollary}

\subsection{$NLS$-detection of multislopes} 
Suppose that $M$ is a compact, irreducible $3$-manifold whose boundary is a union $\partial M = T_1 \sqcup T_2 \sqcup \cdots \sqcup T_r$ of incompressible tori and $[\alpha] = ([\alpha_1], \ldots , [\alpha_r]) \in \mathcal{S}_{rat}(M)$ is a  rational multislope. Let $U_1, \ldots, U_r$ be irreducible boundary-incompressible HF-generalised solid tori. Set 
$$W = M \cup_{T_1}U_1 \cup_{T_2} \cdots \cup_{T_r} U_r,$$ 
where $U_i$ is glued to $M$ along $T_i$ so that $[\lambda_{U_i}] = [\alpha_i]$. Note that $W$ is not uniquely determined by $[\alpha]$, since the gluing condition does not determine the isotopy class of the attaching map. A simple induction using Corollary \ref{cor: gluing defn of nls detection} shows that the condition that $W$ not be an $L$-space is independent of the choice of gluing maps and the $U_i$. 

\begin{definition}
\label{def: nls-detd multislope}
We say that the rational multislope $[\alpha]$ is {\it $NLS$-detected} if $W$, as above, is not an $L$-space.   
\end{definition}

\section{Foliation-detection and gluing knot manifolds}
\label{sec: foln gluing knot mflds}

\subsection{Foliations on tori and slopes}  
\label{subsec: slopes of folns}
In this section, we define what it means for a slope to be detected by a co-oriented taut foliation. Before giving the general definition (Definition \ref{def: ctf detection}), we first point out that a codimension $1$ foliation $\mathcal{F}$ on a torus $T$ without Reeb annuli detects a rational slope $[\alpha]$ if and only if it contains a closed leaf of slope $[\alpha]$. The slope of this leaf is well-defined since any two closed leaves are disjoint essential simple closed curves. 

In general, the absence of Reeb annuli in $\mathcal{F}$ is equivalent to it being a suspension   foliation. In other words, $(T, \mathcal{F})$ is pairwise homeomorphic to $(T_h, \mathcal{F}_h)$ where $T_h$ is the mapping torus of a homeomorphism $h \in \mbox{Homeo}_+(S^1)$, called the {\it holonomy} of $\mathcal{F}$, and $\mathcal{F}_h$ is the associated suspension foliation on $T_h$ (\cite[Proposition 4.3.2]{HH}). Given $h_1, h_2 \in \mbox{Homeo}_+(S^1)$, there exists a fibre preserving homeomorphism $f: (T_{h_1}, \mathcal{F}_{h_1}) \rightarrow (T_{h_2}, \mathcal{F}_{h_2})$ if and only if $h_1$ is conjugate to $h_2$ (\cite[Theorem 3.15]{CC1}). 

We call $\mathcal{F}$ {\it linear} if its holonomy is conjugate to a rotation of $S^1$. In this case $\mathcal{F}$ is a fibration by simple closed curves if the rotation angle is a rational multiple of $\pi$. Otherwise, it is a foliation by real lines each of which is dense in $T$. The converse holds;   a fibration on $T$ by simple closed curves  is linear, as is a foliation on $T$ by real lines, each of which is dense.

Let $\widetilde{\mathcal{F}}$ be the pullback foliation of a suspension foliation $\mathcal{F}$ on $T$ to the universal cover $\widetilde{T}\cong \mathbb{R}^2$. A characterising feature of suspension foliations  is the existence of a fibration of $T$ by simple closed curves everywhere transverse to it. It is easy to see that each component of the inverse image in $\widetilde T$ of any fibre of the fibration is a line transversely intersecting every leaf of $\widetilde{\mathcal{F}}$ once and only once,  from which it follows that the leaf space $\mathcal{L}$ of $\widetilde{\mathcal{F}}$ is a real line.  There is a natural action of $\pi_1(T)$ on $\mathcal{L}$ induced by the covering action, which determines the slope $[\alpha] \in \mathcal{S}(T)$ detected by $\mathcal{F}$ as follows. 

The action of $\pi_1(T)$ on $\mathcal{L}$ has no global fixed points and so as $\pi_1(T)$ is abelian, the stabiliser of a point $x \in \mathcal{L}$ is either trivial or an infinite cyclic summand of $\pi_1(T)$ generated by a primitive element $\alpha \in H_1(T) \equiv \pi_1(T)$. In the latter case, the leaf of $\widetilde{\mathcal{F}}$ corresponding to $x$ projects to a closed leaf of $\mathcal{F}$ of slope $[\alpha]$, which is therefore the slope of $\mathcal{F}$. The reader will verify that the elements of $\pi_1(T) = H_1(T) \subset H_1(T,\mathbb{R})$ which lie to one side of the line on $H_1(T,\mathbb{R})$ determined by the slope $[\alpha]$ act on $\mathcal{L}$ in a strictly increasing fashion and those which lie to the other side act in a strictly decreasing fashion. 

If $\pi_1(T)$ acts freely on $\mathcal{L} \equiv \mathbb R$, then each non-trivial element $\gamma$ of $\pi_1(T)$ acts in either a strictly increasing or strictly decreasing way, and it is not hard to see that there is a line $[\alpha]$ in $H_1(T; \mathbb R)$ such that if $\gamma$ lies to one side of $[\alpha]$, it acts in a strictly increasing fashion, and if it lies to the other, it acts in a strictly decreasing fashion. Here we define the slope of $\mathcal{F}$ to be $[\alpha]$ and note that it is irrational since $H_1(T) \cap [\alpha] = \{0\}$.

\begin{definition}
\label{def: ctf detection}
A slope $[\alpha]$ on the boundary of a knot manifold $M$ is {\it $CTF$-detected} if there is a co-oriented taut foliation $\mathcal{F}$ on $M$ which intersects $\partial M$ transversely in a suspension foliation of slope $[\alpha]$. We say that $\mathcal{F}$  {\it strongly $CTF$-detects} $[\alpha]$ if $\mathcal{F} \cap \partial M$ is a linear foliation of slope $[\alpha]$. 
\end{definition} 

We use $\mathcal{D}_{CTF}(M)$ to denote the set of $CTF$-detected slopes on $\partial M$. It follows from \cite{Gab1} that $\mathcal{D}_{CTF}(M)$ always contains $[\lambda_M]$.

\subsection{Gluing foliated knot manifolds}
Here we prove the foliation case of Theorem \ref{thm: * gluing}.

\begin{theorem}
\label{thm: fln gluing}
Suppose that $M_1$ and $M_2$ are knot manifolds and $W = M_1 \cup_f M_2$ where $f: \partial M_1 \rightarrow \partial M_2$ is a homeomorphism which identifies rational slopes $[\alpha_1] \in \mathcal{D}_{CTF}(M_1)$ and $[\alpha_2] \in \mathcal{D}_{CTF}(M_2)$. Then $W$ admits a co-orientable taut foliation. 
\end{theorem}

\begin{proof} 
We are done when the first  Betti number  of $W$ is $1$ or more by \cite{Gab1}, so without loss of generality we assume that $M_1$ and $M_2$ are rational homology solid tori and $[\alpha_1] \neq [\lambda_1]$.  

For each $i=1,2$, let $\mathcal{F}_i$ be a taut foliation on $M_i$ which detects the slope $[\alpha_i]$ and let $b_1$ be a simple closed curve on $\partial M_1$ everywhere transverse  to $\mathcal{F}_1$. We can assume that the  rational slope $[b_1]$ of $b_1$ has distance $1$ from $[\alpha_1]$. Isotope $f$ so that $b_2 = f(b_1)$ is everywhere transverse  to $\mathcal{F}_2$.   

Our strategy is to alter each $\mathcal{F}_i$ to produce a new foliation $\mathcal{F}'_i$ which detects $\alpha_i$, is transverse to $b_i$ and if $h_i': b_i\rightarrow b_i$ denotes the holonomy of $\mathcal{F}'_i|_{\partial M_i}$, then $f|_{b_1}\circ h'_1 \circ (f|_{b_1})^{-1}$ 
is conjugate to $h'_2$ in $\mbox{Homeo}_+(b_2)$. Since  $f|_{b_1}\circ h'_1 \circ (f|_{b_1})^{-1}$ is the holonomy map of $f(\mathcal{F}_1'|_{\partial M_1})$, there exists an orientation-preserving  homeomorphism $f': (\partial M_2, f(\mathcal{F}_1'|_{\partial M_1})) \rightarrow (\partial M_2, \mathcal{F}_2'|_{\partial M_2})$ which preserves $b_2$ as an oriented curve. On the other hand, our hypotheses imply that the slope of $ f(\mathcal{F}_1'|_{\partial M_1})$ is $[\alpha_2]$, so $f'_*: H_1(\partial M_2) \to H_1(\partial M_2)$ preserves $\alpha_2 \in H_1(\partial M_2)$. Since $[\alpha_2] \ne [b_2]$, it follows that $f'_*$ is the identity on $H_1(\partial M_2)$. Hence $f'$ is isotopic to the identity map, which shows that $f(\mathcal{F}_1'|_{\partial M_1})$ and $\mathcal{F}_2'|_{\partial M_2}$ are isotopic foliations. Consequently, we can obtain a taut foliation on $W$ by gluing the $\mathcal{F}_i'$ together. Here are the details. 

Fix a leaf $L_i$ of $\mathcal{F}_i$, $i = 1, 2$, whose intersection with $\partial M_i$ is non-empty and contains a compact component $a_i$, necessarily a simple closed curve of slope $[\alpha_i]$. By assumption, $b_i$ intersects $a_i$ once and hence their point of intersection $p_i$ is fixed under the holonomy map. 

For each $i=1,2$ we perform a Denjoy blow up on the leaf $L_i$. Topologically this is equivalent to replacing $L_i$ with the product $L_i\times I$ foliated by $\{L_i\times t: t\in I\}$ and injectively immersed in $M_i$ (see \cite[Operation 2.1.1]{Gab2}). Let $\iota_i: L_i\times I\rightarrow M_i$ denote the immersion. The new foliations, which we continue to call $\mathcal{F}_1, \mathcal{F}_2$, detect the same slopes. 

Let $b_i'$ be the intersection $b_i\cap \iota_i(a_i\times I)$ and $b''_i = \overline{b_i\setminus b_i'}$ (see Figure \ref{fig: holonomy case II}) and note that $b_i'$ and $b_i''$ are invariant under the holonomy maps.  For $i=1,2$, let $h''_i: b''_i\rightarrow b''_i$ be the holonomy map of $\mathcal{F}_i|_{\partial M_i}$ over $b''_i$.  Note that $h''_i$ is determined by the restriction of the holonomy of the pre-blow up version $\mathcal{F}_i|_{\partial M_i }$ to the complement of $p_i$. The holonomy map of $\mathcal{F}_i|_{\partial M_i}$ over $b_i'$ is trivial. 

Suppose that $L_1$ and $L_2$ are noncompact and for each $i$ choose an infinite ray $\gamma_i$ properly embedded in $L_i$ starting at a point on $a_i\subset \partial L_i$ and terminating in an end of $L_i$. We alter the 
holonomy map over $b'_1$  by cutting the product foliation $\mathcal{F}_1|_{\iota_1(L_1\times I)}$ open along $\iota_1(\gamma_1 \times I)$ and regluing through the map given by ${\rm id} \times h''_2: \gamma_1\times I\rightarrow \gamma_1\times I$. \footnote{Here we view $h''_2: b_2''\rightarrow b_2''$ as a homeomorphism of the interval $I$ The specific choice of the identification between $b''_2$ and $I$ is insignificant, since different identifications give us the same homeomorphism of $I$ up to conjugation.} We denote the resulting foliation on $M_1$ by $\mathcal{F}_1'$.   Similarly, we can cut the product foliation $\mathcal{F}_2|_{\iota_2(L_2\times I)}$ open along $\iota_2(\gamma_2 \times I)$ and reglue through a map 
given by ${\rm id} \times h''_1$, thus producing a new foliation $\mathcal{F}_2'$ on $M_2$. Since a co-oriented foliation is taut if and only if it does not contain a dead end component \cite[\S 6.3]{CC1}, it follows that both $\mathcal{F}_i'$ remain taut. Moreover, by construction, the holonomy maps of $\mathcal{F}_i'$ over $b_i$ are conjugate to each other, thus completing the proof.

\begin{figure}[ht]
 \centering
 \includegraphics[scale=0.4]{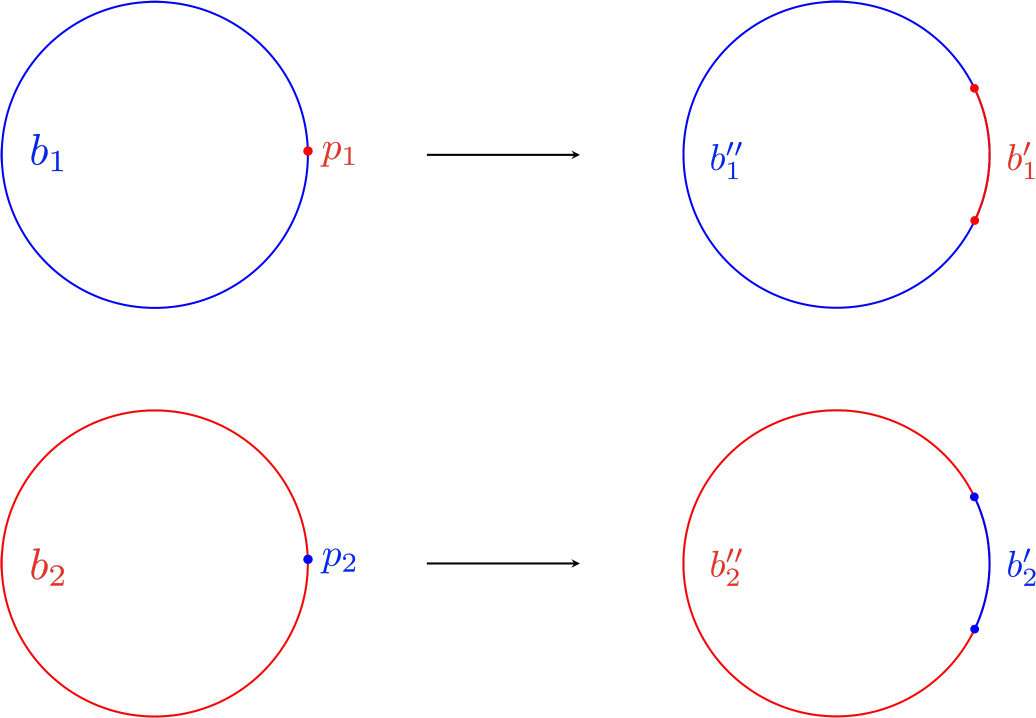}
\caption{In the picture, $p_i = a_i\cap b_i$, $i=1,2$. Circles are oriented counterclockwise and arcs on the circles are equipped with the induced orientations. Holonomy maps of $\mathcal{F}_i'\cap \partial M_i$, for $i=1,2$ over arcs of the same colour are conjugate after the alteration. }
\label{fig: holonomy case II}
\end{figure}

Suppose next that some $L_i$ is compact. Since a co-oriented taut foliation cannot have a compact separating leaf (otherwise, it would contain a dead end component), it follows that $L_i$ must be non-separating. Therefore, by the definition of longitudinal slope, the components of $\partial L_i$ are of slope $[\lambda_i]$. So $i = 2$.

Suppose next that some $L_i$ is compact. Since a co-oriented taut foliation cannot have a compact separating leaf (otherwise it would contain a dead-end component), it follows that $L_i$ must be non-separating. Therefore, by the definition of longitudinal slope, the components of $\partial L_i$ are of slope $[\lambda_i]$. Hence $i = 2$.

If $\partial L_2$ is connected, the assumption that $M_2$ is boundary-incompressible implies that $g(L_2)>0$. There is a standard construction we can 
then use to obtain a taut foliation on $L_2\times I$ tangent to $L_2\times \{0,1\}$ and transverse to $\partial L_2\times I$ so that the holonomy map of the induced foliation along $\partial L_2\times I$ is conjugate to 
$h_1''$ (see, for example, \cite[Lemma 3.1, Lemma 3.2]{Li1}). Replace the product foliation on $\iota_2(L_2\times I)$ with it to produce a foliation $\mathcal{F}_2'$ on $M_2$ with the desired holonomy. Since $L_1$ is noncompact, we can change $\mathcal{F}_1$ to $\mathcal{F}_1'$ so that the holonomy of $\mathcal{F}_1'$ over $b_1'$ is conjugate to $h_2''$ using an infinite ray starting at a point on $a_1$ as before. This deals with the case that $\partial L_2$ is connected. 

If $\partial L_2$ has more than one component, let $c$ be a properly embedded arc in $L_2$ connecting $a_2$ and a different boundary component $a_2'$ of $L_2$. Cut the product foliation $\mathcal{F}_2|_{\iota_2(L_2\times I)}$ open along $\iota_2(c \times I)$ and reglue using ${\rm id} \times 
h''_1: c\times I\rightarrow c\times I$, so that the holonomy map of the resulting foliation, denoted by $\mathcal{F}_2'$, is conjugate to $h_1''$ over $b_2'$. (Since $c$ intersects $a_2'\subset \partial L_1$ nontrivially, we have also changed the holonomy of the foliation over $b_2''$.) Use $\bar{h}_2''$ to denote the new holonomy map of the induced foliation on $\partial M_2$ over $b_2''$ and replace $\mathcal{F}_1$ by a new foliation $\mathcal{F}_1'$ on $M_1$ so that the holonomy of $\mathcal{F}_1' \cap \partial M_1$ over $b_1'$ is conjugate to $\bar{h}_2''$. This completes the proof. 
\end{proof}

\begin{corollary} 
\label{cor: ctf detd implies nls detd}
If $M$ is a knot manifold, then the set of rational slopes in $\mathcal{D}_{CTF}(M)$ is contained in $\mathcal{D}_{NLS}(M)$.
\end{corollary}

\begin{proof}
Assume that $[\beta]$ is a rational slope in $\mathcal{D}_{CTF}(M)$ and let $N$ be the twisted $I$-bundle over the Klein bottle. Fix a homeomorphism $f: \partial N \to \partial M$ which identifies $[\lambda_N]$ with $[\beta]$. Then $W = N \cup_f M$ admits a co-oriented taut foliation by Theorem \ref{thm: fln gluing}, so is not an $L$-space. Since $N$ is an irreducible, boundary-incompressible $HF$-generalised solid torus, Corollary \ref{cor: gluing defn of nls detection} implies that $[\beta] \in \mathcal{D}_{NLS}(M)$, which completes the proof.  
\end{proof}

Theorem \ref{thm: fln gluing}  leads to a characterisation of foliation-detection analogous to that 
of $NLS$-detection given by Corollary \ref{cor: gluing defn of nls detection}.

A {\it $CTF$-generalised solid torus} is a knot manifold $U$ such that $\mathcal{D}_{CTF}(U) = \{[\lambda_U]\}$. Examples include all Seifert fibred $HF$-generalised solid tori such as the twisted $I$-bundle over the Klein bottle. 

\begin{corollary} 
\label{cor: gluing defn of ctf detection}
Suppose that $M$ is a knot manifold and $W = M \cup_f U$ is a manifold obtained by gluing a hyperbolic $CTF$-generalised solid torus $U$ to $M$ using a homeomorphism $f: \partial U \to \partial M$ which identifies the longitudinal slope $[\lambda_U]$ of $U$ with a rational slope $[\alpha] \in \mathcal{S}(M)$. Then $[\alpha]$ is foliation-detected in $M$ if and only 
if $W$ admits a co-oriented taut foliation. 
\end{corollary}

\begin{proof}
Theorem \ref{thm: fln gluing} gives the forward direction of the corollary. For the reverse direction, suppose that $W$ admits a co-oriented taut foliation $\mathcal{F}$.  Since $\mathcal{F}$ is taut,  it cannot have any separating compact leaves. Hence up to isotopy, we can assume that $\mathcal{F}$ is transverse to $\partial M = \partial U$ by \cite{Th}.  We can assume, moreover, that there are no boundary-parallel annular leaves in either $\mathcal{F}_M = \mathcal{F} \cap M$ or $\mathcal{F}_{U} = 
\mathcal{F} \cap U$ (\cite{BR}).  

Suppose that $\mathcal{F}_U$ is not taut. Then there are compact leaves $L_1, L_2, \ldots , L_k$ of $\mathcal{F}_U$ whose union is the frontier of a compact submanifold $U_0$ of $U$. Further, up to changing the co-orientation on $\mathcal{F}|_U$, $U_0$ lies to the positive side of each $L_i$ (\cite[Corollary 6.3.4]{CC1}). The co-orientation of $\mathcal{F}_U$ gives an inward pointing nowhere vanishing vector field over $U$ so it follows from the Poincar\'e-Hopf theorem for manifolds with boundary that $\sum_i \chi(L_i) = 0$. The irreducibility and $\partial$-incompressibility of $U$ imply that $\chi(L_i)\leq 0$ for all $i$ and so $\chi(L_i)=0$ for all $i$. That is, $L_i$ is either a torus, or an annulus. The first possibility is ruled out by the tautness of $\mathcal{F}$ while the second contradicts the fact that $U$ is hyperbolic. Thus $\mathcal{F}_{U}$ is taut. Since $U$ is a $CTF$-generalised solid torus, \cite[Corollary 3.1]{Gab3} rules out the possibility that  $\mathcal{F}_{\partial U} = \mathcal{F}_U\cap \partial U$ contains a Reeb annulus. Then $\mathcal{F}_{\partial U}$ is a suspension foliation and therefore $\mathcal{F}_{U}$ detects a slope on $\partial U$, which is necessarily $[\lambda_{U}]$. 

It remains to show that $\mathcal{F}_M$ is taut. If this is false, then as before there would be a finite set of toral and annular leaves $L_1, L_2, \ldots , L_k$ of $\mathcal{F}_M$ whose union is the frontier of a compact submanifold $M_0$ of $M$. Further, up to changing the co-orientation on $\mathcal{F}_M$, $M_0$ lies to the positive side of each $L_i$. This means that $M_0$ is a {\it dead-end component} of $\mathcal{F}_M$; no transversal to $\mathcal{F}_M$ which enters it can exit it.  If each $L_i$ is disjoint from $\partial M$ then $M_0$ is a submanifold of $\mbox{int}(M)$ and therefore a dead-end component of $\mathcal{F}$ in $W$, which contradicts its tautness. So some $L_i$, say $L_1$, intersects $\partial M$ nontrivially. Since $L_1$ lies on the boundary of a dead-end, there is no closed transversal intersecting it. On the other hand, the fact that $\mathcal{F}_{\partial M} = \mathcal{F}_{\partial U}$ is a suspension foliation implies that there is a closed transversal to 
$\mathcal{F}_M$ lying on $\partial M$ which intersects $L_1$ non-trivially. This contradiction shows that $\mathcal{F}_M$ is taut and therefore $[\lambda_{U}] = [\alpha] \in \mathcal{D}_{CTF}(M)$. 
 
\end{proof}

We currently do not have an example of a hyperbolic $CTF$-generalised solid torus. However, although the hyperbolicity condition was used in the proof of Corollary \ref{cor: gluing defn of ctf detection}, we expect that it can be removed.

\subsection{Foliation-detection of multislopes} 
\label{subsec: fln-detn multislopes}
Suppose that $M$ is a compact, irreducible $3$-manifold whose boundary is a union $\partial M = T_1 \sqcup T_2 \sqcup \cdots \sqcup T_r$ of incompressible tori and $[\alpha] = ([\alpha_1], \ldots , [\alpha_r]) \in \mathcal{S}(M)$ is a multislope. We say that $[\alpha]$ is {\it $CTF$-detected} if there is a co-oriented taut foliation on $M$ which is transverse to $\partial M$ and intersects each $T_i$ in a suspension foliation of slope $[\alpha_i]$. A multislope $[\alpha]$ is {\it strongly $CTF$-detected} if there is a co-oriented taut foliation on $M$ which is transverse to $\partial M$ and intersects each $T_i$ in a linear foliation of slope $[\alpha_i]$. 

The set of $CTF$-detected multislopes is denoted $\mathcal{D}_{CTF}(M)$ and the set of strongly $CTF$-detected multislopes by $\mathcal{D}_{CTF}^{str}(M)$. 

Corollary \ref{cor: gluing defn of ctf detection} can be generalised to provide a characterisation of $CTF$-detected rational multislopes. The proof, which is entirely analogous to that of Corollary \ref{cor: gluing defn of ctf detection}, is left to the reader, though we note that it depends on the mild generalisation of Theorem \ref{thm: fln gluing} detailed in the last three paragraphs of the proof of Theorem \ref{thm: general * gluing}. 

\begin{proposition}
\label{prop: gluing defn of multislope ctf detection}
Suppose that $M$ is a compact, irreducible $3$-manifold whose boundary is a union $\partial M = T_1 \sqcup T_2 \sqcup \cdots \sqcup T_r$ of incompressible tori and $[\alpha] = ([\alpha_1], \ldots , [\alpha_r]) \in \mathcal{S}_{rat}(M)$ is a rational multislope. Let $W = M \cup_{T_1} U_1 \cup_{T_2} \cdots \cup_{T_r} U_r,$ where each $U_i$ is a hyperbolic $CTF$-generalised solid torus glued to $M$ along $T_i$ so that $[\lambda_{U_i}] = [\alpha_i]$. Then $[\alpha] \in \mathcal{D}_{CTF}(M)$ if and only if $W$ admits a co-oriented taut foliation. 
\end{proposition}

\section{Order-detection and gluing knot manifolds}
\label{sec: ord gps}
Much of this section is based on \cite{BC2}, where the various forms of order-detection of slopes are defined (weak, regular and strong) and analysed. 

\subsection{Generalities on left-orders} 
\label{subsec: gens on los}
A {\it left-order} $\mathfrak{o}$ on a non-trivial group $G$ can be described as either a total order $<_{\mathfrak{o}}$ on $G$ invariant under 
left-multiplication or as a semigroup $P(\mathfrak{o}) \subset G$, called 
the {\it positive cone} of $\mathfrak{o}$, which satisfies $G = P(\mathfrak{o}) \sqcup \{1\} \sqcup P(\mathfrak{o})^{-1}$. 
A group is called {\it left-orderable} if it is non-trivial and admits a left-order. 

Let $\mathfrak{o}$ be a left-order on a group $G$. An $\mathfrak{o}$-{\it 
convex} subset of $G$ is a subset $C \subseteq G$ such that if $k, h \in C$ and $g \in G$ satisfy $k <_\mathfrak{o} g <_\mathfrak{o} h$, then $g \in C$. 
It is a simple exercise to verify that if $C$ is a proper $\mathfrak{o}$-convex subgroup of $G$ and $g \in P(\mathfrak{o}) \setminus C$, then  
$$g^{-1} <_\mathfrak{o} c <_\mathfrak{o} g$$
for all $c \in C$. Thus proper $\mathfrak{o}$-convex subgroups of $G$ are $\mathfrak{o}$-bounded above and below. 

We say that a subset $A$ of $G$ is $\mathfrak{o}$-{\it cofinal} if 
$$G= \{g \in G \; | \; \mbox{ there are $a_1, a_2 \in A$ such that } a_1 \leq_\mathfrak{o} g \leq_\mathfrak{o} a_2 \}$$ 
 We say that an element $g \in G$ is $\mathfrak{o}$-cofinal if the cyclic 
subgroup $\langle g \rangle$ of $G$ it generates is $\mathfrak{o}$-cofinal. No element of a proper $\mathfrak{o}$-convex subgroup of $G$ is 
$\mathfrak{o}$-cofinal.

The set $LO(G)$ of left-orders on $G$ becomes a compact, Hausdorff, totally disconnected space when endowed with the Sikora topology \cite{Si} which, moreover, is metrisable when $G$ is countable. Setting 
$$P(g \cdot \mathfrak{o}) = g P(\mathfrak{o}) g^{-1}$$
determines an action of $G$ on $LO(G)$ by homeomorphisms. 

To each left-order $\mathfrak{o}$ on a countable group $G$ we associate its {\it dynamic realisation} 
$$\rho_{\mathfrak{o}}: G \to \mbox{Homeo}_+(\mathbb R),$$ 
which is faithful and well-defined up to conjugation in $\mbox{Homeo}_+(\mathbb R)$ (cf. \cite[Proposition 2.1]{Navas} or \cite[\S 3]{BC2}). 

The map $G \to \mathbb R, g \mapsto \rho_\mathfrak{o}(g)(0)$, is injective, unbounded above and below, and determines an order isomorphism between $(G, \mathfrak{o})$ and $(\mathcal{O}_{\rho_\mathfrak{o}}(o), <)$, where $\mathcal{O}_{\rho_\mathfrak{o}}(o)$ is the orbit of $0 \in \mathbb R$ under $\rho_{\mathfrak{o}}$ and $<$ is the induced order from $\mathbb R$. Hence the action on $\mathbb R$ induced by $\rho_{\mathfrak{o}}$ is nontrivial, 
i.e. there are no global fixed points, and an element $g \in G$ is $\mathfrak{o}$-cofinal if and only if $\rho_{\mathfrak{o}}(g)$ is fixed point free. Equivalently, $g \in G$ is $\mathfrak{o}$-cofinal if and only if $\rho_{\mathfrak{o}}(g)$ is conjugate in $\mathrm{Homeo}_+(\mathbb{R})$ to translation by $\pm 1$, which we denote by $\hbox{sh}(\pm 1)$.  

For a topological space $X$ we will often write $LO(X)$ for $LO(\pi_1(X))$. (The base point will be understood.)  

\subsection{Left-orders on $\mathbb Z^2$ and order-detection in knot manifolds}
\label{subsec: lo z2 and o-d}
The basic properties of positive cones imply that given a left-ordering $\mathfrak{o}$ of $\mathbb Z^2$, there is a line $L({\mathfrak{o}}) \subset \mathbb R^2$ uniquely determined by the fact that all elements of  $\mathbb Z^2$ which lie to one side of it are $\mathfrak{o}$-positive and all elements lying to the other are $\mathfrak{o}$-negative. The element $[L(\mathfrak{o})]$ of the projective space of $P^2(\mathbb R)$ determined by $L({\mathfrak{o}})$ is called the {\it slope} of $\mathfrak{o}$. We say that $[L(\mathfrak{o})]$ is {\it rational} if $L(\mathfrak{o}) \cap \mathbb Z^2 \cong \mathbb Z$ and {\it irrational} otherwise. 

The reader will verify that each element of $\mathbb Z^2 \setminus L(\mathfrak{o})$ is $\mathfrak{o}$-cofinal, while no element of $\mathbb Z^2 \cap L(\mathfrak{o})$ is. Indeed, $\mathbb Z^2 \cap L(\mathfrak{o})$ is an $\mathfrak{o}$-convex subgroup of $\mathbb Z^2$. 
Further, if $\mathfrak{o} \in LO(\mathbb Z^2)$ and $\alpha \in \mathbb Z^2 \setminus \{0\}$ is primitive, the following statements are equivalent
\vspace{-.2cm}
\begin{itemize}

\item $[\alpha] = [L(\mathfrak{o})]$;

\vspace{.2cm} \item $\alpha$ is not $\mathfrak{o}$-cofinal;

\vspace{.2cm} \item  $\langle \alpha \rangle$ is an $\mathfrak{o}$-convex subgroup of $\mathbb Z^2$.   

\end{itemize}

\begin{definition}
We say that a slope $[\alpha]$ on the boundary of a knot manifold $M$ is {\it weakly $LO$-detected} if there is a left-order $\mathfrak{o} \in LO(M)$ such that $[\alpha] = [L(\mathfrak{o}|_{\pi_1(\partial M)})]$. Given such an $\mathfrak{o}$ we say that $\mathfrak{o}$ {\it weakly $LO$-detects} $[\alpha]$. 
\end{definition}

\begin{definition}
We say that a slope $[\alpha]$ on the boundary of a knot manifold $M$ is {\it $LO$-detected} if there is a left-order $\mathfrak{o} \in LO(M)$ such that $[\alpha] = [L((\gamma \cdot \mathfrak{o})|_{\pi_1(\partial M)})]$ for each $\gamma \in \pi_1(M)$. Given such an $\mathfrak{o}$ we say that it {\it $LO$-detects} $[\alpha]$. 
\end{definition}
The invariance of the slope weakly order-detected over the $\pi_1(M)$-orbit of $\mathfrak{o}$ in the definition of order-detection may not seem intuitive, but this is precisely what is needed to invoke the Bludov-Glass theorem \cite[Theorem A]{BG} in the proof of the gluing theorem for left-orders \cite[Theorem 1.3]{BC2}.

We note that (weakly) order-detected slopes can be either rational or irrational. 

The set $\mathcal{D}_{LO}(M)$ of order-detected slopes is contained in the set $\mathcal{D}_{LO}^{wk}(M)$. Both are shown to be closed in $\mathcal{S}(M)$ in \cite{BC2}.

\cite{BC2} also defines the notion of a slope on the boundary of a knot manifold $M$ being {\it strongly $LO$- (or order-) detected}. In particular, a rational slope $[\alpha]$ on $\partial M$ is strongly $LO$-detected if and only if the fundamental group of the Dehn filled manifold $M(\alpha)$ has a left-orderable quotient. It is shown in \cite[Corollary 8.3]{BC2} that the set of strongly $LO$-detected slopes $\mathcal{D}^{str}_{LO}(M)$ is contained in $\mathcal{D}_{LO}(M)$.

\begin{example}
\label{exam: long is ord-det} 
Let $\lambda_M$ be the longitudinal slope of a knot manifold $M$. Since $H_1(M(\lambda_M))$ has quotient $\mathbb Z$, $[\lambda_M]$ is in  $\mathcal{D}^{str}_{LO}(M)$ and hence in $\mathcal{D}_{LO}(M)$.
\end{example}

\subsection{Boundary-cofinality and order-detection}
\label{subsec: bdry cofinality} 
A left-order which weakly order-detects a slope on the boundary of a knot manifold $M$ does not necessarily order-detect it; invariance of slope under the $\pi_1(M)$-action on $LO(M)$ may not hold. Here we discuss a situation where invariance holds automatically, a fact that we will exploit in the proofs of our order-detection results.

\begin{definition}
A left-order $\mathfrak{o} \in LO(M)$ is called {\it boundary-cofinal} if $\pi_1(\partial M)$ is $\mathfrak{o}$-cofinal. 
\end{definition}

\begin{lemma}
\label{lemma: bdry-cofinal}
Let $M$ be a knot manifold, $\mathfrak{o} \in LO(M)$ and $\rho_\mathfrak{o}$ its dynamic realisation. Then the following statements are equivalent.

$(1)$ $\mathfrak{o}$ is boundary-cofinal.  

$(2)$ The action of $\pi_1(\partial M)$ on the reals determined by $\rho_\mathfrak{o}|_{\pi_1(\partial M)}$ is fixed point free.  

$(3)$ Up to conjugation of $\rho_\mathfrak{o}$ in $\mbox{Homeo}_+(\mathbb R)$, there is a primitive element $\beta \in \pi_1(\partial M)$ such that $\rho_\mathfrak{o}(\beta) = \mbox{sh}(1)$ $($and therefore $\rho_\mathfrak{o}(\pi_1(\partial M)) \leq \mbox{Homeo}_\mathbb Z(\mathbb R)$$)$. 

\end{lemma}

\begin{proof}
It follows from the properties of left-orders on $\pi_1(\partial M)$ and the discussion in the penultimate paragraph of \S \ref{subsec: gens on los} that statement (1) implies statement (3) and statement (2) implies statement (1). That statement (3) implies (2) is obvious.  
\end{proof}

\begin{remark} 
\label{rem: bdry cof iff fpf}
It is shown in \cite{BC2} that if $M$ is a knot manifold then every $\mathfrak{o} \in LO(M)$ is boundary-cofinal if and only if any non-trivial action of $\pi_1(M)$ on $\mathbb R$ by orientation-preserving homeomorphisms restricts to a non-trivial action of $\pi_1(\partial M)$ on $\mathbb R$. The latter condition has arisen in the work of Nie on non-left-orderable Dehn filling (\cite{Nie}). 
\end{remark}

\begin{lemma}
\label{lemma: weak implies regular for bdry cofinal} 
Suppose that $M$ is a knot manifold and $\mathfrak{o} \in LO(M)$ is boundary-cofinal. Then $[L(\mathfrak{o}|_{\pi_1(\partial M)})] = [L((\gamma \cdot \mathfrak{o})|_{\pi_1(\partial M)})]$ for each $\gamma \in \pi_1(M)$. Hence $\mathfrak{o}$ order-detects $[L(\mathfrak{o}|_{\pi_1(\partial M)})] $. 
\end{lemma}

\begin{proof}
By Lemma \ref{lemma: bdry-cofinal} we can suppose that if $\rho_\mathfrak{o}$ is the dynamic realisation of $\mathfrak{o}$, then $\rho_\mathfrak{o}(\pi_1(\partial M)) \leq \mbox{Homeo}_\mathbb Z(\mathbb R)$ and there is a primitive $\beta_0 \in \pi_1(\partial M)$ such that $\rho_\mathfrak{o}(\beta_0) = \mbox{sh}(1)$. 

As we mentioned above, the map $\pi_1(M) \to \mathbb R, \gamma \mapsto \rho_\mathfrak{o}(\gamma)(0)$, is injective, unbounded above and below, and determines an order isomorphism between $(G, \mathfrak{o})$ and $(\mathcal{O}_{\rho_\mathfrak{o}}(o), <)$, where $\mathcal{O}_{\rho_\mathfrak{o}}(o)$ is the orbit of $0 \in \mathbb R$ under $\rho_{\mathfrak{o}}$ and $<$ is the induced order from $\mathbb R$. 

Set $[\alpha] = [L(\mathfrak{o}|_{\pi_1(\partial M)})]$ and suppose that $\gamma \in \pi_1(M)$. By the properties of positive cones of left-orders on $\pi_1(\partial M)$ (see \S \ref{subsec: lo z2 and o-d}), the identity $[L((\gamma \cdot \mathfrak{o})|_{\pi_1(\partial M)})] = [\alpha] $ will hold if we can show that 
$$P(\mathfrak{o}|_{\pi_1(\partial M)}) \setminus [\alpha]  \subseteq P((\gamma \cdot \mathfrak{o})|_{\pi_1(\partial M)}) \setminus [L((\gamma \cdot \mathfrak{o})|_{\pi_1(\partial M)})]$$ 
By the definition of $\gamma \cdot \mathfrak{o}$, $1 <_{\gamma \cdot \mathfrak{o}} \beta$ if and only if $1 <_\mathfrak{o} \gamma^{-1} \beta \gamma$, so we are reduced to verifying that if $1 <_\mathfrak{o} \beta \in \pi_1(\partial M) \setminus [\alpha]$, then $1 <_\mathfrak{o} \gamma^{-1} \beta \gamma$. We consider the cases $\gamma >_\mathfrak{o} 1$ and $\gamma <_\mathfrak{o} 1$ separately. 

If $\gamma >_\mathfrak{o} 1$, the $\mathfrak{o}$-cofinality of $\beta$ allows us to choose an integer $n > 0$ such that $\gamma <_\mathfrak{o} \beta^n$. Then $1 <_\mathfrak{o} \gamma^{-1}  \beta^n$ and therefore $1 <_\mathfrak{o} (\gamma^{-1} \beta^n) \gamma = (\gamma^{-1} \beta \gamma)^n$. Thus $1 <_\mathfrak{o} \gamma^{-1} \beta \gamma$. 

If $\gamma <_\mathfrak{o} 1$, choose $n > 0$ such that $\beta^{-n} <_\mathfrak{o}  \gamma$. Then $1 <_\mathfrak{o} \gamma^{-1} (\beta^{n} \gamma)= (\gamma^{-1} \beta \gamma)^{n}$. Then as before, $1 <_\mathfrak{o} \gamma^{-1} \beta \gamma$, which completes the proof.
\end{proof}

We can determine the slope that is order-detected by a boundary-cofinal left-order $\mathfrak{o}$ using the  translation number quasimorphism $\tau: \mbox{{\rm Homeo}}_{\mathbb Z}(\mathbb R) \to \mathbb R$ (\cite[\S 5]{Ghys}) as follows. 

Conjugate $\rho_\mathfrak{o}$ so that $\rho_\mathfrak{o}(\pi_1(\partial M)) \leq \mbox{Homeo}_\mathbb Z(\mathbb R)$. It is known that $\tau$ restricts to a homomorphism on abelian subgroups of $\mbox{{\rm Homeo}}_{\mathbb Z}(\mathbb R)$ so we have a linear map
$$\tau_\mathbb R: H_1(\partial M; \mathbb R) = \pi_1(\partial M) \otimes \mathbb R \xrightarrow{\;\; (\tau \circ \rho_\mathfrak{o}) \otimes 1_\mathbb R\;\;} \mathbb R,$$
necessarily non-zero (Lemma \ref{lemma: bdry-cofinal}). Hence the kernel of $\tau_\mathbb R$ determines a line in $H_1(\partial M; \mathbb R)$ and therefore a slope $[\mbox{ker}(\tau_\mathbb R)]$ on $\partial M$.

\begin{lemma}
\label{lemma: slope by tau}
Suppose that $M$ is a knot manifold and $\mathfrak{o} \in LO(M)$ is boundary cofinal. If $\rho_\mathfrak{o}$ is a dynamic realisation of $\mathfrak{o}$ for which $\rho_\mathfrak{o}(\pi_1(\partial M)) \leq \mbox{Homeo}_\mathbb Z(\mathbb R)$, then the slope order-detected by $\mathfrak{o}$ is $[\mbox{ker}(\tau_\mathbb R)]$. In particular, if $\alpha \in \pi_1(\partial M) \setminus \{1\}$ and $\tau(\rho_\mathfrak{o}(\alpha)) = 0$, then the slope order-detected by $\mathfrak{o}$ is $[\alpha]$. 
\end{lemma}

\begin{proof}
For each $\beta \in \pi_1(\partial M) \setminus \mbox{ker}(\tau_\mathbb R)$ we have $\tau(\rho_\mathfrak{o}(\beta)) \ne 0$. Then $\tau(\rho_\mathfrak{o}(\beta)) > 0$ if and only if $\rho_\mathfrak{o}(\beta)(0) > 0$ (\cite[\S 5]{Ghys}). On the other hand, we have an order isomorphism between $(G, \mathfrak{o})$ and $(\mathcal{O}_{\rho_\mathfrak{o}}(o), <)$, so $\tau(\rho_\mathfrak{o}(\beta)) > 0$ if and only if $\beta >_\mathfrak{o} 1$. Similarly $\tau(\rho_\mathfrak{o}(\beta)) < 0$ if and only if $\beta <_\mathfrak{o} 1$. Thus the elements of $\pi_1(\partial M)$ which lie to one side of $[\mbox{ker}(\tau_\mathbb R)]$ are positive and those which lie to the other side are negative, which completes the proof. 
\end{proof}

\subsection{Order-detection via representations}
\label{subsec: orders via reps} 

The following proposition will be used to convert Theorem \ref{thm: universal circle fixed point} into results on order-detection. 

\begin{proposition}
\label{prop: rep to ord det}
Suppose that $M$ is a knot manifold and $\alpha \in \pi_1(\partial M)$ represents a rational slope. If $\rho: \pi_1(M) \to \mbox{Homeo}_+(\mathbb R)$ is a homomorphism such that $\rho(\alpha)$ has a fixed point but $\rho(\pi_1(\partial M))$ does not, then $[\alpha]$ is order-detected in $M$. 
\end{proposition}

\begin{proof} 
Let $\mbox{Fix}(\alpha) \subseteq \mathbb R$ denote the fixed point set of $\rho(\alpha)$ and choose $x_0 \in \mbox{Fix}(\alpha)$. Any $\beta \in \pi_1(\partial M) \setminus \langle \alpha \rangle$ commutes with $\alpha$ so $\mathcal{O}_\beta(x_0)$, which we define to be $\{\rho(\beta^n)(x_0) \; | \; n \in \mathbb Z\}$, is contained in $\mbox{Fix}(\alpha)$. If $\mathcal{O}_\beta(x_0)$ were bounded above or below,  its supremum, respectively infinum, would be fixed by both $\rho(\alpha)$ and $\rho(\beta)$ and therefore by $\rho(\pi_1(\partial M))$, contrary to our assumptions. Thus $\mathcal{O}_\beta(x_0)$ is unbounded above and below. Up to replacing $\beta$ by $\beta^{-1}$, which doesn't change $\mathcal{O}_\beta(x_0)$, we can suppose that $\rho(\beta)$ is strictly increasing. Then the map $\mathbb Z \to \mathcal{O}_\beta(x_0) \subset \mathbb R, n \mapsto \rho(\beta^n)(x_0)$ is injective, order-preserving, and unbounded above and below. It follows that $\mathcal{O}_\beta(x_0) $ is discrete in $\mathbb R$.

Now $\rho(\beta)$ acts freely on $\mathcal{O}_\beta(x_0)$, as otherwise $\rho(\pi_1(\partial M))$ has a fixed point, and from this we see that it acts freely on 
the reals. Then $\rho(\beta)$ is conjugate to translation by $\pm 1$ and therefore we can suppose that $\rho(\pi_1(\partial M)) \leq \mbox{{\rm Homeo}}_{\mathbb Z}(\mathbb R)$. Composing this restriction with translation number $\tau: \mbox{{\rm Homeo}}_{\mathbb Z}(\mathbb R) \to \mathbb R$  yields a non-trivial homomorphism $\tau_\rho: \pi_1(\partial M) \to \mathbb R$ (\cite[\S 5]{Ghys}). Since $\rho(\alpha)$ has fixed points, $\tau_\rho(\alpha) = 0$, and so the kernel of $\tau_\rho$ is $\langle \alpha \rangle$. 

As a non-trivial subgroup of the left-orderable group $\pi_1(M)$ (\cite[Theorem 1.1]{BRW}), $\mbox{Stab}_\rho(x_0)$ is left-orderable. Moreover, the total order on $\mathcal{O}(x_0) = \{\rho(\gamma)(x_0) \; | \; \gamma \in \pi_1(M)\} \equiv \pi_1(M)/\mbox{Stab}_\rho(x_0)$ induced from $\mathbb R$ is invariant under the $\rho$-action. We claim that there is a left-order $\mathfrak{o} \in LO(M)$ for which 
\begin{equation}
\label{eqn: ineq} 
\gamma_1 \leq_{\mathfrak{o}} \gamma_2 \; \Rightarrow \; \rho(\gamma_1)(x_0) \leq \rho(\gamma_2)(x_0)
\end{equation}
To see this, let $\mathfrak{o}_0$ be a left-order on $\mbox{Stab}_\rho(x_0)$ and consider the left-order on $\pi_1(M)$ given by 
$$\gamma_1 <_{\mathfrak{o}} \gamma_2 \hbox{ if and only if } \left\{ 
\begin{array}{cl} 
1 <_{\mathfrak{o}_0} \gamma_1^{-1} \gamma_2 & \hbox{when } \gamma_1^{-1} \gamma_2 \in \mbox{Stab}(x_0) \\ 
& \\
\rho(\gamma_1)(x_0) < \rho(\gamma_2)(x_0) & \hbox{when } \gamma_1^{-1} \gamma_2 \not \in \mbox{Stab}(x_0)
\end{array}
\right.$$
It is easy to verify that (\ref{eqn: ineq}) holds and since $\mathcal{O}_\beta(x_0)$ is unbounded above and below, it implies that each $\beta \in \pi_1(\partial M) \setminus \langle \alpha \rangle$ is $\mathfrak{o}$-cofinal. Thus $\mathfrak{o}$ is boundary-cofinal. Further, (\ref{eqn: ineq}) implies that 
$\mbox{Stab}_\rho(x_0)$ is a proper $\mathfrak{o}$-convex subgroup of $\pi_1(M)$. Since $\langle \alpha \rangle =  \mbox{Stab}_\rho(x_0) \cap \pi_1(\partial M)$, it is a proper $\mathfrak{o}|_{\pi_1(\partial M)}$-convex subgroup of $\pi_1(\partial M)$. Thus $\mathfrak{o}$ weakly order-detects $[\alpha]$ and since it is boundary-cofinal, Lemma \ref{lemma: weak implies regular for bdry cofinal} implies that it order-detects $[\alpha]$.
\end{proof}

Representations of the fundamental group of a knot manifold $M$ with values in $\mbox{Homeo}_+(\mathbb R)$ often arise as lifts $\widetilde \rho: \pi_1(M) \to \mbox{{\rm Homeo}}_{\mathbb Z}(\mathbb R) \leq \mbox{Homeo}_+(\mathbb R)$ of representations $\rho: \pi_1(M) \to \mbox{{\rm Homeo}}_+(S^1)$. 

\begin{proposition} 
\label{prop: divide tau lambda}
Let $M$ be a knot manifold whose longitude is integrally null-homologous and $\rho: \pi_1(M) \to \mbox{{\rm Homeo}}_+(S^1)$ a homomorphism such that $\rho(\pi_1(\partial M))$ fixes a point of $S^1$. Suppose that $\rho$ lifts to a homomorphism $\widetilde \rho_0: \pi_1(M) \to \mbox{{\rm Homeo}}_{\mathbb Z}(\mathbb R)$ such that $\tau(\widetilde \rho_0(\lambda_M))$, necessarily an integer, is non-zero. Then if $d \in \mathbb Z$ divides $\tau(\widetilde \rho_0(\lambda_M))$, any  rational slope of distance $d$ from 
$[\lambda_M]$ is order-detected. 
\end{proposition}

\begin{proof}
Fix a dual class $\mu$ to $\lambda_M$ in $H_1(\partial M) = \pi_1(\partial M)$ and suppose that $\alpha = d \mu + q \lambda_M$ is a primitive class, where $d$ divides $m = \tau(\widetilde \rho_0(\lambda_M))$. We claim that there is a lift $\widetilde \rho$ of $\rho$ for which $\widetilde \rho(\alpha)$ has a fixed point. 

The centre of $\mbox{{\rm Homeo}}_{\mathbb Z}(\mathbb R)$ is the group $K$ of integer translations of the reals and is also the kernel of the homomorphism 
$\mbox{{\rm Homeo}}_{\mathbb Z}(\mathbb R)  \to \mbox{{\rm Homeo}}_+(S^1)$. Since $\lambda$ is integrally null-homologous in $H_1(M)$, $\mu$ generates $H_1(M)/T_1(M) \cong \mathbb Z$, and hence there is a homomorphism $\phi: \pi_1(M) \to K$ such that $\phi(\mu) = - \tau(\tilde{\rho}_0(\mu)) - q(m/d)$. Then $\widetilde \rho(\gamma) = \phi(\gamma) \widetilde \rho_0(\gamma)$ defines a lift of $\rho$ for which $\tau(\widetilde \rho(\mu)) = -q(m/d)$. Since  $\phi(\lambda_M) = 0$ we have $\tau(\widetilde \rho(\lambda_M)) = \tau(\widetilde \rho_0(\lambda_M)) = m$. Then 
$$\tau(\widetilde \rho(\alpha)) = d \tau(\widetilde \rho(\mu)) + q \tau(\widetilde \rho(\lambda_M)) =  d (-q(m/d)) + qm = 0$$ 
Thus $\widetilde \rho(\alpha)$ has a fixed point in $\mathbb R$. On the other hand, $\widetilde \rho(\pi_1(\partial M))$ has no fixed points in $\mathbb R$ since the translation number of $\widetilde \rho(\lambda_M)$ is non-zero. Proposition \ref{prop: rep to ord det} now shows that $[\alpha]$ is order-detected. 
\end{proof}

\subsection{The gluing theorem for left-orders}

The left-order case of Theorem \ref{thm: * gluing} was proven in \cite[Theorem 1.3]{BC2}. 

\begin{theorem}  {\rm (Boyer-Clay)}
\label{thm: LO-gluing}
Suppose that $W = M_1 \cup_f M_2$ where $M_1, M_2$ are knot manifolds and $f: \partial M_1 \xrightarrow{\; \cong \; } \partial M_2$. 
If $f$ identifies an order-detected slope on the boundary of $M_1$ with an order-detected slope on the boundary of $M_2$, then $W$ has a left-orderable fundamental group. 
\end{theorem}
Note that it is not required that the identified slopes be rational. 

This leads to a characterisation of order-detection analogous to that of $NLS$-detection given by Corollary \ref{cor: gluing defn of nls detection}.

An {\it $LO$-generalised solid torus} is a knot manifold $U$ such that $[L(\mathfrak{o}|_{\pi_1(\partial U)})] = [\lambda_U]$ for each $\mathfrak{o} \in LO(U)$. Equivalently, $\mathcal{D}_{LO}^{wk}(U) = \{[\lambda_U]\}$. Examples include the twisted $I$-bundle over the Klein bottle and, as we shall see in Proposition \ref{prop: hfgst and detection}, the hyperbolic knot manifold $v2503$.

\begin{corollary}  
\label{cor: gluing defn of lo detection}
Suppose that $M$ is a knot manifold and $W = M \cup_f U$ is a manifold obtained by gluing an $LO$-generalised solid torus $U$ to $M$ using a homeomorphism $f: \partial U \to \partial M$ which identifies the longitudinal slope $[\lambda_U]$ of $U$ with a rational slope $[\alpha] \in \mathcal{S}(M)$. Then $[\alpha]$ is $LO$-detected in $M$ if and only if $W$ has a left-orderable fundamental group. 
\end{corollary}

\begin{proof}
The forward direction follows immediately from Theorem \ref{thm: LO-gluing} and Example \ref{exam: long is ord-det}. 

For the reverse direction, suppose that $\mathfrak{o}$ is a left-order on $\pi_1(W)$ and let $\mathfrak{o}_M$ be its restriction to $\pi_1(M)$. Set $T = \partial M = \partial U$. It suffices to show that $[L((\gamma \cdot \mathfrak{o}_M)|_{\pi_1(T)})] = [\alpha]$ for each $\gamma \in \pi_1(M) \leq \pi_1(W)$. But as $(\gamma \cdot \mathfrak{o}_M)|_{\pi_1(T)} = (\gamma \cdot \mathfrak{o})|_{\pi_1(T)} = ((\gamma \cdot \mathfrak{o})_{\pi_1(U)})|_{\pi_1(T)}$, the assumption that $U$ is an $LO$-generalised solid torus implies that $[L((\gamma \cdot \mathfrak{o}_M)|_{\pi_1(T)})] = [\lambda_U] = [\alpha]$, which completes the proof. 
\end{proof}

\subsection{Order-detection of multislopes} 
Suppose that $M$ is a compact, irreducible $3$-manifold whose boundary is a union $\partial M = T_1 \sqcup T_2 \sqcup \cdots \sqcup T_r$ of incompressible tori and $[\alpha] = ([\alpha_1], \ldots , [\alpha_r]) \in \mathcal{S}(M)$ is a multislope. We say that $[\alpha]$ is {\it $LO$-detected} if there is a left-order $\mathfrak{o} \in LO(M)$ such that for each $\gamma \in \pi_1(M)$, the slope of $(\gamma \cdot \mathfrak{o})|_{\pi_1(T_i)}$  is $[\alpha_i]$. The set of order-detected multislopes is denoted $\mathcal{D}_{LO}(M)$.

Here is an extension of Corollary \ref{cor: gluing defn of lo detection} to multislopes, which is an order-detection analogue of the characterisation of foliation-detected rational multislopes given in Proposition \ref{prop: gluing defn of multislope ctf detection}. 

\begin{proposition}
\label{prop: gluing defn of multislope lo detection}
Suppose that $M$ is a compact, irreducible $3$-manifold whose boundary is a union $\partial M = T_1 \sqcup T_2 \sqcup \cdots \sqcup T_r$ of incompressible tori and $[\alpha] = ([\alpha_1], \ldots , [\alpha_r]) \in \mathcal{S}_{rat}(M)$ is a rational multislope. Let $W = 
M \cup_{T_1}U_1 \cup_{T_2} \cdots \cup_{T_r} U_r,$ 
where each $U_i$ is an $LO$-generalised solid torus glued to $M$ along $T_i$ so that $[\lambda_{U_i}] = [\alpha_i]$. Then $[\alpha] \in \mathcal{D}_{LO}(M)$ if and only if $W$ has a left-orderable fundamental group. 
\end{proposition}

\begin{proof}
The proof of the reverse direction of the proposition is entirely analogous to the case $r = 1$ dealt with in Corollary \ref{cor: gluing defn of ctf detection}. The forward direction follows from a multislope extension of Theorem \ref{thm: LO-gluing}, and is dealt with in \cite[\S 5]{BC2}. 
\end{proof}

\section{The generalised gluing theorem for manifolds with multiple boundary components}
\label{sec: detn and gluing} 

\subsection{Multislope detection and $HF$-generalised solid tori} 
\label{subsec: hfgst}
Recall the family of $HF$-generalised solid tori from \S \ref{subsec: NLS 
detn}, which we defined as knot manifolds all of whose non-longitudinal Dehn fillings are $L$-spaces. Proposition \ref{prop: hfgst implies small detd set} shows that for such manifolds, $\mathcal{D}_{NLS}(M) = \{[\lambda_M]\}$.

The following proposition is due to Sarah and Jake Rasmussen (cf. Proposition 1.9 of \cite{RR}) and Gillespie (Corollary 2.7 of \cite{Gi}). 

\begin{proposition} {\rm (Gillespie, Rasmussen-Rasmussen)} 
\label{prop: char hfgst}
Let $M$ be a rational homology solid torus. Then the following statements 
are equivalent.

$(1)$ $M$ is an HF-generalised solid torus. 

$(2)$ Any Thurston norm minimizing surface representing a generator of $H_2(M, \partial M)$ is planar and $M$ admits an $L$-space filling. 

\end{proposition}

\begin{proposition}
\label{prop: hfgst and detection}
The hyperbolic $HF$-generalised solid torus $v2503$ is an $LO$-generalised solid torus.
\end{proposition}

\begin{proof}
Let $U_0 = v2503$, and let $\lambda$ be the longitude of $U_0$. To prove that $U_0$ is an $LO$-generalised solid torus we invoke work of Nie \cite{Nie}.

Let $\mathfrak{o}$ be a left-order on $\pi_1(U_0)$. To show that $U_0$ is an $LO$-generalised solid torus we must show that the slope weakly detected by $\mathfrak{o}$ is $[\lambda]$. Let $\rho_\mathfrak{o}: \pi_1(U_0) \to \mbox{Homeo}_+(\mathbb R)$ be the dynamic realisation of $\mathfrak{o}$; \cite[Lemma 6.1]{Nie} shows that $\rho_\mathfrak{o}(\pi_1(\partial U_0))$ has no fixed points in $\mathbb R$. It follows from Lemma \ref{lemma: bdry-cofinal} that $\mathfrak{o}$ is boundary-cofinal, and hence by Lemma \ref{lemma: weak implies regular for bdry cofinal} that the slope weakly detected by $\mathfrak{o}$  is order-detected by $\mathfrak{o}$. By Lemma \ref{lemma: bdry-cofinal} we can also assume that $\rho_\mathfrak{o}(\pi_1(\partial U_0)) \leq \mbox{Homeo}_{\mathbb Z}(\mathbb R)$. Therefore, by Lemma \ref{lemma: slope by tau}, to show that the slope order-detected by $\mathfrak{o}$ is $[\lambda]$, it suffices to show that $\tau(\rho_\mathfrak{o}(\lambda)) = 0$.

As a first case, suppose that $\tau(\rho_\mathfrak{o}(\lambda)) > 0$, so that $\rho_\mathfrak{o}(\lambda)$ is strictly increasing. 

Following Nie (\cite{Nie}), let $\mathcal{M}(\lambda^{-1})$ be the root-closed, conjugacy-closed submonoid of $\pi_1(U_0)$ generated by $\lambda^{-1}$. In other words, $\mathcal{M}(\lambda^{-1})$ is the minimal subset of $\pi_1(U_0)$ containing $\lambda^{-1}$ which has the following properties:
\vspace{-.2cm}
\begin{itemize}

\item $\eta, \gamma \in \mathcal{M}(\lambda^{-1}) \Rightarrow \eta \gamma \in \mathcal{M}(\lambda^{-1})$;

\vspace{.2cm} \item $\gamma \in \mathcal{M}(\lambda^{-1})$ and $\eta \in \pi_1(U_0) \Rightarrow \eta \gamma \eta^{-1} \in \mathcal{M}(\lambda^{-1})$; 

\vspace{.2cm} \item $\gamma^n \in \mathcal{M}(\lambda^{-1})$ for some integer $n > 0 \Rightarrow \gamma \in \mathcal{M}(\lambda^{-1})$.

\end{itemize}
We  can construct $\mathcal{M}(\lambda^{-1})$ inductively by taking $\mathcal{M}_0 = \{\lambda^{-1}\}$ and defining $\mathcal{M}_{k+1}$ to be the set of all conjugates, positive roots, and products of elements in $\mathcal{M}_{k}$. Then $\mathcal{M}(\lambda^{-1}) = \cup_k  \mathcal{M}_k$.

By assumption, $\tau(\rho_\mathfrak{o}(\lambda^{-1})) < 0$, so $\rho_\mathfrak{o}(\lambda^{-1})$ is strictly decreasing. Further, since being strictly decreasing is preserved by the composition of such homeomorphisms, taking their positive roots, or conjugating them in $\mbox{Homeo}_+(\mathbb R)$, each element of $\rho_\mathfrak{o}(\mathcal{M}(\lambda^{-1}))$ is strictly decreasing. In particular, $\lambda \not \in \mathcal{M}(\lambda^{-1})$. But Proposition \ref{prop: hfst} shows that there is an essential planar surface properly embedded in $U_0$ whose oriented boundary consists of like-oriented curves of slope $[\lambda_{U_0}]$, which implies that $\lambda \in \mathcal{M}(\lambda^{-1})$. Thus the assumption that $\tau(\rho_\mathfrak{o}(\lambda)) > 0$ leads to a contradiction. 

A similar contradiction is reached if $\tau(\rho_\mathfrak{o}(\lambda)) < 
0$, for then each element of  $\rho_\mathfrak{o}(\mathcal{M}(\lambda^{-1}))$ is strictly increasing, including $\rho_\mathfrak{o}(\lambda)$. Hence 
$\tau(\rho_\mathfrak{o}(\lambda)) = 0$, which completes the 
proof.
\end{proof}

\begin{remark}
We expect that all $HF$-generalised solid tori are $LO$-generalised solid tori. This holds in the graph manifold case, including Seifert fibred manifolds, by the solution of the $L$-space conjecture for this class of manifolds (\cite{BC1}, \cite{HRRW}). In general, the argument used in the above analysis of $U_0 = v2503$ shows that to prove that an $HF$-generalised solid torus $U$ is an $LO$-generalised solid torus, it suffices to show that if $\rho: \pi_1(U) \to \mbox{Homeo}_+(\mathbb R)$ is a representation without a global fixed point, then $\rho|_{\pi_1(\partial U)}$ has no global fixed point and $\tau(\rho_\mathfrak{o}(\lambda_U)) = 0$. The former is shown to hold when there is some slope on $\partial U$ which is not weakly order-detected (\cite[\S 4]{BC2}). More concretely, work of Nie (\cite{Nie}) shows that it holds for the exteriors $U$ of Berge-Gabai knots standardly embedded in a Heegaard solid torus of $S^1 \times S^2$, which are $HF$-generalised solid tori by \cite[Proposition 7.8]{RR}. 
\end{remark}

The following proposition combines the definition of $NLS$-detection of  rational multislopes with Propositions \ref{prop: gluing defn of multislope ctf detection} , \ref{prop: gluing defn of multislope lo detection}, and \ref{prop: hfgst and detection}. 

\begin{proposition} 
\label{prop: ctf implies lo detd genl}
Suppose that $M$ is a compact, connected, orientable, irreducible, boundary-incompressible $3$-manifold with $\partial M = T_1 \sqcup \cdots \sqcup T_r$ a non-empty union of tori. Let $[\alpha] \in \mathcal{S}_{rat}(M)$ and set $W 
= M \cup_{T_1} U_1 \cup_{T_2} \cdots \cup_{T_r} U_r$ where each $U_i$ 
is a hyperbolic $\ast$-generalised solid torus which is glued to $M$ along $T_i$ so that $[\lambda_{U_i}]$ is identified with $[\alpha_i]$. Then $W$ has property $\ast$ if and only if $[\alpha] \in \mathcal{D}_{*}(M)$. 

\end{proposition}

\begin{remark}
\label{rem: detd slopes the same}
The $L$-space conjecture predicts that $\mathcal{D}_{NLS}(M) = \mathcal{D}_{LO}(M)$: If we take $U$ to be the hyperbolic $HF$- and $LO$-generalised solid torus $v2503$, then it contends that the $W$ of Proposition \ref{prop: ctf implies lo detd genl} is $CTF$ if and only if it is $LO$. If we also knew that $v2503$ was a $CTF$-generalised solid torus, then the $L$-space conjecture would imply that $\mathcal{D}_{NLS}(M) = \mathcal{D}_{LO}(M) = \mathcal{D}_{CTF}(M)$. The equality of these three sets has been verified when $M$ is a graph manifold (\cite{BC1}, \cite{HRRW}). We note, moreover, that it is expected that their equality is equivalent to the $L$-space conjecture for closed, connected, orientable, irreducible, toroidal $3$-manifolds. By Corollary \ref{cor: ctf detd implies nls detd}, showing that $v2503$ is a $CTF$-generalised solid torus is equivalent to showing that no irrational slope on its boundary is $CTF$-detected.

\end{remark}

\subsection{Gluing coherence and the generalised gluing theorem}
Here we show that Theorem \ref{thm: * gluing} has a natural extension to 
more general gluings which we will use to prove Theorem \ref{thm: results 
GL conjecture}. 

Let $W$ be an irreducible rational homology $3$-sphere and $T_1, T_2, \ldots, T_m$ a disjoint family of essential tori in $W$ which split it into pieces $M_1, M_2, \ldots, M_{m+1}$.  

For each family of rational slopes $([\alpha_1], [\alpha_2], \ldots, [\alpha_m]) \in \mathcal{S}_{rat}(T_1) \times \mathcal{S}_{rat}(T_2) \times \cdots \times \mathcal{S}_{rat}(T_m)$ and $1 \leq j \leq m+1$, let 
$$[\alpha(j)] = ([\alpha_{i_1}], [\alpha_{i_2}], \ldots, [\alpha_{i_{r}}]) \in \mathcal{S}_{rat}(M_j),$$ 
where $\partial M_j = T_{i_1} \sqcup T_{i_2} \sqcup \cdots \sqcup T_{i_{r}}$.
For $\ast \in \{NLS, LO, CTF\}$ we say that $([\alpha_1], [\alpha_2], \ldots, [\alpha_m])$ is {\it $\ast$-gluing coherent} if $[\alpha(j)] \in \mathcal{D}_{\ast}(M_j)$ for each $j$. 

\begin{theorem} 
\label{thm: general * gluing}
Let $W = \cup_j M_j$ be an irreducible rational homology $3$-sphere expressed as a union of submanifolds $M_1, M_2, \ldots, M_{m+1}$ along a disjoint family of essential tori $T_1, T_2, \ldots, T_m$. Suppose that $([\alpha_1], [\alpha_2], \ldots, [\alpha_m]) \in \mathcal{S}_{rat}(T_1) \times \mathcal{S}_{rat}(T_2) \times \cdots \times \mathcal{S}_{rat}(T_m)$ is $\ast$-gluing coherent. Then $W$ has property $\ast$. 
\end{theorem}

\begin{proof}
We induct on $m$, the case $m = 1$ being Theorem \ref{thm: * gluing}.

First consider the cases that $\ast = NLS$ or $LO$. Suppose that $m \geq 2$ and split $W$ open along $T_1$, say $W = M_1' \cup_{T_1} M_2'$. Up to reindexing the tori $T_i$ we can suppose that $T_2, \ldots, T_s \subset \mbox{int}(M_1')$ and $T_{s+1}, \ldots, T_m \subset 
\mbox{int}(M_2')$. Set $W_k = M_k' \cup U_0$ where $U_0$ is attached to 
$M_k'$ by a homeomorphism which identifies the longitudinal slope of $U_0$ with $[\alpha_1]$. Then as $([\alpha_{2}], \ldots, [\alpha_s])$ is $\ast$-gluing coherent in $W_1$ and $([\alpha_{s+1}], \ldots, [\alpha_m])$ is $\ast$-gluing coherent in $W_2$, our induction hypothesis implies that both $W_1$ and $W_2$ have property $\ast$, so $[\alpha_1] \in \mathcal{D}_{\ast}(M_k')$ for both $k$. The case $m = 1$ now implies that $W$ has property $\ast$.  

We proceed differently in the case that $\ast = CTF$, owing to the fact that currently we do not have an example of a hyperbolic $CTF$-generalised solid torus. 

After reindexing the $M_i$ we can suppose that $M_1$ and $M_2$ are incident to $T_1$. We claim that $M_2' = M_1 \cup_{T_1} M_2$ admits a co-oriented taut foliation which detects the slope $[\alpha_i]$ on each $T_i$ lying in its boundary. If this is the case, our induction hypothesis applied to the union $W = M_2' \cup M_3 \cup  \cdots M_{m+1}$ along the tori $T_2, T_3, \ldots , T_m$ will complete the proof. 

The verification that $M_2'$ admits a co-oriented taut foliation as claimed is entirely similar to the proof of Theorem \ref{thm: fln gluing}. Fix, for $i = 1, 2$, a co-oriented taut foliation $\mathcal{F}_i$ on $M_i$ which detects the multislope $[\alpha(i)]$ and let $L_i$ be a leaf of $\mathcal{F}_i$ whose intersection with $T_{1}$ is non-empty and has a compact component $a_i$, necessarily a simple closed curve of slope $[\alpha_{1}]$. Since $W$ is a rational homology $3$-sphere, $[\alpha_1]$ is not rationally null-homologous in at least one of $M_1, M_2$, say in $M_1$. 

For each $i$ perform a Denjoy blow up on the leaf $L_i$ and as before, call the new foliation $\mathcal{F}_i$. Alter the holonomy of $\mathcal{F}_2|_{T_{1}}$ in the thickened neighbourhood of $a_2$ as in the proof of Theorem \ref{thm: fln gluing}, so that it matches that of $\mathcal{F}_1|_{T_{1}}$ outside the thickened neighbourhood of $a_1$. Since $[\alpha_{1}]$ is not rationally null-homologous in $M_1$, $L_1$ is either non-compact or there is a $j \geq 2$ such that $L_1 \cap T_{j} \ne \emptyset$. In the former case we alter the holonomy of $\mathcal{F}_1|_{T_{1}}$ in the thickened neighbouhood of $L_1$ along a properly embedded ray in $L_1$ based on $a_1$ so that it matches that of $\mathcal{F}_2|_{T_{1}}$ outside the thickened neighbourhood of $a_2$. In the latter, we do a similar operation along a properly embedded arc in $L_1$ connecting $a_1$ to a point in $T_{j} \; (\ne T_1)$. The foliations on the $M_i$ now match up along $T_1$ and glue together to give a  co-oriented taut foliation on $M_2'$. Further note that the operations of thickening leaves and changing holonomy do not change the slopes detected on the components of $\partial M_i \setminus T_1$, so there is a foliation on $M_2'$ as required.  
\end{proof}

\begin{remark}
That the converse of Theorem \ref{thm: general * gluing} holds when $\ast = NLS$ is an easy consequence of \cite[Theorem 13]{HRW1}. We expect it to hold when $\ast = LO$ or $CTF$.  
\end{remark}

\section{Order-detection and the dynamics of universal circle actions} 
\label{sec: univ circles and knot mflds}

The goal of this section is to prove Theorem \ref{thm: universal circle fixed point}.

\subsection{Universal circle actions on knot manifolds}
\label{subsec: existence uca}
The usual setting for a universal circle representation is a closed $3$-manifold which admits a taut foliation with hyperbolic leaves (Thurston; Calegari-Dunfield \cite{CaD}). Here we justify their existence in the case of the tautly foliated rational homology solid torus $M$. 

Throughout we take $M \not \cong S^1 \times D^2$ to be an irreducible rational homology solid torus whose longitude is integrally null-homologous in $M$. 
We denote by $F$ an oriented Thurston norm-minimizing surface with connected 
boundary representing a generator $[F]$ of $H_2(M, \partial M) \cong \mathbb Z$. Then 
$\chi(F)<0$. 

Our fundamental groups will be based at a point $p \in \partial F$, though we will usually suppress $p$ from our notation. 

By \cite[Theorem 5.5]{Gab1}, $M$ admits a co-orientable taut foliation $\mathcal{F}_M$ of finite depth for which $F$ is a compact leaf. Moreover, the foliation is transverse to $\partial M$ and the induced foliation on $\partial M$ has no Reeb component. We assume that leaves of $\mathcal{F}_M$ are smooth. See \cite{Ca1, KR}.

Consider the double $W = M \cup_{\rm id} -M$ of $M$ and let $\mathcal{F} = \mathcal{F}_M \cup_{\rm id} -\mathcal{F}_M$ be the resulting finite depth foliation on $W$ which contains $S=F\cup_{\rm id} -F$ as a compact leaf. Glue the foliations $\mathcal{F}_M$ and $\mathcal{F}_{-M}$ together so that leaves of $\mathcal{F}$ are smooth along $T$.  Let $\iota: (W,\mathcal{F}) \rightarrow (W,\mathcal{F})$ be the involution on $(W,\mathcal{F})$ switching the two components of $W\setminus T$ and preserving the foliation $\mathcal{F}$, where $T=\partial M \subset W$ and $\iota|_T = {\rm id}_T$. We claim that there is an $\iota$-invariant Riemannian metric on $W$, for which the leaves of $\mathcal{F}$ become hyperbolic with respect to the induced metric. 

To see this, first let $g_0$ be any smooth Riemannian metric on $W$ for which $\iota$ is an isometry. This can be easily obtained by doubling the Riemannian manifold $(M,g_M)$ where $g_M$ is a Riemannian metric on $M$, for which leaves of $\mathcal{F}_M$ are perpendicular to $\partial M$. Now suppose that $\mathfrak{m}$ is a transverse invariant measure of $\mathcal{F}$ (\cite[Definition 7.5]{Ca2}). Since $M$ is a rational homology solid torus, $b_1(W)=1$ and therefore by \cite[Theorem 6.3]{Pl}, all leaves 
of $\mathcal{F}$ that are in the support of $\mathfrak{m}$ are compact. It then follows from the definition of the Euler characteristic of an invariant measure (\cite[Definition 7.7]{Ca2}) that the Euler characteristic of $\mathfrak{m}$ is a positively-weighted sum of the Euler characteristics of compact leaves, which is negative by our assumption. (Note that all compact leaves are homologous to $F$ and norm-minimizing.) By \cite[Theorem 4.3]{Can}, all leaves of $\mathcal{F}$ are (conformally) hyperbolic. By amalgamating the uniformization map for each leaf of $\mathcal{F}$, we obtain a uniformization map $\eta: W\rightarrow \mathbb{R}$. By \cite[Proposition 4.11, Theorem 4.14]{Can}, $\eta$ is continuous over $M$ and smooth in tangential directions. Set $g = \eta g_0$. Then with the induced leafwise metric from $g$, each leaf of $\mathcal{F}$ becomes hyperbolic with constant curvature $-1$. We note that since $(W, \mathcal{F}, g_0)$ is invariant under $\iota$, so is $\eta$, and therefore $g$ is symmetric with respect to $\iota$. That is, $\iota$ is an isometry of $(W,g)$. 

With the metric $g$ on $W$, we can construct universal circle actions $\rho: \pi_1(W)\rightarrow {\mbox{Homeo}}_+(S_{univ}^1)$  by \cite{CaD}\footnote{Though Calegari and Dunfield assume their manifolds are atoroidal, all that is needed for the construction is that the foliation's leaves are hyperbolic.}. We define the universal circle action of $\pi_1(M)$ associated to $\mathcal{F}_M = \mathcal{F}|_M$ to be the restriction $\rho_M = \rho|_{\pi_1(M)}$. 

The restriction of $\rho_M$ to $\pi_1(F)$ is semi-conjugate to a discrete 
faithful representation $\rho_0: \pi_1(F) \to PSL(2, \mathbb R)$ for which $\rho_0(\lambda)$ is hyperbolic. If $\rho_M(\lambda)$ is conjugate to $\rho_0(\lambda)$, it would have a unique attracting fixed point $u_0$ on $S^1$, and so as each $\gamma \in \pi_1(\partial M)$ commutes with $\lambda$, $\rho_M(\pi_1(\partial M))$ fixes $u_0$. Showing that we still have a fixed point when $\rho_M(\lambda)$ is only semiconjugate to $\rho_0(\lambda)$ requires that we review some of the key ingredients in the construction of $\rho$, which we do in \S \ref{subsec: universal circle actions} 
and \S \ref{subsec: sawblades markers sections}.

\subsection{From special sections to the universal circle action}
\label{subsec: universal circle actions}
Let $(\widetilde{\mathcal{F}}, \tilde{g})$ be the lift of $(\mathcal{F}, g)$ to the universal cover $\pi: \widetilde{W}\rightarrow W$ and fix a lift 
$\tilde p \in \widetilde{W}$ of the base point $p \in T$. We use $\tilde p$ to identify $\pi_1(W; p)$ with the group of deck transformations of this cover, in the usual way.

The {\it leaf space} $\mathcal{L}$ of $\widetilde{\mathcal{F}}$ is the quotient of $\widetilde{W}$ obtained by collapsing each leaf of $\widetilde{\mathcal{F}}$ to a point. In general,  $\mathcal{L}$ is a connected, oriented, possibly non-Hausdorff, $1$-manifold \cite[Corollary D.1.2]{CC2}.

Since $\mathcal{F}$ is taut, each leaf $L$ of $\widetilde{\mathcal{F}}$ is simply connected and isometric to the hyperbolic plane $\mathbb{H}^2$ with respect to the pull-back metric $\tilde g$. Hence each leaf $L$ has a 
circle at infinity, or ideal boundary, which we denote by $\partial_\infty L$, and, with a properly defined topology, these circles piece together 
to form an oriented $S^1$-bundle $E_\infty \to \mathcal{L}$. The deck transformations of the cover $\widetilde W$ are isometries when restricted to the leaves of $\widetilde{\mathcal{F}}$ and from this it can be shown that they induce an action of $\pi_1(W)$ on $E_\infty$ by bundle maps.

The construction of a universal circle action from a taut foliation is originally due to Thurston, though largely unwritten. A different construction was given by Calegari and Dunfield in \cite{CaD}. This same construction with  some minor alterations is also described in \cite{Ca2}. 

The key element in the definition of a universal circle $S^1_{univ}$ of $\mathcal F$ is the construction of a set of {\it special sections} $\mathcal{S}$ of the circle bundle $E_\infty \to \mathcal{L}$, described below. It turns out that given any $x\in E_\infty$, there is a unique special section  $\sigma_x: \mathcal{L}\rightarrow E_\infty$ in $\mathcal{S}$ passing through $x$, where $x$ is called the {\it base point of $\sigma_x$}. Of course sections based at different 
points in $E_\infty$ may be the same. The fundamental group $\pi_1(W)$ acts on $\mathcal{S} = \{\sigma_x : x\in E_\infty\}$ via $g\cdot \sigma_x 
= \sigma_{g\cdot x}$.

It is possible to endow $\mathcal{S}$ with a circular order which is natural with respect to the constructions. For instance, the action of $\pi_1(W)$ on $\mathcal{S}$ is order-preserving. Also, as each leaf $L$ of $\widetilde{\mathcal{F}}$ represents a point of $\mathcal{L}$, there is surjective evaluation map $e_L: \mathcal{S} \rightarrow \partial_\infty L$ which has degree one (i.e. point inverses are order intervals) and, apart from collapsing, is order preserving with respect to the circular order on the circle $\partial_\infty L$ determined by its orientation. 

We can embed $\mathcal{S}$, endowed with the order topology, as a dense subspace of an oriented universal circle $S^1_{univ}$. Further, the action of $\pi_1(W)$ on $\mathcal{S}$ extends to a {\it universal circle action} $\rho: \pi_1(W)\rightarrow {\mbox{Homeo}}_+(S_{univ}^1)$.

\subsection{Sawblades, markers, and special sections}
\label{subsec: sawblades markers sections}
Briefly, special sections are obtained by assembling {\it markers} via a left-most up, right-most down rule, and markers are constructed from {\it 
sawblades}.

By \cite[Lemma 2.4]{CaD}, there exists a separation constant $\epsilon >0$ such that each leaf of $\widetilde{\mathcal{F}}$ is quasi-isometrically 
embedded in its $\epsilon$-neighbourhood in $\widetilde{W}$. A {\it marker} for $\widetilde{\mathcal{F}}$ is a map $m: I\times \mathbb{R}^+ \rightarrow \widetilde{W}$ such that $m(\{s\} \times \mathbb{R}^+)$ is a geodesic ray in a leaf of $\widetilde{\mathcal{F}}$ for each $s \in [0,1]$ and $m(I\times \{t\})$ is a tight transversal of $\widetilde{\mathcal{F}}$ of length  less than $\epsilon/3$ for each $t\in \mathbb{R}^+$ \cite[Definition 5.1]{CaD}\footnote{We will not (explicitly) use the condition on the lengths of tight transversals in the exposition below. However, this requirement is the key to many crucial properties of the markers in constructing special sections.}. Each geodesic ray defines a unique point on the ideal boundary of a leaf and hence a marker defines a partial section of the circle bundle $E_\infty$. We  call the image of the partial section given by a marker $m$ the {\it endpoints} or the {endpoint set} of the marker, and denote it by $e(m) \subset E_\infty$. As in \cite{CaD}, we will sometimes abuse notation by referring to the endpoint set as a marker. 

To construct the set of special sections $\mathcal{S}$, one needs a $\pi_1(W)$-invariant set of markers $\mathfrak{M}$ whose endpoint sets intersect the circle at infinity of each leaf of $\widetilde{\mathcal{F}}$ in a 
dense subset. The existence of such a set of markers is the content of the leaf pocket theorem \cite[Theorem 5.2]{CaD}. There is a lot of flexibility in the construction of $\mathfrak{M}$ and we briefly describe how to do this given a simple closed geodesic in a leaf of $\mathcal{F}$. For the details, see \cite[\S 5.3]{CaD}. Also see \cite[Construction 7.19, Construction 7.20]{Ca2}. 

\begin{figure}[ht]
\includegraphics[scale=0.7]{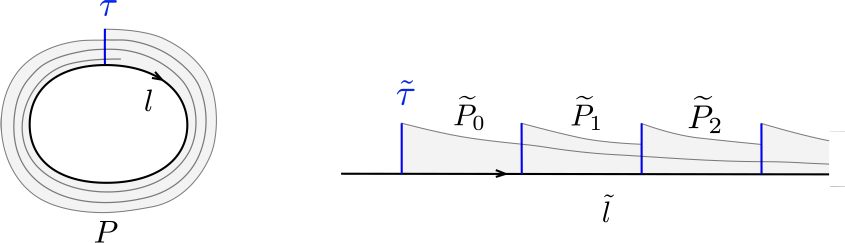}
\caption{The shaded region in the left-hand picture is a sawblade. The right-hand picture shows a chain of lifts $\widetilde{P}_n$ of a sawblade which 
lies along a ray contained in a copy of $\tilde{l}$. There is a unique marker starting at $\tilde{\tau}$ which lies in an apriori bounded neighbourhood of the union $\widetilde{P}_\infty = \cup \widetilde{P}_n$.}
\label{fig: sawblade}
\end{figure}

Let $l$ be an oriented simple closed geodesic on a leaf $L$ of $\mathcal{F}$. Taking a sufficiently small transversal $\tau: [0,1] \rightarrow W$ of $\mathcal{F}$ with $\tau(0) \in l$ and reversing the orientation of $l$ if necessary, we may assume that the holonomy map $h_l$ along $l$ is non-expanding on $\tau$, meaning that $\tau^{-1}(h_l(\tau(s)))$ is defined for all $s\in [0,1]$ and bounded above by $s$. See Figure \ref{fig: sawblade}. Note that the holonomy map $h_l$ is trivial on $\tau$ if and only if $h_l(\tau(s))= \tau(s)$ for all $s\in [0,1]$. We slide $\tau$ along $l$  and obtain an immersed rectangular region $P: [0,1]\times[0,1]\rightarrow W$, called a {\it sawblade}, where by construction we have $P(s,0) = \tau(s)$ and $P(s,1)= h_l(P(s,0))$ for all $s\in [0,1]$. We may also assume that the length of each transversal $P|_{[0,1]\times \{t_0\}}$ for $t_0\in [0,1]$ is less than $\epsilon/3$, where $\epsilon > 0$ is the separation constant. 

Given a lift $\tilde{\tau}$ of $\tau$ to $\widetilde{W}$, there exists a unique lift $\widetilde{P}_0:[0,1]\times [0,1] \rightarrow \widetilde{W}$ of $P$ with $\widetilde{P}_0|_{[0,1]\times \{0\}} = \tilde{\tau}$. For each $n\in\mathbb{N}$, let $\widetilde{P}_n$ be the lift of $P$ satisfying $\widetilde{P}_{n-1}|_{[0,1]\times \{1\}}\subseteq \widetilde{P}_{n}|_{[0,1]\times \{0\}}$. We use $\widetilde{P}_\infty$ to denote the union of 
the $\widetilde{P}_n$ for $n\geq 0$. The bottom edge of  $\widetilde{P}_\infty$ projects to $l$ and hence is a geodesic ray on the leaf $\widetilde{\mathcal{F}}$ containing it. Moreover, one can shrink $\tau$ further so that for each $s>0$, the horizontal ray starting at $\tau(s)$ in $\widetilde{P}_\infty$ is a quasi-geodesic ray\footnote{Here the horizontal ray starting at $\tau(s)$ is defined as the intersection of $\widetilde{P}_\infty$ and the leaf of $\widetilde{\mathcal{F}}$ containing $\tau(s)$. These horizontal rays stay inside $\widetilde{P}_\infty$  because the holonomy map along $l$ is non-expanding.}.  By a leaf-wise straightening of these rays we obtain a marker $m$ with $m|_{[0,1]\times \{0\}} = \tilde{\tau}$ and the image of $m$ lies in a small neighbourhood of $\widetilde{P}_\infty$ (see \cite[\S 5.3]{CaD}).

Finally, we obtain a $\pi_1(W)$-invariant collection of markers by applying the construction to all possible lifts of $\tau$ to $\widetilde{W}$. For instance, in Figure \ref{fig: sawblade} we have a marker starting at each blue vertical transversal.

\begin{remark}{\label{rem: transversal pointing down}}
We can construct another set of markers by using a sufficiently short transversal $\tau: [0,1]\rightarrow W$ lying to the other side of $l$ (i.e. $\tau(1)\in l$). As before, we want to slide $\tau$ along $l$ in the holonomy non-expanding direction to obtain a sawblade. 
\end{remark}

\subsection{The action of peripheral subgroup on the leaf space}
\label{subsec: peripheral acts on the leaf space}
 We use $T\subset W$ to denote $\partial M$ and take $\widetilde{T}$ to be the cover of $T$ that contains our lifted base point $\tilde{p}$. Since 
$\mathcal{F}\cap T$ has no Reeb component and contains $\lambda = \partial F$ as a compact leaf, up to isotopy we may assume that a given dual curve to $\lambda$, $\mu$ say, is a closed transversal of $\mathcal{F}$ containing $p$. By an abuse of notation we will also use $\mu$ and $\lambda$ to denote the elements of $\pi_1(W)$ they carry.

Given a transversal $\nu$ of $\widetilde{\mathcal{F}}$, since $\widetilde{W}$ is simply connected, $\nu$ intersects each leaf of $\widetilde{\mathcal{F}}$ at most once. Hence, it corresponds to an embedded line in the leaf space of $\widetilde{\mathcal{F}}$, which we will denote by  $\mathcal{L}_{\nu}$. Consider $\tilde{\mu} \subset \widetilde{W}$, the cover of $\mu$ which passes through $\tilde{p}$.  Notice that the leaves in $\mathcal{L}_{\tilde{\mu}}$ are precisely those that intersect  $\widetilde{T}$ nontrivially, so as $\pi_1(T)$ acts on $\widetilde{T}$,  it also acts on $\mathcal{L}_{\tilde{\mu}}$, and we will show that it  fixes a special 
section whose base point is on $E_\infty|_{\mathcal{L}_{\tilde{\mu}}}$.

For each $n \in \mathbb Z$, let $\tilde{p}_n = \mu^n \cdot \tilde{p}$ be the lifts of $p$ along $\tilde{\mu}$ and note that they divide $\tilde \mu$ into a sequence of closed subintervals $[\tilde p_n, \tilde p_{n+1}]$. We use $\widetilde{S}_n$ to denote the cover of $S = F\cup -F$ containing $\tilde{p}_n$, and use $\tilde s_n$ to denote the corresponding point on $\mathcal{L}$. Hence $[\widetilde{s}_n, \widetilde{s}_{n+1}]$ is the subinterval of ${\mathcal{L}_{\tilde{\mu}}}$ consisting of the leaves of $\widetilde{\mathcal{F}}$ containing points of $[\tilde p_n, \tilde p_{n+1}]$.

It is easy to see that the action of $\mu$ on $\mathcal{L}_{\tilde\mu}$ identifies the interval $[\tilde{s}_i, \tilde{s}_{i+1}]$ with $[\tilde{s}_{i+1}, \tilde{s}_{i+2}]$, while $\lambda$ leaves each $[\tilde{s}_i, \tilde{s}_{i+1}]$ invariant. Moreover, if $l$ is a (closed) longitudinal leaf of $\mathcal{F}\cap T$, then $\lambda$ fixes all leaves in $\mathcal{L}_{\tilde{\mu}}$ that contain a cover of $l$ in $\widetilde{T}$. 

\begin{lemma} 
 \label{lem: action of lambda}
Let $l$ be a longitudinal leaf of $\mathcal{F}\cap T$. If $\tilde{l}$ is a component of the inverse image of $l$ in $\widetilde{T}$ and $\widetilde{L}$ is the leaf of $\widetilde{\mathcal{F}}$ containing $\tilde{l}$, then the subgroup of $\pi_1(W)$ which stabilises $\widetilde{L}$ contains $\lambda$ and acts on $\widetilde{L}$ by  hyperbolic elements of $\mbox{Isom}_+(\widetilde{{L}})$.  
 \end{lemma}
 
\begin{proof}
Our discussion above shows that $\lambda$ preserves $\widetilde{L}$. Moreover, since the deck transformation group acts on $(\widetilde{W}, \tilde{g})$ by isometries, a non-trivial element $\gamma$ of $\pi_1(W)$ which stabilises $\widetilde{L}$ is either an elliptic, parabolic, or hyperbolic isometry of  the hyperbolic plane $(\widetilde{L}, \tilde{g}|_{\widetilde{L}})$. As it acts fixed point freely, it cannot be elliptic. Nor can it be parabolic since $(W,g)$ has positive injectivity radius. This means that there exists a real number $R_0 > 0$ such that given any $x \in W$, the exponential map $\exp_x: T_xW \rightarrow W$ is injective when restricted to tangent vectors of length less than $R_0$. Pulling back to $\widetilde W$ we see that for any $\tilde x \in \widetilde L$, $d_{\widetilde{L}}(\tilde{x}, \gamma \cdot \tilde{x}) \geq d_{\widetilde W} (\tilde{x}, \gamma \cdot \tilde{x}) > 2R_0$, which rules out parabolicity. It is therefore a hyperbolic isometry.
\end{proof}

\subsection{Proof of Theorem \ref{thm: universal circle fixed point}(1)} 
By \cite[Lemma 6.11]{CaD}, the endpoints of two markers are either disjoint or can be amalgamated into one larger interval transverse to the fibres of $E_\infty$. One can take maximal unions of the overlapping endpoint sets of markers to obtain a $\pi_1(W)$-invariant collection of intervals on $E_\infty$ transverse to the circle fibres, which we call the {\it maximal markers}. By construction, special sections follow these maximal markers as long as they are defined; the complications occur when a section leaves a marker. In the latter case one needs to apply a ``leftmost up, rightmost down" rule to determine which marker to follow next, thus guaranteeing that $\mathcal{S}$ is $\pi_1$-invariant and circularly orderable (see \cite[\S 6.10 - \S 6.20]{CaD}). However, in our case, we will specify a set of markers which allows us to define special sections by simply concatenating their endpoints.

As above, we take $W=M \cup_{{\rm id}} -M$, $\mathcal{F}$ to be the foliation on $W$ obtained by doubling $\mathcal{F}_M$, $T = \partial M \subset W$, and $\iota: (W,\mathcal{F},g)\rightarrow (W,\mathcal{F},g)$ to be the involution defined in \S \ref{subsec: existence uca}.

By construction, $\iota$ preserves leaves that intersect $T$ nontrivially\footnote{In fact, all leaves of $\mathcal{F}$ intersect $T$ nontrivially.} and acts isometrically on each of these leaves. Let $l$ denote a (compact or non-compact) leaf of $\mathcal{F}\cap T$. Since it is fixed by $\iota$ pointwise, it follows from the uniqueness of geodesics that $l$ is a geodesic on the hyperbolic leaf of $\mathcal{F}$ containing it. Therefore, each leaf of $\widetilde{\mathcal{F}}\cap \widetilde{T}$ is also a geodesic of the hyperbolic leaf of $\widetilde{\mathcal{F}}$ containing it.

We first prove the existence of a global fixed point of $\rho|_{\pi_1(T)}$ in two special cases. The strategies used will allow us to establish the general statement.

\subsubsection{Special case I: $\mathcal{F}\cap T$ is a fibration by longitudinal curves.}
\label{subsec: all longitudes} 
Let $l$ be a leaf of $T\cap \mathcal{F}$.  By assumption, the holonomy along $l$ is trivial. Applying the construction in \S \ref{subsec: sawblades markers sections} to $l$, we obtain a $\pi_1(W)$-invariant set of markers $\mathfrak{M}_l$. We assume that the small transversal $\tau: [0,1]\rightarrow W$, $\tau(0)\in l$, and the sawblade $P$ used in constructing $\mathfrak{M}_l$ are on $T$.

We take $\mathfrak{M}$ to be the union of $\mathfrak{M}_l$ where $l$ ranges over all components of $\mathcal{F}\cap T$.  Since $\mathcal{F}$ is of 
finite depth, its minimal set is the union of compact leaves. In particular, they all intersect $T$. It follows from the proof of the leaf pocket theorem (\cite[Theorem 5.2]{CaD}) that the union of all endpoints of the markers in $\mathfrak{M}$ intersects the ideal boundary of each leaf of $\widetilde{\mathcal{F}}$ in a dense subset. As such, we can use $\mathfrak{M}$ to construct a set of special sections $\mathcal{S}$, from which we 
obtain a universal circle $S^1_{univ}$ and universal circle action $\rho: 
\pi_1(W)\rightarrow {\mbox{Homeo}}_+(S^1_{univ})$.

Next we identify a special section that is invariant under the action of $\pi_1(T)$, which gives us the desired point on $S^1_{univ}$ fixed under $\rho|_{\pi_1(T)}$.

Given a leaf $\widetilde{L} \in \mathcal{L}_{\tilde{\mu}}$, we let $\tilde{l} = \widetilde{L}\cap \widetilde{T}$, $L = \pi(\widetilde{L})$ and $l = \pi(\tilde{l})$, where $\pi: \widetilde{W}\rightarrow W$ is the universal covering. Let $\widetilde{P}_\infty$ denote the union of lifts of a sawblade $P \subset T$ along a positive geodesic ray contained in $\tilde{l}$ as in \S \ref{subsec: sawblades markers sections} (cf.  Figure \ref{fig: sawblade}). 

In $\mathfrak{M}$, there exists a marker $m_{\widetilde{L}}:[0,1]\times \mathbb{R}^+\rightarrow \widetilde{W}$, unique up to the action of $\lambda$, where $m_{\widetilde{L}}|_{0\times \mathbb{R}^+}$ is $\tilde l$. As noted above, leaves of $\widetilde{T}\cap \widetilde{\mathcal{F}}$ are geodesics with respect to the leafwise hyperbolic metric, so the image of $m_{\widetilde{L}}$ in fact lies in $\widetilde{P}_\infty$ (not just in a neighbourhood of $\widetilde{P}_\infty$) and each geodesic ray $m_{\widetilde{L}}|_{\{s\}\times \mathbb{R}^+}$, $s\in [0,1]$, is a horizontal ray in $\widetilde{P}_\infty$.  Because the horizontal rays in $\widetilde{P}_\infty$ are invariant under $\lambda$ and points on the ray are shifted forward by $\lambda$, the endpoint of the geodesic ray $m_{\widetilde{L}}|_{\{s\}\times \mathbb{R}^+}$ is the unique attracting fixed point of the action of $\lambda$ on the ideal boundary of the leaf containing it (cf. Lemma \ref{lem: action of lambda}).  

 Therefore, the endpoint sets of markers in the set $\{m_{\widetilde{L}}: \widetilde{L}\in \mathcal{L}_{\tilde{\mu}}\} \subset \mathfrak{M}$ overlap each other. The maximal marker obtained by taking the union of the endpoints of these makers  determines a section of $E_\infty|_{\mathcal{L}_{\tilde{\mu}}}$ whose value at each leaf in ${\mathcal{L}_{\tilde{\mu}}}$ is the attracting fixed point of $\lambda$.

 Now, let $q_n$, $n\in \mathbb{Z}$, be the attracting  fixed point of $\lambda$ on $\widetilde{S}_n$. By the discussion above, the special sections $\sigma_{q_n}$ all follow the maximal marker and hence agree over $\mathcal{L}_{\tilde{\mu}}$. It then follows from the construction of the special sections  that the $\sigma_{q_n}$ are identical over the entire $\mathcal{L}$. That is, they represent the same section in $\mathcal{S}$, which we denote by $\sigma$. Moreover, since $\lambda \cdot \sigma_{q_n} = \sigma_{\lambda \cdot q_n} = \sigma_{q_n}$
 and $\mu \cdot \sigma_{q_n} = \sigma_{\mu \cdot q_n} = \sigma_{q_{n+1}}$, we have that $\sigma$ is fixed under the action $\rho|_{\pi_1(T)}$, as required.

\subsubsection{ Special case II: $\mathcal{F}\cap T$ has a unique compact leaf} 
\label{subsubsec: one longitude}
\begin{figure}[ht]
\includegraphics[scale=0.6]{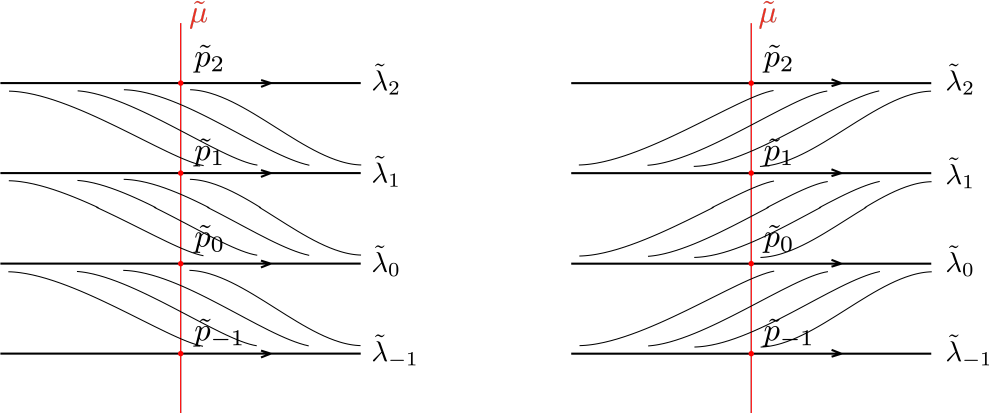}
\caption{In the figure, $\tilde{p}_0 = \tilde{p}$, $\tilde{p}_n = \mu^n\cdot \tilde{p}_0$ and $\tilde{\lambda}_n = \mu^n\cdot \tilde{\lambda}_0 = \widetilde{T}\cap \widetilde{S}_n$. }
\label{fig: foliation tilde T}
\end{figure}

In this case, $S = F\cup -F$ is the unique compact leaf of $\mathcal{F}$ and $\widetilde{\mathcal{F}}\cap \widetilde{T}$ is conjugate to one of the foliations in Figure \ref{fig: foliation tilde T}. We continue to use the notation established previously. 

Assume that the foliation $\widetilde{F}\cap \widetilde{T}$ is as depicted in the left-hand picture of Figure \ref{fig: foliation tilde T}. The argument for the case that $\widetilde{\mathcal{F}}\cap \widetilde{T}$ is as pictured on the right-hand side of Figure \ref{fig: foliation tilde T} is similar.

Since $S$ is the unique compact leaf, the set of markers obtained by applying \S \ref{subsec: sawblades markers sections} to $\lambda$ intersects the ideal boundary of every leaf in a dense subset. We pick two transversals $\tau:[0,1] \rightarrow W$ and $\tau':[0,1] \rightarrow W$ with $\tau(0)\in \lambda$ and $\tau'(1)\in \lambda$ of $\mathcal{F}$ (see Remark \ref{rem: transversal pointing down}).  By \S \ref{subsec: sawblades markers sections}, from each transversal, we obtain a $\pi_1(W)$-invariant set of markers, denoted by $\mathfrak{M}_\tau$ and $\mathfrak{M}_{\tau'}$ respectively. We take $\mathfrak{M} = \mathfrak{M}_\tau \cup \mathfrak{M}_{\tau'}$.  This guarantees that on each side of the ideal boundary of a lift $\widetilde{S}$ of the compact leaf $S$ there is a collection of markers in $\mathfrak{M}$ that intersects $\partial \widetilde{S}$ in a dense subset\footnote{This density property of the set of makers $\mathfrak{M}$ is needed in the construction of the special sections in general (see \cite[\S 6.7]{CaD} and the proof of \cite[Lemma 6.12]{CaD}), but we won't need it in our proof.}. 

In what follows, we will show the existence of a fixed point on the universal circle under the action of $\rho|_{\pi_1(T)}$ using markers in $\mathfrak{M}_\tau$. (The same argument applied to $\mathfrak{M}_{\tau'}$ shows the existence of a second fixed point.)

\begin{figure}[ht]
\begin{center}
\includegraphics[scale=0.5]{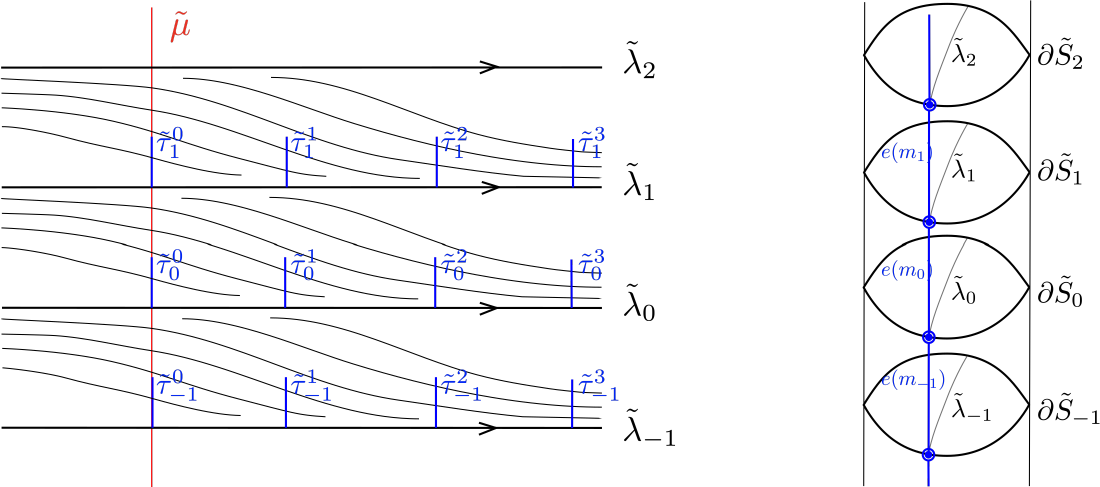}
\end{center}
\caption{For each $n\in \mathbb{Z}$,  $\tilde{\lambda}_n=\widetilde{T}\cap \widetilde{S}_n$ is the cover of $\lambda$ on $\widetilde{S}_n$ and $\{\tilde{\tau}_n^i\}_{i\in\mathbb{N}}$ are the lifts of $\tau$ along $\tilde{\lambda}_n$ satisfying that $\tilde\tau_n^0\subset \tilde\mu$ and $\lambda \tilde{\tau}_n^i = \tilde{\tau}_n^{i+1}$. The union of the maximal markers $\{e(m_n)\}_{n\in \mathbb{Z}}$ defines a section of $E_\infty|_{\mathcal{L}_{\tilde\mu}}$.}
\label{fig: markers and sections case II}
\end{figure}

For each $n\in \mathbb{Z}$, let $\tilde{\lambda}_n=\widetilde{T}\cap \widetilde{S}_n$ be the cover of $\lambda$ on $\widetilde{S}_n$ and $\{\tilde{\tau}_n^i\}_{i\in\mathbb{N}}$ be the lifts of $\tau$ along $\tilde{\lambda}_n$ satisfying that $\tilde\tau_n^0\subset \tilde\mu$ and $\lambda \tilde{\tau}_n^i = \tilde{\tau}_n^{i+1}$ (Figure \ref{fig: markers and sections case II}). By the construction described in \S \ref{subsec: sawblades markers sections}, for each $n\in \mathbb{Z}$ there exists a unique sequence of markers $\{m_n^i\}_{i\in \mathbb{N}}$, starting at $\tilde{\tau}_n^i$  satisfying the following: 
\vspace{-.2cm} 
\begin{enumerate}

\item For all $i\in \mathbb{N}$, the geodesic rays $m_n^i|_{\{0\}\times \mathbb{R}^+}$ are contained in $\tilde{\lambda}_n$.
 
\vspace{.2cm} \item For each $n$, we have $\lambda \cdot m_n^i =m_n^{i+1}$ and $e(m_n^i) \subset e(m_n^{i+1})$, where $e(m^i_n)$ is the endpoint set of the marker $m^i_n$ (see Figure \ref{fig: markers and sections case II}). 
\end{enumerate}
\vspace{-.2cm} 
By Property (2), the union $e(m_n) =\cup_i e(m_n^i)$ is a maximal marker which defines a partial section $\sigma_n: [\tilde{s}_n, \tilde{s}_{n+1}) \rightarrow E_\infty|_{[\tilde{s}_n, \tilde{s}_{n+1})}$ with $\sigma_n(\tilde{L}) = \partial \tilde{L} \cap e(m_n)$ for any $\tilde{L}\in [\tilde{s}_n, \tilde{s}_{n+1})$. Since the endpoint sets of markers are dense in the ideal boundary of every leaf, $\sigma_n$ can be extended to a section over  the closed interval $[\tilde{s}_n, \tilde{s}_{n+1}]$, denoted by $\bar \sigma_n$. Moreover, it is clear that $\bar \sigma_n(\tilde s_{n+1})$ is the attracting fixed point of $\lambda$ on $\partial_\infty \tilde{s}_{n+1}$. In fact, by the definition of the topology on $E_\infty$, the endpoint map $e: UT\widetilde{\mathcal{F}}|_{\tilde \mu}\rightarrow E_\infty|_{\mathcal{L}_{\tilde \mu}}$ is a homeomorphism, where $UT\widetilde{\mathcal{F}}$ denotes the unit tangent bundle to the foliation $\widetilde{\mathcal{F}}$ over $\widetilde{W}$. Since leaves of $\widetilde{T}\cap \widetilde{\mathcal{F}}$ are geodesics with respect to the leafwise hyperbolic metric, the section of $UT\widetilde{\mathcal{F}}|_{[\tilde p_n, \tilde p_{n+1})}$ given by $e^{-1}\circ \sigma_n$ evaluated at each point $x\in [\tilde p_n, \tilde p_{n+1}) \subset \tilde{\mu}$ is the unit vector tangent to both $\widetilde{\mathcal{F}}$ and $\widetilde{T}$, pointing in the positive $\lambda$-direction. Hence, the limit vector at $\tilde p_{n+1}$ is the positive tangent vector of $\tilde\lambda_{n+1}$.

As shown in the right-hand picture in Figure \ref{fig: markers and sections case II}, the union of the maximal markers $\{e(m_n)\}_{n\in \mathbb{Z}}$ defines a section of $E_\infty|_{\mathcal{L}_{\tilde\mu}}$, passing through $q_n$, $n\in \mathbb{Z}$, where as before, $q_n$ is the attracting  fixed point of $\lambda$ on $\widetilde{S}_n$. Hence, the sections $\sigma_{q_n}$ of $E_\infty$ all follow the set of maximal markers $\{e(m_n)\}_{n\in \mathbb{Z}}$ over $\mathcal{L}_{\tilde\mu}$ (see \cite[Lemma 6.12]{CaD}), which implies that they are the same section of $E_\infty$ by the construction of the special sections \cite[\S 6.10, \S 6.20]{CaD}. As in the previous case, this section is fixed under the action $\rho|_{\pi_1(T)}$, so it corresponds to a fixed point on $S^1_{univ}$ under the action of $\pi_1(T)$.

\subsubsection{Proof of Theorem \ref{thm: universal circle fixed point}(1): the general case} In general, the foliation $\mathcal{F}\cap T$ consists of longitudinal leaves and open annuluar regions foliated by lines. Note that there are countably many (finite or countably infinite) such line-foliated open annular regions; we denote them by $\{A_i\}_{i\in \mathbb{N}}$. 

If a compact leaf $l$ of $\mathcal{F}\cap T$ is not a component of any $\bar{A}_i\setminus A_i$, we let  $\mathfrak{M}_l$ be the set of markers constructed as in Special case I in \S \ref{subsec: all longitudes}. Suppose that $l$ is a component of $\overline{A}_i \setminus A_i$ for some $i$. We require $\mathfrak{M}_l$ to contain two sets of markers, constructed from two transversals, one on each side of $l$. On a side of $l$ where $l$ is not the limit of a sequence of compact leaves, 
we construct markers in $\mathfrak{M}_l$ as in Special Case II, \S\ref{subsubsec: one longitude}.  If there exists a sequence of compact leaves $l_n$ of $\mathcal{F}\cap T$ limiting to $l$ as $n \to \infty$, 
one takes a transversal $\tau: [0,1] \to W$ with $\tau(0)$ on $l$ and $\tau(1)$ on $l_N$ for some sufficiently large $N$ 
(so all $l_n$ with $n > N$ also intersect $\tau$). 
Then we obtain a marker by lifting the sawblade from $\tau$, 
whose endpoint set is invariant under the action of $\lambda$. 
Moreover, since $\tau(1)$ also lies on a compact leaf $l_N$ of $\mathcal{F}\cap T$, as in Special Case~I in \S\ref{subsec: all longitudes} the endpoint set of the marker 
defines a section of $E_\infty$ over a closed sub-interval of the leaf space.
 
We let $\mathfrak{M} = \cup_l \mathfrak{M}_l$, where $l$ ranges over the set of longitudinal leaves of $\mathcal{F}\cap T$. 

\begin{figure}[ht]
\centering
\includegraphics[scale=0.62]{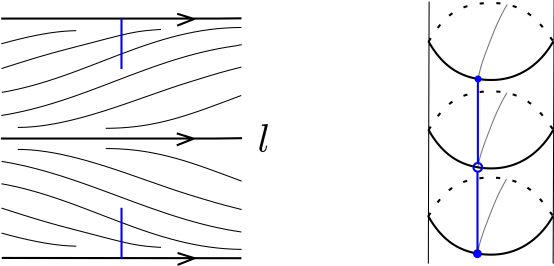}
\caption{A compact leaf $l$ of $\mathcal{F}\cap T$ is adjacent to two line-foliated annular regions, and the holonomy maps along $l$ are repelling on both sides. This results in a puncture in the union of the endpoints of the markers.}
\label{fig: puncture_marker}
\end{figure}

As before, we want to construct a $\pi_1(T)$-invariant section of $E_\infty|_{\mathcal{L}_{\tilde{\mu}}}$ whose intersection with each $\widetilde{S}_n$ is the attracting fixed point of $\lambda$ by taking the union of the endpoints of markers in $\mathfrak{M}$. By our discussion above, the only situation in which we cannot simply take the union of the endpoints of the markers in $\mathfrak{M}$ to obtain a section of $E_\infty|_{\mathcal{L}_{\tilde{\mu}}}$ 
occurs when a compact leaf $l$ of $\mathcal{F}\cap T$ is adjacent to two line-foliated annular regions, and the holonomy maps along $l$ are repelling on both sides. See Figure \ref{fig: puncture_marker}. In this case, the open ends of the set of the endpoints of the two markers coincide, which results in a puncture in the middle. Though it is  possible to show that a special section would pass through the puncture and continue to follow the markers, to simplify matters, we avoid the problem by performing a series of Denjoy blow-ups on $\mathcal{F}$ so that in the resulting foliation, no two line-foliated annular regions are adjacent. In particular, the configuration in Figure \ref{fig: puncture_marker} will never occur.

Here are the details. Let $\Lambda$ denote the union of leaves of $\mathcal{F}$ satisfying: given $L\in \Lambda$, there is a closed component of $L\cap T$ that is a boundary component of some $\overline{A}_i$. Since there are countably many leaves in $\Lambda$, we can perform a Denjoy blowup on $\Lambda$ by thickening the leaves in $\Lambda$ and foliating the I-bundles by parallel leaves \cite[Theorem 7.3]{KR}.

  \medskip
\begin{figure}[ht]
\centering
\includegraphics[scale=0.48]{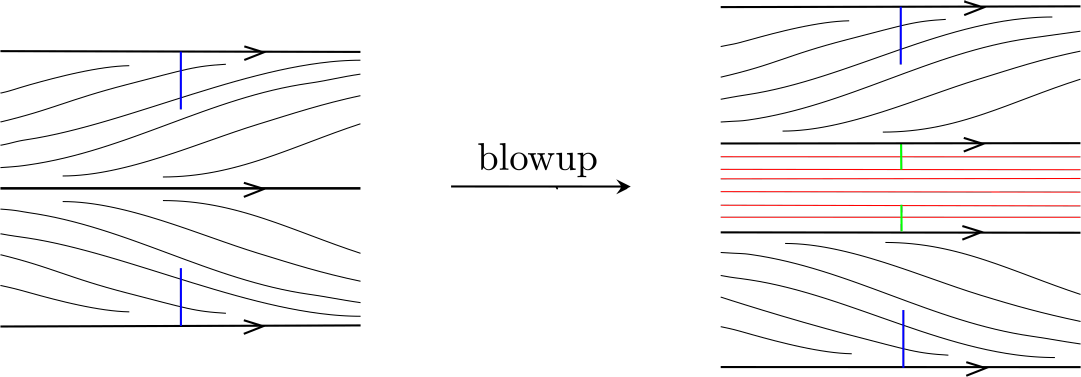}
\caption{The foliation on $T$ after blowup.}
\label{fig: blowup_leaves}
\end{figure}

\medskip 

\begin{figure}[ht]
\centering
\includegraphics[scale=0.5]{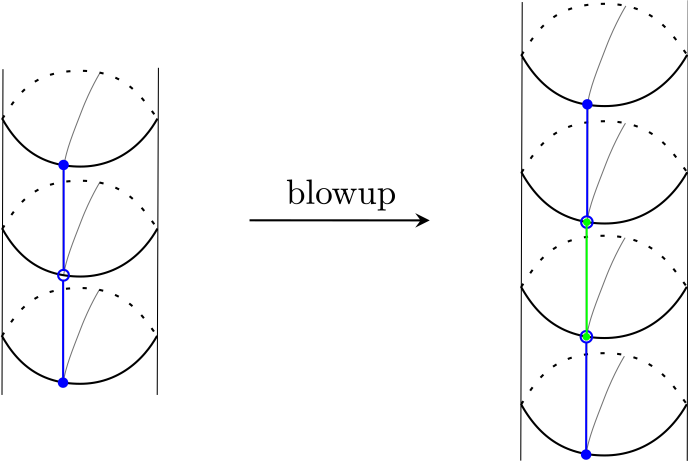}
\caption{After the blowup, the blue puncture is split into two, each of which is filled by the endpoint of a marker in green. }
\label{fig: after_blowup_boundary}
\end{figure}

 The local effect of this blowup on $T\cap \mathcal{F}$ is illustrated in Figure \ref{fig: blowup_leaves}. Figure \ref{fig: after_blowup_boundary} shows how the endpoints of the set of markers for the new foliation overlap.  Then the rest of the argument is the same as before.

\subsection{Proof of Theorem \ref{thm: universal circle fixed point}(2)}
\label{subsec: pf thm 5.1(2)} 
We continue to work under the assumptions of the statement of Theorem \ref{thm: universal circle fixed point} and to use the notation developed in the proof of Theorem \ref{thm: universal circle fixed point}(1). In particular, $F$ is a properly embedded surface in $M$ with connected boundary representing a generator of $H_2(M, \partial M) \cong \mathbb Z$.
We constructed a universal circle action $\rho: \pi_1(W) \rightarrow \mbox{{\rm Homeo}}_+(S^1)$ in \S \ref{subsec: existence uca} associated to a co-oriented taut foliation on the double $W$ of $M$ which has $S = F\cup_{\rm id} -F$ as a compact leaf. We assume that the restriction $\rho_M$ of $\rho$ to $\pi_1(M)$ lifts to $\widetilde \rho_M: \pi_1(M) \rightarrow \mbox{{\rm Homeo}}_{\mathbb Z}(\mathbb R)$ and show that the absolute value of the translation number of $\widetilde{\rho}_M(\lambda)$ is 
$2g - 1$, where $g = g(F)$ is the genus of $M$. 

By construction, the restriction of $\rho$ to $\pi_1(S)$ is semi-conjugate to a discrete faithful representation $\rho_S: \pi_1(S) \to PSL(2, \mathbb R)$, and since translation numbers are invariant under semi-conjugacy, it suffices to show that the absolute value of the translation number of $\widetilde{\rho}_F(\lambda)$ is $2g - 1$, where $\widetilde{\rho}_F: \pi_1(F) \to \widetilde{SL_2}$ is a lift of $\rho_F = \rho_S|_{\pi_1(F)}$.  

Write 
$$\pi_1(S) = \langle a_1, b_1, \ldots, a_{2g}, b_{2g} \; | \; [a_1, b_1] [a_2, b_2] \cdots [a_{2g}, b_{2g}] \rangle$$ 
and 
$$\pi_1(F) = \langle a_1, b_1, \ldots, a_{g}, b_{g} \; | \; [a_1, b_1] [a_2, b_2] \cdots [a_{g}, b_{g}] \rangle$$ 
and 
$$\pi_1(-F) = \langle a_{g+1}, b_{g+1}, \ldots, a_{2g}, b_{2g} \; | \; [a_{g+1}, b_{g+1}] [a_{g+2}, b_{g+2}] \cdots [a_{2g}, b_{2g}] \rangle$$
Since $\rho_S$ is discrete and faithful, if $A_i, B_i$ are arbitrary lifts of $\rho_S(a_i), \rho_S(b_i)$ to $\widetilde{SL_2}$,  then
$$[A_1, B_1] [A_2, B_2] \cdots [A_{2g}, B_{2g}] = \mbox{sh}(\pm (2g(S) - 2))$$
(See \cite[\S 2.2 and Theorem 2.8]{Mann}.) On the other hand, 
$$[a_{g+1}, b_{g+1}] [a_{g+2}, b_{g+2}] \cdots [a_{2g}, b_{2g}] = ([a_1, b_1] [a_2, b_2] \cdots [a_{g}, b_{g}])^{-1} \in \pi_1(S)$$
so $[A_{g+1}, B_{g+1}] [A_{g+2}, B_{g+2}] \cdots [A_{2g}, B_{2g}]$ commutes with $[A_1, B_1] [A_{2}, B_{2}] \cdots [A_{g}, B_{g}]$. Hence 
\begin{eqnarray} 
2(2g - 1) = 2g(S) - 2 & = & |\tau([A_1, B_1] [A_2, B_2] \cdots [A_{2g}, B_{2g}])| \nonumber \\ 
& = &  |\tau([A_1, B_1] \cdots [A_{g}, B_{g}]) + \tau([A_{g+1}, B_{g+1}] \cdots [A_{2g}, B_{2g}])|  \nonumber \\
& \leq&  |\tau([A_1, B_1]  \cdots [A_{g}, B_{g}])| + |\tau([A_{g+1}, B_{g+1}]\cdots [A_{2g}, B_{2g}])|  \nonumber  
\end{eqnarray} 
As products of $g$ commutators both $|\tau([A_1, B_1]  \cdots [A_{g}, B_{g}])|$ and $|\tau([A_{g+1}, B_{g+1}]\cdots [A_{2g}, B_{2g}])|$ are bounded above by 
$2g - 1$ (see e.g. \cite[\S 2]{JN}). Hence $|\tau([A_1, B_1]  \cdots [A_{g}, B_{g}])| = 2g - 1$. Finally, taking $A_i = \tilde \rho_M(a_i)$ and $B_i = \tilde \rho_M(b_i)$ yields
$$|\tau(\tilde \rho_F(\lambda))| = |\tau([A_1, B_1]  \cdots [A_{g}, B_{g}])| = 2g - 1,$$
which completes the proof.

\begin{proof}[Proof of Theorem \ref{thm: meridional detn}(2)]
Theorem \ref{thm: meridional detn}(2) will follow from Proposition \ref{prop: divide tau lambda} and Theorem \ref{thm: universal circle fixed point}(2) once we verify that the representation $\rho_M$ constructed in Theorem \ref{thm: universal circle fixed point}(1) lifts to $\mbox{{\rm Homeo}}_{\mathbb Z}(\mathbb R)$. The obstruction to lifting $\rho_M$ is its Euler class $e(\rho_M) \in H^2(M) \cong T_1(M)$, which coincides with the Euler class of the tangent bundle of the foliation $\mathcal{F}_M$ of \S \ref{subsec: existence uca}. The latter was shown to vanish whenever $H^2(M)$ is a $\mathbb Z/2$ vector space in \cite[Proposition 2.1]{Hu2}.  
\end{proof}

\subsection{Families of non-order-detected slopes} 
Corollary \ref{cor: ord detd fibration} states that if $K$ is a non-trivial knot in the $3$-sphere,  each rational slope whose distance from the longitude of $K$ divides $2g(K) - 1$ is order-detected in the exterior $M$ of $K$. In particular this holds for the slopes $2g(K) - 1$. Using the work of Nie (\cite{Nie}) we prove a complement to this result by producing positive $L$-space knots $K$ for which no  rational slope $m/n > 2g(K) - 1$ is order-detected. We take $\mu, \lambda \in \pi_1(\partial M)$ to be a meridian, longitude pair for $K$. 

The root-closed, conjugacy-closed submonoid of $\pi_1(M)$ generated by $S \subseteq \pi_1(M)$ is the minimal subset  $\mathcal{M}(S)$ of $\pi_1(M)$ containing $S$ which has the following properties:
\vspace{-.2cm}
\begin{itemize}

\item $\eta, \gamma \in \mathcal{M}(S) \Rightarrow \eta \gamma \in \mathcal{M}(S)$;

\vspace{.2cm} \item $\gamma \in \mathcal{M}(S)$ and $\eta \in \pi_1(M) \Rightarrow \eta \gamma \eta^{-1} \in \mathcal{M}(S)$; 

\vspace{.2cm} \item $\gamma^n \in \mathcal{M}(S)$ for some integer $n > 0 \Rightarrow \gamma \in \mathcal{M}(S)$.

\end{itemize}
(cf. the proof of Proposition \ref{prop: hfgst and detection}). Following Nie (\cite{Nie}) we say that $K$ has property (D) if
\vspace{-.2cm}
\begin{itemize}

\item[(a)]  for any homomorphism $\rho: \pi_1(M) \to \mbox{Homeo}_+(\mathbb R)$, if $s \in \mathbb R$ is a common fixed point of $\rho(\mu)$ and $\rho(\lambda)$, then $s$ is a fixed point of every element in $\pi_1(M)$\footnote{Condition (a) is equivalent to requiring that each $\mathfrak{o} \in LO(M)$ be boundary-cofinal; see Remark \ref{rem: bdry cof iff fpf}.}, and,

\vspace{.2cm} \item[(b)] $\mu$ is in the root-closed, conjugacy-closed submonoid $\mathcal{M}(\mu^{2g(K)-1} \lambda, \mu^{-1})$ generated by $\mu^{2g(K)-1} \lambda$ and $\mu^{-1}$.

\end{itemize}
Nie showed that manifolds obtained by $p/q$-surgery on a knot $K$ satisfying property (D) have non-left-orderable fundamental groups when $p/q \geq 2g(K) - 1$. We show, 

\begin{proposition}  
\label{prop: (d) implies not ord-det}
Suppose that $K$ is a knot in the $3$-sphere which satisfies property $(D)$. Then no rational slope $p/q > 2g(K)-1$ is weakly order-detected.
\end{proposition}

\begin{proof} 
Suppose that $\mathfrak{o} \in LO(M)$ weakly detects the slope $[p\mu + q\lambda]$, where $p, q > 0$ are integers and $p/q > 2g(K)-1$. 
Hence if $\rho_\mathfrak{o}: \pi_1(M) \to \mbox{Homeo}_+(\mathbb R)$ is a dynamic realisation of $\mathfrak{o}$,  
condition (a) of the definition of property (D) shows that  $\rho_\mathfrak{o}(\pi_1(\partial M))$ has no fixed points in $\mathbb R$. Then as $\pi_1(\partial M)$ is abelian,  the image by $\rho_\mathfrak{o}$ of one of its primitive elements has no fixed points in $\mathbb R$, and so is conjugate to $\mbox{sh}(\pm 1)$. Hence, up to 
replacing $\rho_\mathfrak{o}$ by a conjugate representation we can suppose that $\rho_\mathfrak{o}(\pi_1(\partial M)) \leq \mbox{Homeo}_{\mathbb 
Z}(\mathbb R)$. Since we have assumed that $\mathfrak{o} \in LO(M)$ detects $[p\mu + q\lambda]$, the composition $\tau \circ (\rho_\mathfrak{o}|_{\pi_1(\partial M)})$ is a non-trivial homomorphism $\pi_1(\partial M) \to \mathbb R$ whose kernel is generated by $p\mu + q\lambda$ (Lemma \ref{lemma: slope by tau}). In particular, $\rho_\mathfrak{o}(\mu)$ acts without fixed point on $\mathbb R$, so we can assume that $\rho_\mathfrak{o}(\mu) = \mbox{sh}(\pm 1)$.

First suppose that $\rho_\mathfrak{o}(\mu) = \mbox{sh}(1)$, so that $\tau(\rho_\mathfrak{o}(\mu)) = 1$. 
Then 
\begin{eqnarray} 
0 = \tau(\rho_\mathfrak{o}(\mu^p\lambda^q)) & = & (p - (2g(K)-1)q)\tau(\rho_\mathfrak{o}(\mu)) + q\tau(\rho_\mathfrak{o}(\mu^{2g(K)-1}\lambda)) \nonumber \\ 
& = & (p - (2g(K)-1)q) + q\tau(\rho_\mathfrak{o}(\mu^{2g(K)-1}\lambda)) 
\nonumber \\ 
& > & q\tau(\rho_\mathfrak{o}(\mu^{2g(K)-1}\lambda)) \nonumber   
\end{eqnarray} 
Thus $\tau(\rho_\mathfrak{o}(\mu^{2g(K)-1}\lambda)) < 0$ and therefore $\rho_\mathfrak{o}(\mu^{2g(K)-1}\lambda)$ is a strictly decreasing homeomorphism. This is also true of $\rho_\mathfrak{o}(\mu^{-1}) = \mbox{sh}(-1)$ and so as $\mu$ is in the root-closed, conjugacy-closed submonoid $\mathcal{M}(\mu^{2g(K)-1} \lambda, \mu^{-1})$, $\mbox{sh}(1) = \rho_\mathfrak{o}(\mu)$ is strictly decreasing, which is false. Thus $\rho_\mathfrak{o}(\mu)\ne \mbox{sh}(1)$.

A similar argument shows that if $\rho_\mathfrak{o}(\mu) = \mbox{sh}(-1)$, then $\rho_\mathfrak{o}(\mu^{-1})$ and $\rho_\mathfrak{o}(\mu^{2g(K)-1}\lambda)$ are strictly increasing which, as $\mu \in \mathcal{M}(\mu^{2g(K)-1} \lambda, \mu^{-1})$, leads to the impossible conclusion that $\rho_\mathfrak{o}(\mu) = \mbox{sh}(-1)$ is as well. Thus no $\mathfrak{o} \in LO(M)$ detects the slope $[p\mu + q\lambda]$, when $p/q > 2g(K)-1$.
\end{proof}

Nie has shown that many $L$-space knots satisfy property (D) including nontrivial $(1,1)$-knots which are positive $L$-space knots (\cite{Nie}). This family includes $1$-bridge braid knots and certain families of twisted torus knots.

\begin{corollary}  
\label{cor: nie 1bb}
If $K$ is a nontrivial $(1,1)$-knot which is a positive $L$-space knot, then no  rational slope $p/q > 2g(K) - 1$ on the boundary of the exterior of $K$ is order-detected. 
 
\end{corollary}

\section{Foliation-detection in fibred knot manifolds}
\label{sec: ctf on hyperbolic fibred knot}

The goal of this section is to prove Theorem \ref{thm: meridional detn}(3). In this section we use $X$ to denote a fibred rational homology solid torus. We prove that 
each  rational slope on $\partial X$ that is distance $1$ from the longitudinal slope is foliation-detected. 
 
The fibre $F$ of $X$ is a connected compact surface satisfying $\partial F\cong S^1$ and $\chi(F)<0$. The main step in the proof is to extend a lamination on $\partial X$ that detects the desired slope into $X$. This is done by applying a mild modification of Tao Li's argument in the proof of \cite[Theorem 1]{Li1}  to the branched surfaces constructed in \cite{Rob2}.  

\subsection{Branched surfaces}
\label{subsubsec:branched surface}
A {\it branched surface} $B$ is a compact subspace of a compact $3$-manifold $M$ locally modeled on the spaces depicted in Figure \ref{fig:branched surfaces} \cite{FO}. The set of points $\Gamma$ of $B$ which do not have $2$-disk neighbourhoods is called the {\it branch 
locus} of $B$ and the components of its complement $B\setminus \Gamma$ are called the {\it sectors} of $B$.  In Figure \ref{fig:branched surfaces}, the arrow  crossing an edge  of a  sector points in the direction in which two sectors merge into one.  A {\it sink disk} of $B$ is a disk sector with all arrows pointing into it. A sink disk that intersects $\partial 
M$ nontrivially is called a {\it half sink disk}. We call a branched surface {\it sink disk free} if it contains no sink disk or half sink disk.  The condition that $B$ be sink disk free is crucial in Tao Li's construction in \cite{Li1, Li2}.

\begin{figure}[ht]
\begin{tikzpicture}[scale=0.65]
\draw (1,1) -- (4,1) -- (5,3) -- (2,3) -- (1,1);
\draw (7,1) -- (10,1) -- (11,3) -- (8,3) -- (7,1);
\draw (8.5,1) to [out=0, in=235] (9.5, 1.5) to  (10.5,3.5);
\draw (8.5,1) -- (9.5, 3);
\draw (9.5, 3) to [out=0, in=235] (10.5,3.5);
\draw [-{>[scale=2.5, length=1, width=1]}] (9.2, 2) -- (8.7, 2);
\draw (13,1) -- (16,1) -- (17,3) -- (14,3) -- (13,1);
\draw (14.5,1) to [out=0, in=235] (15.5, 1.5) to  (16.5,3.5);
\draw (14.5,1) -- (15.5, 3);
\draw (15.5, 3) to [out=0, in=235] (16.5,3.5);
\draw [-{>[scale=2.5, length=1, width=1]}] (14.9, 1.5) -- (14.5, 1.5);
\draw [-{>[scale=2.5, length=1, width=1]}] (14.3, 1.85) -- (14.3, 2.2);
\draw [dashed] (13.5,2) -- (16.5,2);
\draw [dashed] (13.4,1.15) -- (16.05,1.15);
\draw (16.05,1.15) -- (16.4,1.15);
\draw (16.5,2) to [out = 250, in=85] (16.4,1.15);
\draw [dashed] (13.5,2) to [out = 250, in=85] (13.3,1.15);
\end{tikzpicture}
 \caption{Local models of branched surfaces}
 \label{fig:branched surfaces}
\end{figure} 

A branched surface $B$ properly embedded in a $3$-manifold $M$ has a regular neighbourhood $N(B)$ in $M$ fibred by intervals $I$.  (See Figure \ref{fig:N of branched surface}.) The {\it vertical boundary of $N(B)$}, denoted $\partial_v N(B)$, is that part of $\partial N(B)$ which is tangent to the fibres. The {\it horizontal boundary of $N(B)$}, denoted $\partial_h N(B)$, is the part that is transverse to the fibres.  

\begin{figure}[ht]
\begin{tikzpicture}[scale=0.65]
\draw (1,1) -- (4,1) -- (5,3) -- (2,3) -- (1,1);
\draw (1, 0.8) -- (4, 0.8);
\draw (4, 0.8) -- (5,2.8);
\foreach \x in {1, 1.1, 1.2, ..., 4}
	\draw (\x, 0.8) -- (\x, 1);
\foreach \n in {1,2,...,10}
	\draw (4+0.1*\n, 0.8+0.2*\n) -- (4+0.1*\n, 1+0.2*\n);
\draw (10,1) -- (11,3);
\draw (8,3) -- (7,1);
\draw (7,1) -- (8.5,1);
\draw (8,3) -- (9.5,3);
\draw (9.3,1) -- (10,1);
\draw (8.5,1) to [out=0, in=235] (9.5, 1.5) to  (10.5,3.5);
\draw (9.5, 3) to [out=0, in=235] (10.5,3.5);
\draw (7, 0.8) -- (10, 0.8);
\draw (10, 0.8) -- (11,2.8);
\foreach \x in {1, 1.1, 1.2, ..., 4}
	\draw (6+\x, 0.8) -- (6+\x, 1);
\foreach \n in {1,2,...,10}
	\draw (10+0.1*\n, 0.8+0.2*\n) -- (10+0.1*\n, 1+0.2*\n);
\foreach \n in {0,1,2,...,7,8}
	\draw (8.5+0.1*\n, 1) -- (8.5+0.1*\n, 1+0.00055*\n*\n*\n);
\draw (9.3, 1.1) to [out=25, in=225] (9.55, 1.3);
\draw (9.45, 1.18) -- (9.41, 1.36);
\draw (10.4, 3) -- (11,3);
\draw (9.55, 1.3) --  (10.55,3.3);
\foreach \n in {0,1,2,...,10}
	\draw (9.55 + 0.1*\n , 1.3+ 0.2*\n) -- (9.5+0.1*\n , 1.5+ 0.2*\n);
\draw (16,1) -- (17,3);
\draw (14,3) -- (13,1);
\draw (13,1) -- (14.5,1);
\draw (14,3) -- (15.5,3);
\draw (15.3,1) -- (16,1);
\draw (14.5,1) to [out=0, in=235] (15.5, 1.5) to  (16.5,3.5);
\draw (15.5, 3) to [out=0, in=235] (16.5,3.5);
\foreach \n in {0,1,2,...,7,8}
	\draw (14.5+0.1*\n, 1) -- (14.5+0.1*\n, 1+0.00055*\n*\n*\n);
\draw (15.3, 1.1) to [out=25, in=225] (15.55, 1.3);
\draw (15.45, 1.18) -- (15.41, 1.36);
\draw (16.4, 3) -- (17,3);
\draw (15.55, 1.3) --  (16.55,3.3);
\foreach \n in {0,1,2,...,10}
	\draw (15.55 + 0.1*\n , 1.3+ 0.2*\n) -- (15.5+0.1*\n , 1.5+ 0.2*\n);
\draw (13, 0.8) -- (16, 0.8);
\draw (16, 0.8) -- (16.5, 1.8);
\foreach \x in {1, 1.1, 1.2, ..., 4}
	\draw (12+\x, 0.8) -- (12+\x, 1);
\draw (17,2.4) -- (16.7, 1.7) to [out=245, in=90] (16.5,0.9);
\foreach \n in {0, 1, 2, ...,10}
	\draw (16+ 0.1*\n, 1+0.2*\n) to (16+ 0.1*\n, 0.8+0.2*\n);
\foreach \n in {0, 1, 2, ...,5}
	\draw (16.5+ 0.1*\n, 1.8+0.2*\n) to [out = 270, in = 135-9*\n] (16.6+0.08*\n, 1.5+0.18*\n);
\draw (16.38, 1.56) -- (16.55, 1.56);
\draw (16.05,0.9) -- (16.5, 0.9);
\draw (16.55, 1.56) to [out= 235, in =85] (16.38,1.05) -- (16.12, 1.05);
\draw (16.38, 1.05) to (16.5,0.9);
\draw (16.25, 1.05) -- (16.25, 0.9);
\draw (16.17, 1.05) -- (16.15, 0.9);
\draw (16.33, 1.05) -- (16.35, 0.9);
\draw (16.39, 1.16) -- (16.51, 1.06);
\draw (16.41, 1.28) -- (16.52, 1.19);
\draw (16.45, 1.37) -- (16.55, 1.29);
\draw (16.50, 1.49) -- (16.58, 1.4);
\end{tikzpicture}
 \caption{Local models of normal neighbourhood of branched surfaces}
 \label{fig:N of branched surface}
\end{figure} 

The branched surfaces we consider are all transversely oriented. That is, 
the $I$-fibres of $N(B)$ can be oriented coherently. Consequently, the sectors of $B$ inherit an orientation in which each interval intersects $B$ 
positively in $M$.

A lamination of a $3$-manifold $M$ is a co-dimension one foliation on a closed subset of $M$. We say a lamination $\mathcal{L}$ is {\it carried} by a branched surface $B$ if $\mathcal{L}$ can be isotoped to lie in $N(B)$ transverse to the fibres; and $\mathcal{L}$ is {\it fully carried} by $B$ if it can be isotoped to lie in $N(B)$ so that it intersects each fibre transversely in a non-empty set. 

\subsection{Extending a lamination on $\partial M$ into $M$}
\label{subsubsec:extend lamination}
Let $B$ be a branched surface properly embedded in a $3$-manifold $M$ with torus boundary and sink disk free. 

The intersection $\partial B = B\cap \partial M$ is a train track embedded in $\partial M$. Let $N(\partial B) =N(B)\cap \partial M $ denote the normal neighbourhood of $\partial B$ in $\partial M$. Suppose that $\tau$ is a lamination on $N(\partial B)$ that is fully carried by $\partial B$. We assume below that the intersection of $\tau$ with each $I$-fibre of $N(\partial B)$ is the standard Cantor set $\mathcal{C}$. 

Next we describe Tao Li's ideas in \cite{Li1}, which we will use to extend $\tau$ to a lamination on $N(B)$ fully carried by $B$. Of key importance are a graph $L'$ in $B$ and a branched subsurface $P(L')$ of $B$ which are carefully defined at the beginning of \cite[\S 4]{Li1}. Our description will be  largely pictorial. For the reader's convenience, we will use the same notation as in \cite[\S 4]{Li1}. We also refer the reader to \S 
\ref{subsubsec: figure eight} below, where the procedure is implemented to construct the desired lamination in the exterior of the figure eight knot.  

\begin{figure}[ht]
 \centering 
 \includegraphics[scale=0.16]{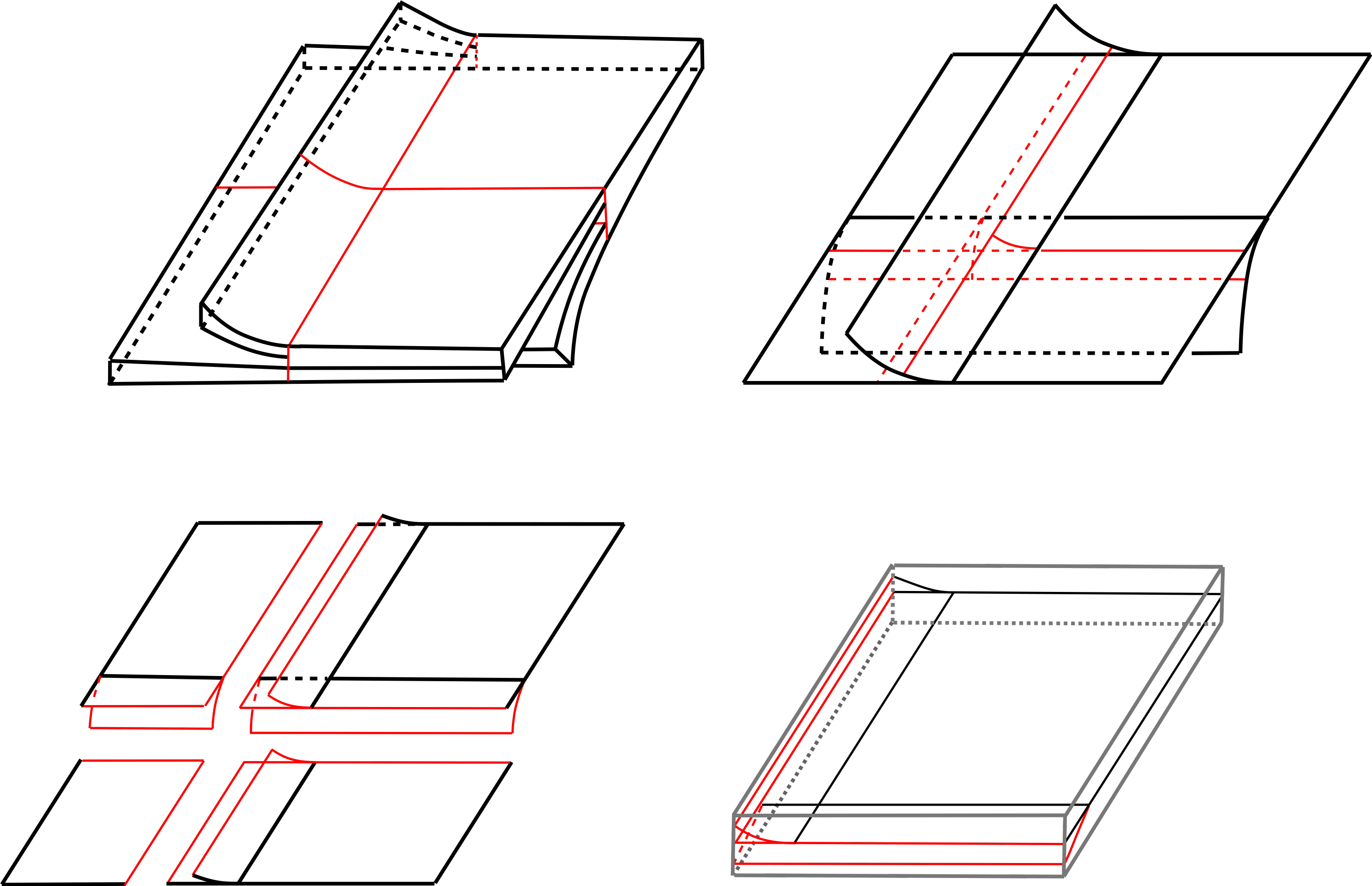}
 \caption{Decompositions of $B$ and $N(B)$}
  \label{fig:decompose branched surface}
\end{figure}

We begin by cutting open the normal neighbourhood $N(B)$ of $B$ along the $I$-fibres which intersect $\overline{\partial_v N(B) \setminus \partial M}$ (depicted by the red lines in Figure \ref{fig:decompose branched surface}). This simultaneously decomposes the branched surface $B$ in such a way that the components of the decomposed $N(B)$ correspond bijectively to the components of the decomposed $B$. Since those parts of $\partial_v N(B)$ contained in the interior of $M$ are ``parallel'' to the branch locus $\Gamma$ of $B$, these components correspond naturally to the sectors of $B$. Given such a sector $\mathcal{D}$, we use $\mathcal{D}^B$ and $N_B(\mathcal{D})$ to denote the corresponding 
pieces of $B$ and $N(B)$ respectively. By construction, the branched surface $\mathcal{D}^B$ is properly embedded inside $N_B(\mathcal{D})$ (see Figure \ref{fig:decompose branched surface}). 

The construction of the lamination on $M$ now proceeds in two parts: 
\vspace{-.2cm}
\begin{enumerate}
 \item For each sector $\mathcal{D}$, construct a lamination on $\partial_v N_B(\mathcal{D})$ which extends $\tau\cap \partial_v N_B(\mathcal{D})$ 
so that if $\mathcal{D}'$ is another sector, the extended laminations match on $\partial_v N_B(\mathcal{D}) \cap \partial_v N_B(\mathcal{D}')$.

\vspace{.2cm} \item For each sector $\mathcal{D}$, extend the lamination on $\partial_v N_B(\mathcal{D})$ into the interior of $N_B(\mathcal{D})$.

\end{enumerate}

\begin{figure}[ht]
\centering 
\includegraphics[scale=0.47]{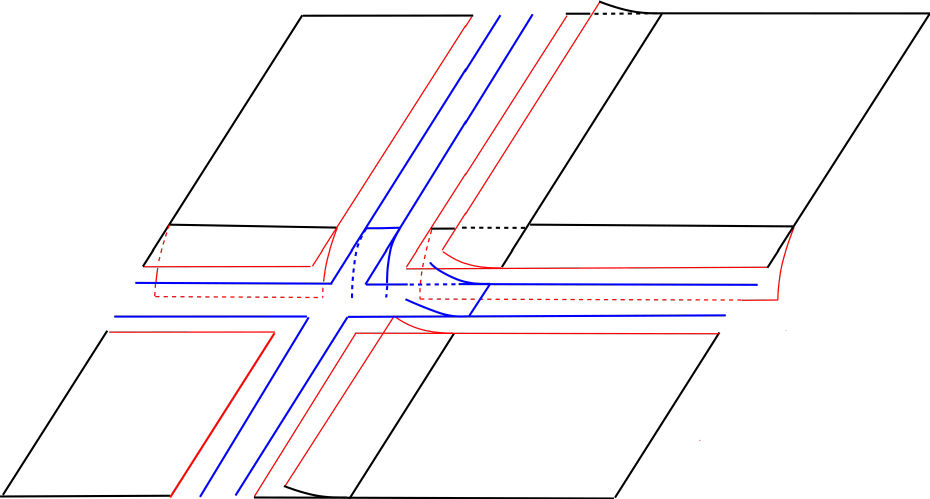}
\caption{$P(L')$ is coloured in blue.}
\label{fig:PL}
\end{figure}

Part (1) is relatively straightforward. Let $L'$ denote the graph on $B$ consisting of the red arcs shown in the picture in the top-right corner of Figure \ref{fig:decompose branched surface} and let $N(L')$ be an appropriate regular neighbourhood of $L'$ in $M$. Then $P(L') = N(L')\cap B$ is a branched surface embedded in $B$. Note that the branch locus of $P(L')$ is a disjoint union of arcs (see Figure \ref{fig:PL} and \cite[Figure 
4.1]{Li1}). 

By a slight abuse of the notation $\mathcal{D}^B$, we shall view $B$ as the union of the $\mathcal{D}^B$ and the branched surface $P(L')$. This is 
illustrated in Figure \ref{fig:PL}.  Similarly,  $N(B)$ can be viewed as the union of the $N_B(\mathcal{D})$ and the normal neighbourhood  $N(P(L'))$ of $P(L')$. 

Endow the normal neighbourhood of each sector of $P(L')$ with the lamination obtained from the product of the sector with the Cantor set $\mathcal{C}$. Since $\mathcal{C} \cong \mathcal{C} \sqcup \mathcal{C}$ and the branch locus of $P(L')$ has no double point, we can merge these laminations 
along neighbouring sectors of $P(L')$ to obtain a lamination on $N(P(L'))$ (cf. Step 1 in \S \ref{subsubsec: figure eight}). 

On $\partial M$, we can identify the part of the lamination on $\partial_v N(P(L')) \cap \partial M$ with the corresponding portion of $\tau$ using a homeomorphism of Cantor sets. This ensures that the resulting lamination carried by $P(L')$ extends $\tau \cap N(P(L'))$. 

Now fix a sector $\mathcal{D}$ of $B$ and note that $\partial_v N_B(\mathcal{D}) = (\partial_v N_B(\mathcal{D}) \cap \partial_v N(P(L'))) \cup (\partial_v N_B(\mathcal{D}) \cap \partial M)$. Define a lamination on $\partial_v N_B(\mathcal{D})$ to be that induced by the lamination on $\partial_v N_B(\mathcal{D}) \cap \partial_v N(P(L'))$ and that induced by $\tau$ on $\partial_v N_B(\mathcal{D}) \cap \partial M$. This completes Part (1). 

Part (2) of the construction is more delicate and the precise details, which depend on the branched surface in question, will occupy the rest of this section.

\subsection{Branched surfaces in a fibred knot exterior}
\label{subsubsec:branched surface fibred knot}
Let $X$ be an oriented fibred rational homology solid torus
$$X = F \times [0,1]/ (x,1)\sim (h(x), 0)$$
where $h$ is the monodromy whose restriction to $\partial F$ is the identity.

Let $\alpha: [0,1]\rightarrow F$ be an oriented properly embedded non-separating arc on $F$.  We freely isotope $h(\alpha)$ to an arc $\alpha_-$ so that  $\alpha$  and $\alpha_-$ have minimal intersection, possibly empty. Since $b_1(X)=1$, we have that $\alpha$ and $\alpha_-$ are not isotopic. Otherwise, there exists a non-separating surface in $X$, obtained essentially by taking the quotient of $\alpha \times [0,1]$. 
Together with the fibre surface $F$ in $X$, this leads to a contradiction to the assumption that $H^1(X) \cong H_2(X,\partial X)$ has rank $1$.

Let $D=\alpha\times [0,1]/\sim$, a (singular) disk in $X$. After isotopy, we assume that the top side of $D$ is attached to the negative side of $F$ along $\alpha_-$ (instead of $h(\alpha)$). We fix an orientation on $D$ and then tilt the disk slightly so that $TD = TF$ at their intersection $D\cap 
F $, so the union of $F$ and $D$ defines a branched surface embedded in $X$ whose branch locus is $D\cap F = \alpha \cup \alpha_-$. We denote it 
by $B = \langle F; D\rangle$. Figure \ref{fig: figure eight train track} shows an example of $\alpha$ and $\alpha_-$ on a fibre of the exterior of the figure-eight knot. 

\subsubsection{Types of   train  tracks on $\partial X$}
\label{subsubsec:boundary train track}
After possibly reversing the orientation of $X$ and replacing $h$ by $h^{-1}$, we can orient the product disk $D$ so that the boundary train track 
$B\cap \partial X$ belongs to one of the three types illustrated in Figure \ref{fig:types of boundary train track}, which correspond to the ``negatively oriented'' train tracks in \cite[Figure 4]{Rob2}.  We use $\mu_0$ to denote the meridional slope of $X$ that is disjoint from $D\cap \partial X$, which is called the canonical meridian in 
\cite[\S 3]{Rob2}. 

\medskip
\begin{figure}[ht]
\centering
\begin{tikzpicture}[scale=0.65]
\draw [white] (11, 10) -- (11.2,10);
\draw [thick, gray] (2.5, 12.5) to [bend left] (2.5, 14.5);
\draw [thick, gray, dashed] (2.5, 12.5) to [bend right] (2.5, 14.5);
\draw [thick, gray] (8, 12.5) to [bend left] (8, 14.5);
\draw [thick, gray] (8, 12.5) to [bend right] (8, 14.5);
\draw [thick, gray] (2.5, 12.5) -- (8,12.5);
\draw [thick, gray] (2.5, 14.5) -- (8,14.5);
\draw [thick] (2.2, 13.5) -- (7.7, 13.5);
\draw [thick, {<[scale=2.5, length=1, width=2]}-] (3, 13.5) -- (3.2, 13.5);
\draw [thick, gray, -{>[scale=2.5, length=1, width=2]}] (7.703,13.7) to  (7.706,13.73);
\node [gray] at (8, 13.5) {\tiny{$\mu_0$}};
\node [thick] at (0,13.5) {Type 1};
\draw [thick] (7.2, 13.5) .. controls (6.7,13.5) and (7,14.5) .. (6.5, 14.5);
\draw [thick, dashed] (6.5, 14.5) .. controls (6,14.5) and (6.5,12.5) .. (6,12.5);
\draw [thick] (6,12.5) .. controls (5.5,12.5) and (6,13.5) .. (5.3, 13.5);
\draw [thick] (4.5,13.5) .. controls (5,13.5) and (5,14.5) .. (4.5, 14.5);
\draw [thick,dashed] (4.5, 14.5) .. controls (3.8,14.5) and (4.5,12.5) .. 
(3.5,12.5);
\draw [thick] (3.5,12.5) .. controls (3,12.5) and (3.2,13.5) .. (3.7,13.5);

\draw [thick, gray] (2.5, 9) to [bend left] (2.5, 11);
\draw [thick, gray, dashed] (2.5, 9) to [bend right] (2.5, 11);
\draw [thick, gray] (8, 9) to [bend left] (8, 11);
\draw [thick, gray] (8, 9) to [bend right] (8, 11);
\draw [thick, gray] (2.5, 9) -- (8,9);
\draw [thick, gray] (2.5, 11) -- (8,11);
\draw [thick] (2.2, 10) -- (7.7, 10);
\draw [thick, {<[scale=2.5, length=1, width=2]}-] (3, 10) -- (3.2, 10);
\draw [thick, gray, -{>[scale=2.5, length=1, width=2]}] (7.703,10.2) to  (7.706,10.23);
\node [gray] at (8, 10) {\tiny{$\mu_0$}};
\node [thick] at (0,10) {Type 2(a)};
\draw [thick] (7.2, 10) .. controls (6.7,10) and (7,11) .. (6.5, 11);
\draw [thick, dashed] (6.5, 11) .. controls (6,11) and (6,9) .. (5.3,9);
\draw [thick] (5.3, 10) .. controls (5.8,10) and (5.8,11) .. (5.3, 11);
\draw [thick] (5.3,9) .. controls (4.5,9) and (5,10) .. (4.3,10);
\draw [thick,dashed] (3.5,9) .. controls  (4.5,9) and (4.8, 11) ..(5.3, 11) ;
\draw [thick] (3.5,9) .. controls (3,9) and (3.2,10) .. (3.7,10);
\draw [thick, gray] (2.5, 5.5) to [bend left] (2.5, 7.5);
\draw [thick, gray, dashed] (2.5, 5.5) to [bend right] (2.5, 7.5);
\draw [thick, gray] (8, 5.5) to [bend left] (8, 7.5);
\draw [thick, gray] (8, 5.5) to [bend right] (8, 7.5);
\draw [thick, gray] (2.5, 5.5) -- (8,5.5);
\draw [thick, gray] (2.5, 7.5) -- (8,7.5);
\draw [thick] (2.2, 6.5) -- (7.7, 6.5);
\draw [thick, {<[scale=2.5, length=1, width=2]}-] (3, 6.5) -- (3.2, 
6.5);
\draw [thick, gray, -{>[scale=2.5, length=1, width=2]}] (7.703,6.7) 
to  (7.706,6.73);
\node [gray] at (8, 6.5) {\tiny{$\mu_0$}};
\node [thick] at (0,6.5) {Type 2(b)};
\draw [thick] (4.3,6.5) .. controls (4.8,6.5) and (5,7.5) .. (4.4, 7.5);
\draw [thick,dashed] (4.3, 7.5) .. controls (3.6,7.5) and (4.5,5.5) .. (3.5,5.5);
\draw [thick] (3.5,5.5) .. controls (3,5.5) and (3.2,6.5) .. (3.7,6.5);
\draw [thick] (6.2,6.5) .. controls (5.7,6.5) and (5.5,7.5) .. (6.1, 7.5);
\draw [thick,dashed] (6.2, 7.5) .. controls (6.9,7.5) and (6,5.5) .. (7,5.5);
\draw [thick] (7,5.5) .. controls (7.5,5.5) and (7.3,6.5) .. (6.8,6.5);
\end{tikzpicture}
\caption{Three possible configurations of the boundary train track $B \cap \partial X$.}
\label{fig:types of boundary train track}
\end{figure}


With respect to the canonical meridian $\mu_0$, in \cite{Rob2}, Roberts showed that there exists a taut foliation $\mathcal{F}$ on $X$ intersecting $\partial X$ transversely in simple closed curves of slope $p/q$ for any $p/q$ in $(-\infty, 1)$ if $B\cap \partial X$ is of Type 1 or Type 2(a), and for any $p/q$ in $(-\infty, \infty)$ if $B\cap \partial X$ is of Type 2(b) (\cite[Theorem 4.1, Corollary 4.3, Corollary 4.4]{Rob2}).  Hence, 
to complete the proof of Theorem \ref{thm: meridional detn}(3), it remains to show that $\mu_0$ is foliation-detected and when $B\cap \partial X$ is of Type 1 or Type 2(a), the slope $1/1$ is also foliation-detected.

In \S \ref{subsubsec: figure eight}, we first use an example to illustrate the idea of the proof. The full proof splits into three parts depending on the type of the train track which arises and are contained in \S \ref{subsec: type 1}, \S \ref{subsec: type 2a} and \S \ref{subsec: type 2b} respectively.

\subsection{An example}
\label{subsubsec: figure eight}
Let $X$ be the mapping torus of $h = T_a T_b^{-1}: F\rightarrow F$, where $F$ is a genus one surface with one boundary component, and $T_a$ and $T_b$ are the right-handed Dehn twists along curves $a$ and $b$ depicted in Figure \ref{fig: figure eight train track}. So $X$ is the exterior of the figure eight knot.

{\underline{Step 0: Understanding the branched surface $B =\langle F; D 
\rangle$}}

As described in \S \ref{subsubsec:branched surface fibred knot}, the branched surface $B =\langle F; D \rangle$ is obtained by attaching the product disk $D$ to the positive side of $F$ along $\alpha$ and to the negative side of $F$ along $\alpha_-$, where $\alpha_-$ is isotopic to $h(\alpha)$. In Figure \ref{fig: figure eight train track}, the arrows crossing $\alpha$ and $\alpha_-$ indicate the branch direction along each of them. 
Following $\partial F$, one can verify that the train track $B\cap \partial X$ is of Type 1. 

The branched surface $B$ has two disk sectors, as shown in  Figure \ref{fig: figure eight sectors}. One is the product disk $D$ and the other, $F'$, is obtained by cutting $F$ open along $\alpha \cup \alpha_-$. 

\begin{figure}[ht]
\centering
\begin{tikzpicture}[scale=0.8]
\draw [gray!80, dashed] (3.5, 2.1) to [bend left] (3.55, 3.3);
\draw [gray!80] (3.5, 2.1) to [bend right] (3.55, 3.3);
\node [gray!80] at (3.8, 2.6) {\small $b$};
\draw [gray!80] (3.45,3.5) ellipse (0.85 and 0.35);
\node [gray!80] at (4,4) {\small $a$};
\draw [thick, org] (1.67, 3.1) ..controls (2, 3.1) and (2.6,3.1)  .. (2.8, 3.5);
\draw [thick, org, dashed] (1.3, 3.83) ..controls (2.1, 3.85) and (2.6,3.9)  .. (2.8, 3.5);
\draw [thick, org, -{>[scale=2.5, length=0.7, width=1]}] (2.5, 3.1) 
-- (2.4, 3.35);
\node [org, above] at (2,3.1) {$\alpha$};
\draw [thick, blue] (1.55, 2.81) .. controls (3,2.5) and (6.3, 3.6) .. (3.55, 3.69);
\draw [thick, blue, dashed] (1.38, 4.1) .. controls (2, 4.1) and (3, 4.2) 
.. (3.5, 3.7);
\node [blue, below] at (4.8,3.4) {$\alpha_-$};
\draw [thick, blue, -{>[scale=2.5, length=0.7, width=1]}] (2.48, 2.88) -- (2.6, 2.6);
\draw [thick] (1.5,2.7) to [bend left] (1.5, 4.3);
\draw [thick] (1.5,2.7) to [bend right] (1.5, 4.3);
\draw [thick] (1.5, 2.7) .. controls (2.5,2) and (5.5, 1.5) .. (5.5, 3.5);
\draw [thick] (1.5, 4.3) .. controls (2.5,5) and (5.5, 5.5) .. (5.5, 3.5);
\draw [thick] (2.8, 3.5) to [bend left] (4.1, 3.5);
\draw [thick] (2.8, 3.5) to [bend right] (4.1, 3.5);
\draw [thick,  -{>[scale=2.5, length=1.5, width=1.2]}] (1.735,3.5) -- (1.735,3.45);
\draw [thick, gray] (7.5, 2.5) to [bend left] (7.5, 4.5);
\draw [thick, gray, dashed] (7.5, 2.5) to [bend right] (7.5, 4.5);
\draw [thick, gray] (13, 2.5) to [bend left] (13, 4.5);
\draw [thick, gray] (13, 2.5) to [bend right] (13, 4.5);
\draw [thick, gray] (7.5, 2.5) -- (13,2.5);
\draw [thick, gray] (7.5, 4.5) -- (13,4.5);
\draw [thick,LimeGreen] (12.2, 3.5) .. controls (11.7,3.5) and (12,4.5) .. (11.5, 4.5);
\draw [thick,LimeGreen, dashed] (11.5, 4.5) .. controls (11,4.5) and (11.5,2.5) .. (11,2.5);
\draw [thick,LimeGreen] (11,2.5) .. controls (10.5,2.5) and (11,3.5) .. (10.3, 3.5);
\draw [thick,LimeGreen] (9.5,3.5) .. controls (10,3.5) and (10,4.5) .. (9.5, 4.5);
\draw [thick, LimeGreen,dashed] (9.5, 4.5) .. controls (8.8,4.5) and (9.5,2.5) .. (8.5,2.5);
\draw [thick,LimeGreen] (8.5,2.5) .. controls (8,2.5) and (8.2,3.5) .. (8.7,3.5);
\draw [thick] (7.2, 3.5) -- (12.7, 3.5);
\filldraw [org] (12.2, 3.5) circle (1pt);
\filldraw [org] (9.5,3.5) circle (1pt);
\filldraw [blue] (10.3, 3.5) circle (1pt);
\filldraw [blue] (8.7,3.5) circle (1pt);
\draw [thick, {<[scale=2.5, length=1, width=2]}-] (8, 3.5) -- (8.2, 
3.5);
\draw [thick, gray, -{>[scale=2.5, length=1, width=2]}] (12.703,3.7) to  (12.706,3.73);
\node [gray] at (13, 3.5) {\tiny{$\mu_0$}};
\end{tikzpicture}
\caption{The orange and blue dots on the train track $B\cap \partial X$ are the endpoints of $\alpha$ and $\alpha_-$ respectively.}
\label{fig: figure eight train track}
\end{figure}

\begin{figure}[ht]
\begin{tikzpicture}[scale=0.8]
\draw [thick, org] (1.2, 6.2) -- (3.8,6.2);
\draw [thick, blue] (1.2, 8.8) -- (3.8, 8.8);
\draw [thick, LimeGreen] (1.2,6.2) -- (1.2, 8.8);
\draw [thick,LimeGreen] (3.8,6.2) -- (3.8, 8.8);
\node at (2.5, 7.5) {$D$};
\node at (7.5, 7.5) {$F'$};
\begin{scope}[thick,decoration={
    markings,
    mark=at position 0.5 with {\arrow{>}}}
    ] 
    \draw[postaction={decorate}] (7,6) -- (8,6);
    \draw[postaction={decorate}] (8,9) -- (7,9);
    \draw[postaction={decorate}] (9, 7) -- (9, 8);
    \draw[postaction={decorate}] (6,8) -- (6,7);
\end{scope}
\draw [thick, org] (6,7) -- (7,6);
\draw [thick, org] (9,8) -- (8,9);
\draw [thick, blue] (8,6) -- (9,7);
\draw [thick, blue] (6,8) -- (7,9);
\draw [thick, org] (1.2, 6.2) -- (3.8,6.2);
\draw [thick, blue] (1.2, 8.8) -- (3.8, 8.8);
\draw [thick, LimeGreen] (1.2,6.2) -- (1.2, 8.8);
\draw [thick,LimeGreen] (3.8,6.2) -- (3.8, 8.8);
\draw [thick,blue, -{[scale=2.5, length=1, width=2]>}] (8.62,6.38) -- (8.35,6.65);
\draw [thick,blue, -{[scale=2.5, length=1, width=2]>}] (6.62,8.38) -- (6.35,8.65);
\draw [thick,blue, -{[scale=2.5, length=1, width=2]>}] (2.5, 8.65) -- (2.5, 9);
\draw [thick,org, -{[scale=2.5, length=1, width=2]>}] (2.5, 6.35) -- (2.5, 6);
\draw [thick,org, -{[scale=2.5, length=1, width=2]>}] (6.38,6.38) -- (6.65,6.65);
\draw [thick,org, -{[scale=2.5, length=1, width=2]>}] (8.38,8.38) -- (8.65,8.65);
\end{tikzpicture}
\caption{Two disk sectors of the branched surface $B$.}
\label{fig: figure eight sectors}
\end{figure}

{\underline{Step 1: Building laminations on $\partial X$ fully carried by 
$B\cap \partial X$}}

\begin{figure}[ht]
 \centering
 \includegraphics[scale=0.46]{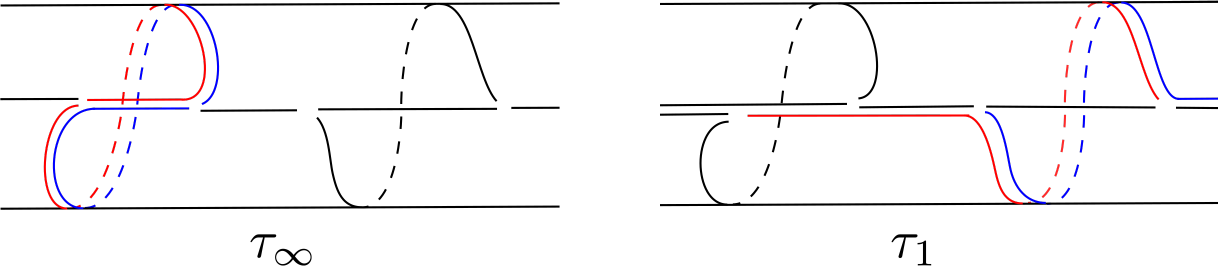}
 \caption{The laminations $\tau_{\infty}$ and $\tau_1$ are obtained by replacing each arc in the figure with its product with the standard Cantor set $\mathcal{C}$. 
}
 \label{fig:boundary lamination type1}
\end{figure}

\medskip
In Figure \ref{fig:boundary lamination type1}, we replace each arc by its 
product with the standard Cantor set $\mathcal{C}$ and merge the products 
using a homeomorphism $\mathcal{C} \sqcup \mathcal{C} \cong \mathcal{C}$. 

 It is easy to verify that the resulting lamination in the left-hand picture, denoted by $\tau_{\infty}$, detects the meridional slope $\mu_0$ and the lamination in the right-hand picture, denoted by $\tau_1$, detects the 
slope $1/1$. 
 
{\underline{Step 2: Extending laminations $\tau_\infty$ and $\tau_0$ into 
$X$}}

In this section, we will extend  $\tau_\infty$ and $\tau_1$ into $X$ to obtain laminations on $X$ detecting $\mu_0$ and the slope $1/1$  respectively. Since the constructions for $\tau_\infty$ and $\tau_1$ are identical, to simplify the notation, we will use $\tau_*$ to represent both $\tau_\infty$ and $\tau_1$.   

We continue to use the notation established in \S \ref{subsubsec:extend lamination}.

\begin{figure}[ht]
\centering 
\includegraphics[scale=0.5]{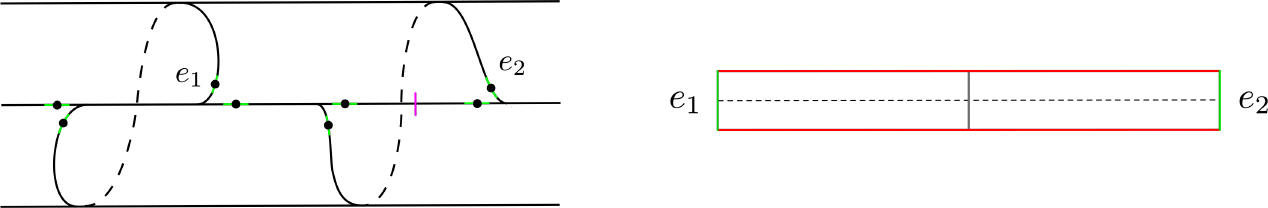}
\caption{The right-hand picture shows a component of $P(L')$ that connects arcs $e_1$ and $e_2$ in $B\cap \partial X$. This component is ``parallel'' to the $\alpha$ arc.}
\label{fig:figure eight L'}
\end{figure}

In this example, the graph $L'$ is the disjoint union of four arcs with endpoints on $\partial X$, two of which are parallel to  $\alpha$ and two are parallel to $\alpha_-$. The left-hand picture in Figure \ref{fig:figure eight L'} shows the endpoints of each component of $L'$.  It is clear that $P(L')$ consists of four rectangular sectors, each with two sides attached to $\partial B$ indicated in green in Figure \ref{fig:figure eight L'}. The lamination on $N(P(L'))$ is the product of the Cantor set and $P(L')$, which connects the corresponding parts of $\tau_\infty$ and $\tau_1$ 
in $N(\partial B)$ through homeomorphisms between Cantor sets.

\begin{figure}[ht]
\centering 
\includegraphics[scale=0.6]{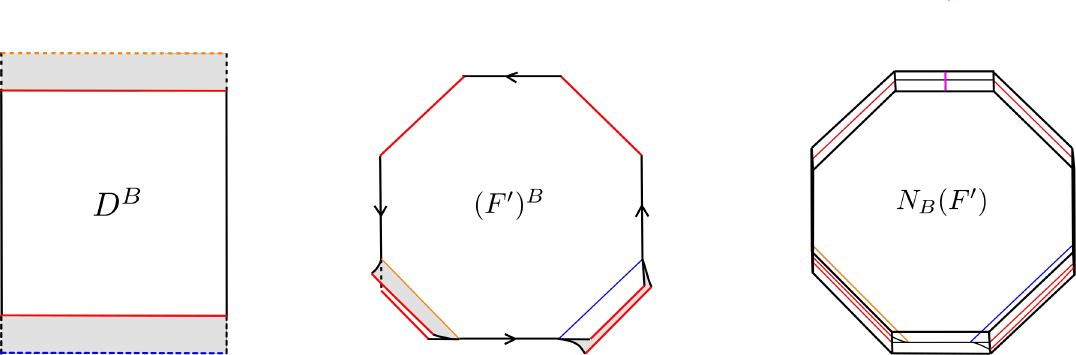}
\caption{The components of $B\setminus P(L')$, namely $D^B$ and $(F')^B$, and  $N_B(F')$  with $(F')^B$ embedded inside.  }
\label{fig:figure eight lamination}
\end{figure}

There are two disk sectors $D$ and $F'$ of $B$, as was shown in Figure \ref{fig: figure eight sectors}. Figure \ref{fig:figure eight lamination} shows the  corresponding pieces of $B\setminus P(L')$, namely $D^B$ and $(F')^B$. Note that the shaded portion of the product disk $D$ is part of $(F')^B$ and {\it not} part of $D^B$.  Following the discussion in \S \ref{subsubsec:extend lamination}, we have that $D^B$ is contained in the piece $N_B(D)$ of $N(B)$, which in our case is just $D^B\times I$. Figure \ref{fig:figure eight lamination} shows the piece $N_B(F')$ of $N(B)$ with $(F')^B$ embedded inside. 

The vertical boundaries $\partial_v N_B(D)$ and $\partial _v N_B(F')$ are 
endowed with the laminations induced by the lamination on $\partial_v N(P(L'))$ and by the lamination $\tau_\ast$ on $\partial X$. To complete the 
construction, we need to extend these laminations across $N_B(D)$ and $N_B(F')$. We will do this by applying \cite[Proposition 4.1]{Li1}.

 We first cut open and reglue $N(P(L'))$ along the $I$-fibres that are above the grey interval indicated in the right-hand picture of Figure  \ref{fig:figure eight L'}. Note that this component of $N(P(L'))$ is adjacent to $\partial_v N_B(D)$. We choose the regluing map so that the resulting lamination on $\partial_v N_B(D)$ is a union of circles, so we ``killed'' the potentially non-trivial holonomy map along $\partial D$ through this cutting-and-regluing operation by making the holonomy map of the new lamination on $\partial_v N_B(D)$ trivial.  We then extend the lamination into $N_B(D)$ with disks bounded by these circles.  
 
During the previous step, we didn't change the lamination $\tau_*$ on $\partial X$ but we did change  the lamination on $\partial_v N_B(F')$. To extend the newly created lamination on $\partial_v N_B(F')$ into the interior of $N_B(F')$, we will apply the same trick, though this time with respect to a vertical arc contained in $\partial_v N_B(F') \cap \partial X$.

More precisely, cut open $\partial_v N_B(F')$ across the pink arc, indicated in both the left-hand picture of Figure \ref{fig:figure eight L'} and 
the rightmost picture of Figure \ref{fig:figure eight lamination}, and reglue it so that the lamination on $\partial_v N_B(F')$ is a union of circles.  Since  the pink arc is disjoint from the closed leaves of $\tau_*$ as  illustrated in Figure \ref{fig:boundary lamination type1}, the new  $\tau_*$ still contains closed leaves of the same  rational slopes as before. We then extend the laminations into $N_B(F')$ using disk leaves.
 
We denote the laminations in $X$ constructed above by $\mathcal{L}_\infty$ and $\mathcal{L}_1$. They detect $\mu_0$ and the slope $1/1$ respectively. 

{\underline{Step 3: Constructing the desired foliations on $X$}}

The laminations $\mathcal{L}_\infty$ and $\mathcal{L}_1$ can be extended to foliations of $N(B)$ that are tangent to $\partial_h N(B)$ and transverse to $\partial_v N(B)$ 
by extending the product laminations on $N(P(L'))$, $N_B(D)$, $N_B(F')$ to product foliations and then gluing these foliations together by the same gluing maps that we  used in Step 2. 

Since $\overline{X \setminus N(B)}$ is a product sutured manifold, homeomorphic to a thickened annulus $S^1 \times I \times I$ where $S^1 \times I 
\times \{0, 1\}$ corresponds to $\partial_h N(B)$, it follows that we can 
extend the foliations of $N(B)$ obtained above to the entire manifold $X$ 
by filling in $\overline{X \setminus N(B)}$ with the parallel leaves $S^1 
\times I \times \{t\}$.  One can verify that the resulting foliations on $X$ are co-oriented and have no compact leaves. As such, they are taut. This completes the construction. 

\subsection{Case I: $\partial B$ is of Type 1}
\label{subsec: type 1}
\begin{figure}[ht]
 \centering
 \includegraphics[scale=0.55]{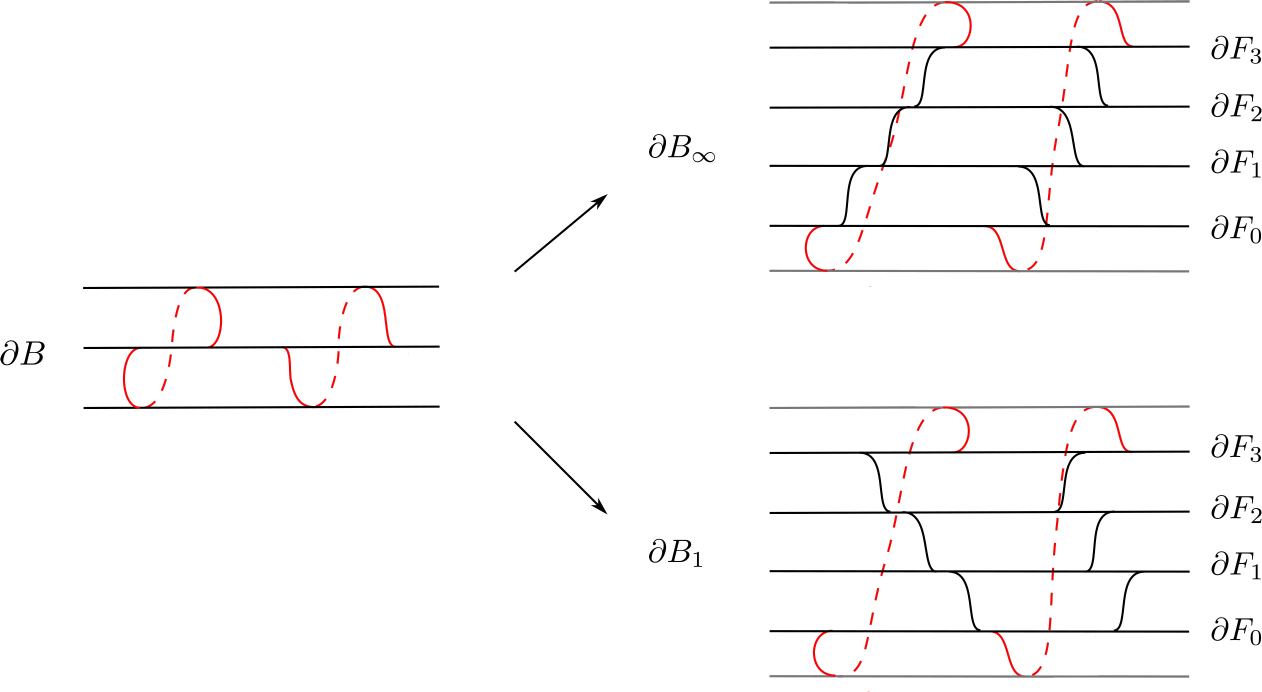}
 \caption{This figure shows the corresponding splittings of $\partial B$ to $\partial B_\infty$ and $\partial B_1$ for $n=4$, where  $\alpha = 
\alpha_4$ and $\alpha_- = \alpha_0$. The product disks in $B_\ast$ are given by 
$D_0 = \alpha_1\times [0, \frac{1}{4}]$, $D_1 = \alpha_2\times [\frac{1}{4}, \frac{1}{2}]$, $D_2 = \alpha_3\times [\frac{1}{2},\frac{3}{4}]$ 
and $D_3 = \alpha_4\times [\frac{3}{4}, 1]$ in $X = F\times [0,1]/\sim$. The sequence of arcs $\{\alpha_i: i=0,\dots, 4\}$ in $B_\infty$ (resp. $B_1$) is called a good positive (resp. negative) oriented sequence of winding number $0$ (resp. $1$) in \cite{Rob1}.}
 \label{fig: train track split type 1}
\end{figure}

By \cite[Theorem 4.1, Corollary 6.4]{Rob2} (also see \cite[\S 4]{Rob1}), a branched surface $B = \langle F; D\rangle$ whose boundary $\partial B$ is of Type 1 admits  two splittings, which we denote by $B_\infty$ and 
$B_1$. They carry co-orientable taut foliations that strongly detect  rational slopes in the intervals $(-\infty, 0]$ and $[0,1)$ respectively. We use $B_\ast$ to represent both $B_\infty$ and $B_1$. By \cite{Rob2}, $B_\ast$ satisfies the following: 
\begin{equation}
B_{\ast} = \langle F_0, F_1, \cdots, F_{n-1}: D_0, \cdots, D_{n-1} \rangle,
\end{equation}
where for each $i=0, \cdots, n-1$,  $F_i = F\times \frac{i}{n}$ is a copy  of $F$ in $X = F\times [0,1]/\hspace{-1.5mm} \sim$, $D_i = \alpha_{i+1}\times [\frac{i}{n}, \frac{i+1}{n}]$ is a product disk between $F_{i}$ and $F_{i+1}$, given by a collection of properly embedded non-separating simple arcs $\{\alpha_k\}_{k=0}^{n}$ on $F$ satisfying: 
\vspace{-.2cm} 
\begin{enumerate}
\item $\alpha_0 = \alpha_-$ and $\alpha_n = \alpha$.
\item $\alpha_i\cap \alpha_{i+1} = \emptyset$ for $i = 0, \cdots, n-1$; 

\vspace{.2cm} \item $\alpha_i\cup \alpha_{i+1}$ is non-separating  for $i 
= 0, \cdots, n-1$;  

\vspace{.2cm} \item for $i = 0, \cdots, n-1$, the endpoints of $\alpha_i$ and $\alpha_{i+1}$ are linked on $\partial F$ and the boundary train track $B_* \cap \partial X$ is  as depicted in Figure \ref{fig: train track split type 1}. 
\end{enumerate} 

We will work with the splittings $B_{\ast}$, since the sectors of $B_{\ast}$ are simple and we have better control of the changes to the lamination on $\partial X$ during the construction.

\begin{figure}[ht]
\centering
\includegraphics[scale=.61]{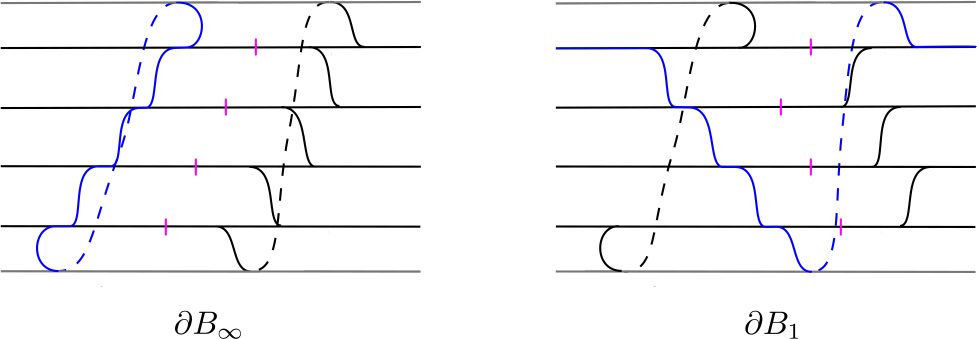}
\caption{The blue portion of the train track carries the detected slope.}
\label{fig:boundary lamination split type 1}
\end{figure}

As in Step 1 of \S \ref{subsubsec: figure eight} (see Figure \ref{fig:boundary lamination type1}), by replacing each sector of the train track $\partial B_\ast$ with its product with the Cantor set and merging the ends of the sectors together through homeomorphisms of the Cantor set, we can construct a lamination $\tau_\infty$ on $\partial X$ fully carried by $\partial B_\infty$ detecting the slope $\mu_0$, and a lamination $\tau_1$ on $\partial X$ fully carried by $\partial B_1$ detecting the slope $1/1$. 
We can also require that the closed leaves of $\tau_\ast$ are carried by the blue portion of the train track $\partial B_\ast$ indicated in Figure 
\ref{fig:boundary lamination split type 1} and disjoint from the rest of the train track. 

\begin{figure}[ht]
\centering
\includegraphics[scale=0.5]{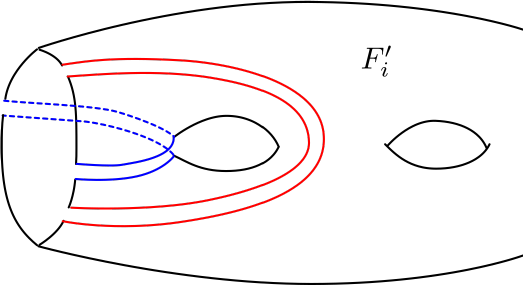}
\caption{The branch locus on $F_i$ consists of two disjoint arcs $\alpha_i$ and $\alpha_{i+1}$, and $F_i'$ is obtained by cutting $F_i$ open along 
$\alpha_i\sqcup\alpha_{i+1}$.}
\label{fig: sector type 1}
\end{figure}

For $\ast\in\{\infty, 1\}$, the branched surface $B_*$ has $n$ disk sectors corresponding to the product disks $D_i$ for $i = 0, \cdots , n -1$. 
It also has an additional $n$ sectors, denoted by $F'_i$, coming from cutting open $F_i$ along $\alpha_i \sqcup \alpha_{i+1}$. 
From Properties (1) - (3) of the arcs $\alpha_i$, it follows that $F_i'$ are connected surfaces with connected boundary of genus $g(F_i) - 1$ as 
shown in Figure \ref{fig: sector type 1}. 

Similar to the case dealt with in \S \ref{subsubsec: figure eight}, the graph $L'$ on $B_*$ is a disjoint union of $4n$ arcs. There are four arcs near each $F_i$ and parallel to the branch locus $\alpha_i \sqcup \alpha_{i+1}$ for $i = 0, \cdots, n-1$. Hence $P(L') = N(L') \cap B_*$ consists of $4n$ disjoint rectangular sectors with two sides of each rectangle attached to $B_*\cap \partial X$. The lamination on  $N(P(L'))$ is the product of the Cantor set and $P(L')$ and connects the parts of $\tau_\ast$ 
lying over $P(L') \cap \partial X$. 

For each sector $\mathcal{D}$ of $B_\ast$ we endow the vertical boundary $\partial_v N_{B_\ast}(\mathcal{D})$ with the lamination induced from its 
intersections with $\partial_v N(P(L'))$ and $\tau_\ast$. As in Step 2 of 
\S \ref{subsubsec: figure eight}, we can extend the lamination from $\partial_v N_{B_\ast}(D_i)$ and $\partial_v N_{B_\ast}(F'_i)$ into $N_{B_\ast}(D_i)$ and $N_{B_\ast}(F_i')$ using \cite[Proposition 4.1]{Li1}. Examples of places where we can cut open and reglue $\partial_v N_{B_\ast}(F'_i)$ without altering the closed leaves carried by $\tau_*$ are indicated by 
pink arcs in Figure \ref{fig:boundary lamination split type 1}.  This alternation of the laminations on $\partial_v N_{B_\ast}(F'_i)$ is only necessary when $g(F)=1$ and $g(F'_i) =0$. Otherwise, one can extend the laminations into $N_{B_\ast}(F'_i)$ without making any changes to the lamination on  $\partial_v N_{B_\ast}(F'_i)$ (see \cite[Lemma 3.2]{Li1}).

The rest of the proof is identical with Step 3 in \S \ref{subsubsec: figure eight}: the laminations constructed above can be extended to a foliation on $N(B)$ that is tangent to $\partial_h N(B)$ and transverse to $\partial_v N(B)$. Since the 
complement of $\overline{X\setminus N(B_*)}$ is a disjoint union of product sutured manifolds, we extend the foliations on $N(B)$ to taut foliations on $X$ by the product foliation on each component of $\overline{X\setminus N(B_*)}$, which detect the desired slopes. 

\subsection{Case II: $\partial B$ is of Type 2a} 
\label{subsec: type 2a}
\begin{figure}[ht]
\centering 
\includegraphics[scale=0.55]{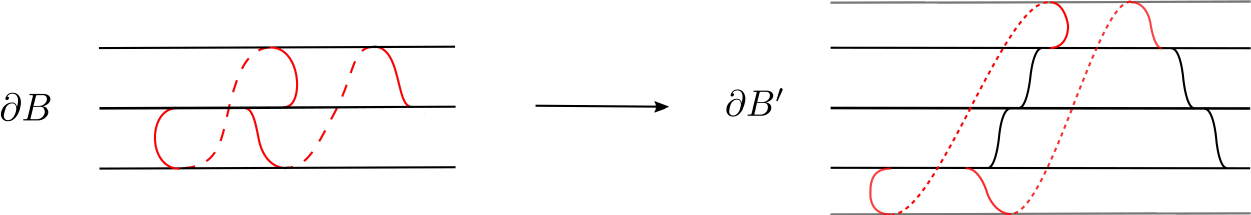}
\caption{An example of the boundary train track $\partial B'$ with $n=3$.}
\label{fig:train track type 2a}
\end{figure}

According to the proof of \cite[Corollary 4.4]{Rob2}, a branched surface $B = \langle F; D\rangle$ whose boundary $\partial B$ is of Type 2(a) admits a splitting 
\begin{equation}
B' = \langle F_0, F_1, \cdots, F_{n-1}: D_0, \cdots, D_{n-1} \rangle,
\end{equation}
where for each $i=0, \cdots, n-1$,  $F_i = F\times \frac{i}{n}$ is a copy  of $F$ in $X = F\times [0,1]/\sim$, $D_i = \alpha_{i+1}\times [\frac{i}{n}, \frac{i+1}{n}]$ is a product disk between $F_{i}$ and $F_{i+1}$, given by a collection of properly embedded non-separating simple arcs $\{\alpha_k\}_{k=0}^{n}$ on $F$ satisfying: 
\vspace{-.2cm} 
\begin{enumerate}

\item $\alpha_i\cap \alpha_{i+1} = \emptyset$ for $i = 0, \cdots, n-1$;

\vspace{.2cm} \item $\alpha_i\cup \alpha_{i+1}$ is non-separating  for $i 
= 0, \cdots, n-1$; 

\vspace{.2cm} \item for $i = 0, \cdots, n-1$, the endpoints of $\alpha_i$ and $\alpha_{i+1}$ are unlinked and the boundary train track $\partial B'$ when $n=3$ is shown in Figure \ref{fig:train track type 2a}. In \cite{Rob2}, the configuration of $\alpha_0 = \alpha_-$ and $\alpha_1$  on $F_0$ is called of Type 2a, and  $\alpha_i$ and $\alpha_{i+1}$ is called of Type 2b for $i = 1, \cdots, n-1$ (cf.Figure \ref{fig:types of 
boundary train track}).
\end{enumerate} 

\begin{figure}[ht]
\centering 
\includegraphics[scale=0.6]{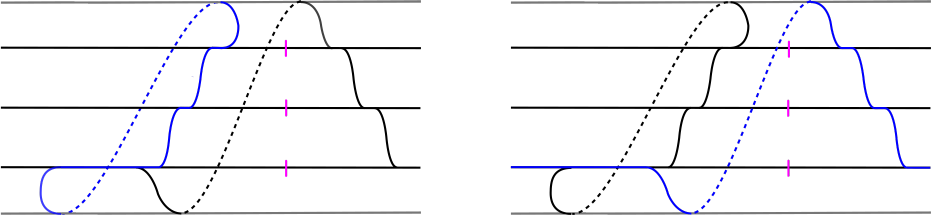}
\caption{The detected slope is carried by the blue portion of the train track.}
\label{fig:boundary lamination type 2a}
\end{figure}

In a similar fashion to the previous two cases, we can construct two laminations $\tau_\infty$ and $\tau_1$ on $\partial X$ fully carried by $\partial B'$ so that  $\tau_\infty$ and $\tau_1$ contains closed leaves of slopes $\mu_0$ and $1/1$ respectively, carried by the blue portion of the train track depicted in Figure \ref{fig:boundary lamination type 2a}. 

\begin{figure}[ht]
\centering
\includegraphics[scale=0.5]{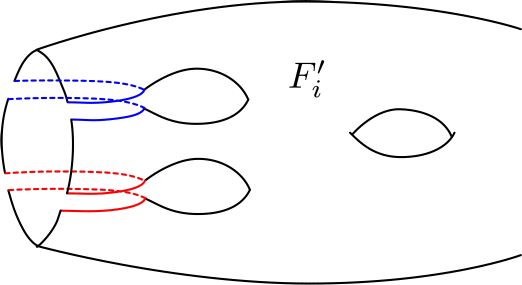}
\caption{Sector $F_i'$ of the branched surface $B'$.}
\label{fig:sector type 2a}
\end{figure}

The branched surface $B'$ has $n$ disk sectors corresponding to the product disks $D_i$ ($i = 0, \cdots , n -1$) and an additional $n$ sectors, denoted $F'_i$, obtained by cutting open the $F_i$ along $\alpha_i \sqcup \alpha_{i+1}$. 
From Properties (1) - (3) of the arcs $\alpha_i$, it follows that $F_i'$ is a connected surface with three boundary components whose genus is equal to $g(F)-2$ (Figure \ref{fig:sector type 2a}).

The graph $L'$ is a disjoint union of $4n$ arcs. We can equip $N(P(L'))$ with the product lamination that extends $\tau_\ast$ on $\partial X$. The 
vertical boundary $\partial_v N_{B'}(\mathcal{D})$ for  each sector $\mathcal{D}$ of $B'$ is then endowed with the lamination induced from $\partial_v N(P(L'))$ and $\tau_\ast$. For each product disk $D_i$, as before we 
apply \cite[Proposition 4.1]{Li1} to kill the holonomy along $\partial D_i$ without altering $\tau_*$ on $\partial X$ so that we can extend the lamination into $N_{B'}(D_i)$ by disk leaves.

\begin{figure}[ht]
\includegraphics[scale=0.48]{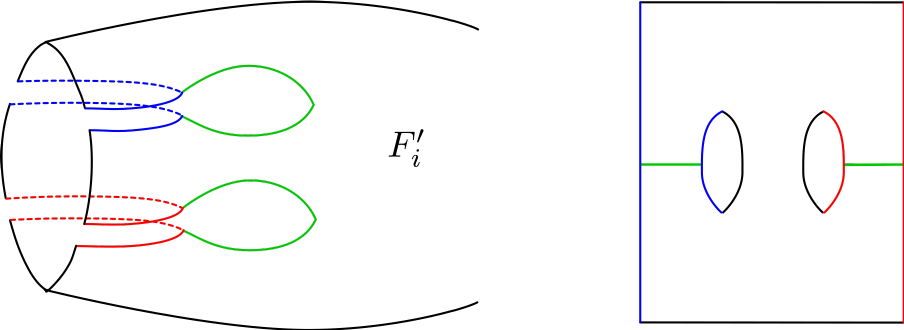}
\caption{Sector $F_i'$ and its boundary.}
\label{fig:partial S'_i}
\end{figure}

The right-hand picture in Figure \ref{fig:partial S'_i} shows the boundary of $F'_i$.  Since the lamination on $\partial_v N_{B'}(F_i')$ intersects each $I$-fibre in the standard Cantor set,  it can be trivially extended to the $I$-bundle over the regular neighbourhoods of the green arcs on $F'_i$ shown in Figure \ref{fig:partial S'_i}. This is the same as  attaching two thickened $1$-handles to $\partial_v N_{B'}(F'_i)$. The remaining 
part of $F_i'$, denoted by $F_i''$, is a connected surface of genus $g(F_i')$ with connected boundary.  As before,  we can cut open and reglue the laminations on $\partial_v N(F''_i)$ across the pink arcs on $\partial_v N(F''_i)\cap \partial X$ shown in Figure \ref{fig:train track type 2a} so 
that the new laminations on $\partial_v N(F''_i)$  are unions of circles which can be extended into  $N(F_i'')$ by parallel leaves homeomorphic to 
$F_i''$. Since the pink arcs are away from the closed leaves on $\partial 
X$, we have obtained laminations on $X$ that detect the same slopes as $\tau_0$ and $\tau_\infty$.

The rest of the construction is completed as in the previous cases.

\subsection{Case III: $\partial B$ is of Type 2b}
\label{subsec: type 2b}
By \cite[Corollary 4.3]{Rob2}, for any  rational slope $\alpha\in (-\infty,\infty)$, there exists a co-orientable taut foliation on $X$ strongly detecting the slope $\alpha$. So in this case, we just need to show that $\mu_0$ is foliation-detected.

\begin{figure}[ht]
\centering 
\includegraphics[scale=0.6]{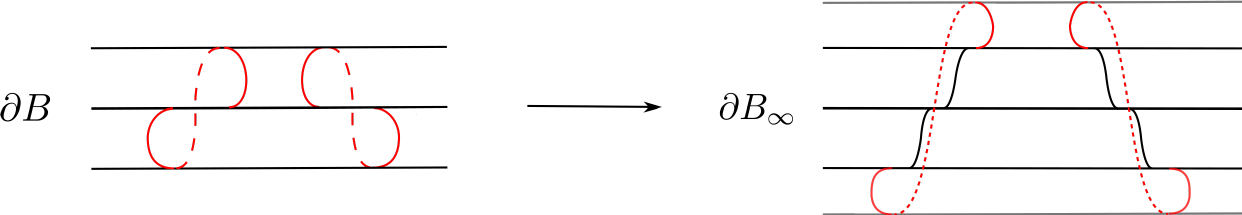}
\caption{An example of the boundary train track $\partial B_\infty$ with $n=3$.}
\label{fig:train track type 2b}
\end{figure}

According to the proof of \cite[Corollary 4.3]{Rob2}, a branched surface $B = \langle F; D\rangle$ whose boundary $\partial B$ is of Type 2(b) admits a splitting 
\begin{equation}
B_\infty = \langle F_0, F_1, \cdots, F_{n-1}: D_0, \cdots, D_{n-1} \rangle,
\end{equation}
where for each $i=0, \cdots, n-1$,  $F_i = F\times \frac{i}{n}$ is a copy  of $F$ in $X = F\times [0,1]/\sim$, $D_i = \alpha_{i+1}\times [\frac{i}{n}, \frac{i+1}{n}]$ is a product disk between $F_{i}$ and $F_{i+1}$, given by a collection of properly embedded non-separating simple arcs $\{\alpha_k\}_{k=0}^{n}$ on $F$ satisfying: 
\vspace{-.2cm} 
\begin{enumerate}

\item $\alpha_i\cap \alpha_{i+1} = \emptyset$ for $i = 0, \cdots, n-1$;

\vspace{.2cm} \item $\alpha_i\cup \alpha_{i+1}$ is non-separating  for $i 
= 0, \cdots, n-1$; 

\vspace{.2cm} \item for $i = 0, \cdots, n-1$, the endpoints of $\alpha_i$ and $\alpha_{i+1}$ are unlinked and the boundary train track $B_\infty\cap \partial X$ when $n=3$ is exactly as in Figure \ref{fig:train track type 2b}.

\end{enumerate}

\begin{figure}[ht]
\centering
\includegraphics[scale=0.43]{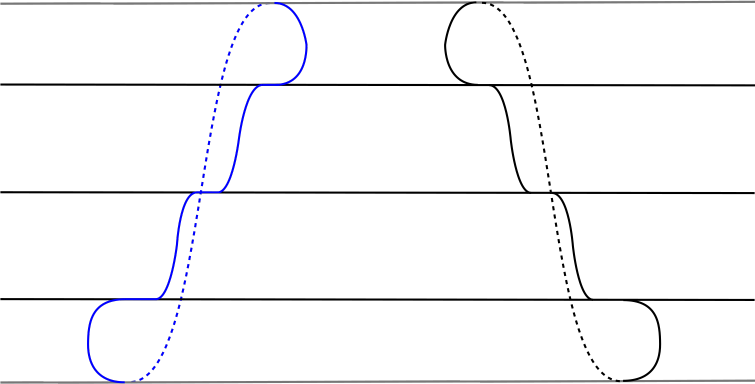}
\caption{The detected slope is carried by the blue portion of the train track.}
\label{fig:lamination type 2b}
\end{figure}

As before, we can construct a lamination $\tau_\infty$ on $\partial X$ carried by $\partial B_\infty$, which contains closed leaves of slope $\mu_0$ carried by the blue portion of the train track shown in Figure \ref{fig:lamination type 2b}.  

The branched surface $B_\infty$ has $n$ product disk sectors $D_i$ ($i = 
0, \cdots , n -1$) and an additional $n$ sectors, denoted by $F'_i$, obtained by cutting open $F_i$ along $\alpha_i \sqcup \alpha_{i+1}$. Properties (1) - (3) imply that $F_i'$ is a genus $g(F)-2$ connected surface with three boundary components, exactly the same as Case II in  \S \ref{subsec: type 2a}, Figure \ref{fig:sector type 2a}. The rest of the proof for this case is identical with the case in \S \ref{subsec: type 2a}. We omit the details.

\section{Proofs of results in Section 2} 
\label{sec: cbc and (*)} 

In this section we apply the material developed above to prove the results on toroidal $3$-manifolds discussed in Section \ref{sec: results}. 

\subsection{Toroidal manifolds with small order first homology} 
\label{subsec: tmwsofh}
Here we prove Theorem \ref{thm: small order lo}.

\begin{proof}[Proof of Theorem \ref{thm: small order lo}] Suppose that $W$ is a closed, irreducible $3$-manifold which contains an essential torus $T$ and for which $|H_1(W)| \in  \{1, 2, 3, 4, 6\}$.
As a rational homology $3$-sphere, $T$ separates $W$ into two rational homology solid tori $M_1, M_2$. The strategy is to show that there is a rational slope on $T$ which is $LO$-detected, respectively $CTF$-detected, in both $M_1$ and $M_2$ and then apply the gluing theorem (Theorem \ref{thm: * gluing}). To that end, let
\begin{itemize}

\item $T_1(M_i)$ be the torsion subgroup of $H_1(M_i)$ and $t_i = |T_1(M_i)|$;

\vspace{.2cm} \item $\lambda_i \in H_1(T)$ be the longitude of $M_i$ and $d_i \geq 1$ its order in $H_1(M_i)$;

\vspace{.2cm} \item $D = |\lambda_1 \cdot \lambda_2|$ be the distance between $\lambda_1$ and $\lambda_2$ on $T$. 

\end{itemize}
Without loss of generality we assume that $d_1 \leq d_2$ with equality implying that $t_1 \leq t_2$.

We claim that 
\begin{equation}
\label{eqn: order of homology} 
|H_1(W)| = D d_1 d_2 t_1 t_2
\end{equation}  
To see this, first note that $D \ne 0$ as otherwise the first Betti number of $W$ would be positive. Thus $\lambda_1$ and $\lambda_2$ generate a subgroup of $H_1(T)$ isomorphic to $\mathbb Z \oplus \mathbb Z$. Since $\lambda_1$ is primitive in $H_1(T)$, $H_1(T)/\langle \lambda_1,  \lambda_2 \rangle \cong \mathbb Z / D$. Also, $\langle \lambda_1,  \lambda_2 \rangle/\langle d_1 \lambda_1, d_2 \lambda_2 \rangle \cong \mathbb Z/d_1 \oplus \mathbb Z/d_2$, so the exact sequence
$$0 \to \langle \lambda_1, \lambda_2 \rangle/\langle d_1 \lambda_1, d_2 \lambda_2 \rangle \to H_1(T)/ \langle d_1 \lambda_1, d_2 \lambda_2 \rangle \to H_1(T)/\langle \lambda_1,  \lambda_2 \rangle\to 0$$
shows that 
\begin{equation}
\label{eqn: step2} 
|H_1(T)/ \langle d_1 \lambda_1, d_2 \lambda_2 \rangle| = D d_1d_2
\end{equation}  
Finally, consider the exact sequence of the pair $(W, T)$ 
\begin{equation}
\label{seq: wt} 
0 \to H_2(W, T) \xrightarrow{\; \partial \;}  H_1(T) \to H_1(W) \to H_1(W, T) \to 0
\end{equation}
The isomorphism $H_2(W, T) \cong H_2(M_1, T) \oplus H_2(M_2, T)$ implies that the image of the homomorphism $\partial$ in (\ref{seq: wt}) is $\langle d_1 \lambda_1, d_2 \lambda_2 \rangle$. Further, $H_2(M_i) \cong \{0\}$ since $M_i$ is a rational homology solid torus, so by Lefschetz duality and universal coefficients we see that $H_1(W, T) \cong H_1(M_1, T) \oplus H_1(M_2, T) \cong H^2(M_1) \oplus H^2(M_2) \cong T_1(M_1) \oplus T_1(M_2)$. Plugging these two facts and (\ref{eqn: step2}) into (\ref{seq: wt}) then implies that (\ref{eqn: order of homology}) holds. 

Since $d_i$ divides $t_i$, $D d_1^2 d_2^2 $ divides $|H_1(W)|$, so our bounds on the latter imply that $d_1 = d_2 = 1$ when $|H_1(W)| \ne 4$. Further, it is easy to verify that if $d_1 < d_2$, then $D = d_1 = t_1 = 1 < d_2 = t_2 = 2$. It follows that for all values of $|H_1(W)|$ we have, 
\begin{itemize}

\item  $d_1 = 1$ and $t_1 \leq 2$, and 

\vspace{.2cm} \item $D > 1$ implies that $d_2 = 1$ and $t_2 \leq |H_1(W)|/2$.

\end{itemize}
Note that  $d_i = 1$ if and only if $\lambda_i$ is integrally nullhomologous in $M_i$, and $t_i \leq 2$ implies that $T_1(M_i)$ is a $\mathbb Z/2$ vector space.  

To prove Theorem \ref{thm: small order lo}(1), suppose that $|H_1(W)| \leq 4$. In the case that $D = 1$, the fact that $d_1 = 1$ and $t_1 \leq 2$ combines with Theorem \ref{thm: meridional detn}(2) to show that $[\lambda_2] \in \mathcal{D}_{LO}(M_1)$. As $[\lambda_2] \in \mathcal{D}_{LO}(M_2)$ always holds, Theorem \ref{thm: * gluing} implies that $W$ has a left-orderable fundamental group. If $D > 1$, then $d_1 = d_2 = 1$, and we leave it to the reader to verify that our conditions on $|H_1(W)|$ imply that there is a rational slope $[\alpha]$ of distance $1$ from both $[\lambda_1]$ and $[\lambda_2]$. Since $t_1 \leq 2$, $[\alpha]  \in \mathcal{D}_{LO}(M_1)$, and as $t_2 \leq |H_1(W)|/2 \leq 2$, $[\alpha]  \in \mathcal{D}_{LO}(M_2)$, so $\pi_1(W)$ is left-orderable.

The proof of Theorem \ref{thm: small order lo}(2)(a), where we allow the possibility that $|H_1(W)| =6$ but assume that both $M_1$ and $M_2$ fibre over the circle, is similar. If $D=1$, Theorem \ref{thm: meridional detn}(3) shows that $[\lambda_2]\in \mathcal{D}_{CTF}(M_1)$. If $D>1$, there is a rational slope $[\alpha]$ of distance $1$ from both $[\lambda_1]$ and $[\lambda_2]$. By Theorem \ref{thm: meridional detn}(3), $[\alpha]\in \mathcal{D}_{CTF}(M_i)$ for both values of $i$, so $W$ admits a co-orientable taut foliation by Theorem \ref{thm: * gluing}. 

For Theorem \ref{thm: small order lo}(2)(b), we assume that $M_1$ is fibred. Since $W$ is an integer homology $3$-sphere, $D=1$ and $[\lambda_2]$ is in $\mathcal{D}_{CTF}(M_1)$ by Theorem \ref{thm: meridional detn}(3). Since $[\lambda_2]\in \mathcal{D}_{CTF}(M_2)$ always holds, Theorem \ref{thm: * gluing} implies that $W$ admits a co-oriented taut foliation. 
\end{proof}

\begin{remark} 
\label{rem: why not 6} 
A similar analysis shows that the fundamental group of $W$ is left-orderable when $W$ is toroidal and $|H_1(W)| = 6$ except, perhaps, when $d_1 = d_2 = 1, D = 2, t_1 = 1, t_2 = 3$. Much of the argument goes through in this case as well: There is a rational slope $[\alpha]$ of distance $1$ from both $[\lambda_1]$ and $[\lambda_2]$ and Theorem \ref{thm: meridional detn}(2) shows that $[\alpha]  \in \mathcal{D}_{LO}(M_1)$, but we cannot apply it to conclude $[\alpha]  \in \mathcal{D}_{LO}(M_2)$ since $T_1(M_2) \cong \mathbb Z/3$, which potentially prevents us from lifting the universal circle representation $\rho_{M_2}: \pi_1(M_2) \to \mbox{Homeo}_+(S^1)$ to $\mbox{Homeo}_{\mathbb Z}(\mathbb R)$ (cf. Theorem \ref{thm: universal circle fixed point}). 
\end{remark}

\begin{proof}[Proof of Proposition \ref{prop: gao toroidal}]
Let $M_1$ be the trefoil exterior and $M_2$ the hyperbolic integer homology solid torus $m137$ (\cite{CHW}). We can write $\pi_1(M_1) = \langle a,b : a^2 = b^3 \rangle$, where $\pi_1(\partial M_1)$ has basis $\mu_1$ and $a^2 = \mu_1^{6} \lambda_1$, where $\mu_1$ is a meridional class and $\lambda_1$ a longitudinal class. Fix a basis $\mu_2, \lambda_2$ for $\pi_1(\partial M_2)$, where $\lambda_2$ is a longitudinal class. Gao has shown that there is a (negative) integer $N$ such that if $n < N$, the fundamental group of the homology sphere $M_2(\mu_2 + n \lambda_2)$ has no non-trivial $PSL(2, \mathbb R)$ representation (\cite[Theorem 1]{Gao}). Fix $n < N$ and let $W$ be the homology sphere $M_1 \cup M_2$, where the boundaries are glued so that $\mu_1$ is identified with $\lambda_2$ and $a^2$ with $\mu_2 \lambda_2^n$, and suppose that $\rho : \pi_1(W) \to PSL(2, \mathbb R)$ is a homomorphism. We show that $\rho$ is the trivial representation. 

Let $\Gamma$ be the subgroup of $PSL(2, \mathbb R)$ consisting of the classes of matrices $\pm (a_{ij})$ with $a_{21} = 0$. By conjugating $\rho$ we may suppose that $\rho(b) \in \Gamma$. Since $a^2 = b^3$, it follows that either $\rho(a) \in \Gamma$ or $\mbox{trace}(\rho(a)) = 0$, which implies $\rho(a^2) = 1$. In the latter case, $\rho(\mu_2 \lambda_2^n) = 1$, hence $\rho(\pi_1(M_2)) = 1$ by Gao's result, which implies that $\rho$ is trivial. In the former, $\rho(\pi_1(M_1)) \leq \Gamma$, and since the second derived subgroup $\Gamma''$ of $\Gamma$ is trivial and $\lambda_1 \in \pi_1(M_1)'', \rho(\lambda_1) = 1$, i.e. $\rho(a^2 \mu_1^{-6}) = 1$. But then $\rho(\mu_2 \lambda_2^{n-6}) = 1$. As before, this implies that $\rho$ is trivial.
\end{proof}

\subsection{Irreducibility and boundary-incompressibility of branched covers}
A link $L$ in a $3$-manifold $W$ is {\it trivial} if it is connected and bounds a disk in $W$.

For $i = 1,2$, let $L_i$ be a link in $W_i$. By choosing a 3-ball $B_i \subset W_i$ such that $(B_i, B_i \cap L_i)$ is homeomorphic to the standard pair $(B^3, B^1)$, we can define a {\it connected sum} $L_1 \# L_2$ in $W_1 \# W_2$ in the usual way. A link is {\it prime} if it is not a connected sum of two non-trivial links.

Let $\widetilde W$ be a 3-manifold with an orientation-preserving $\mathbb {Z}/n$ action, generated by $h$, say, such that the fixed point set Fix$(h)$ is a link $\tilde L$ and the action is free on $\widetilde W \setminus \tilde L$. Then $\tilde L$ has an $h$-invariant tubular neighbourhood $N(\tilde L)$ such that, for each component $\tilde K_i$ of $\tilde L$, $h$ acts on a normal disk in $N(\tilde K_i)$ by rotation through $2\pi n_i/n$ for some integer $n_i$ coprime to $n$. Let $p : (\widetilde W, \tilde L) \to (W, L)$ be the quotient map of the action; then $W$ is a 3-manifold and $L$ is a link in $W$. We say that $\widetilde W$ is an {\it $n$-fold cyclic branched cover} of $(W, L)$.

\begin{proposition}
\label{prop: irr bccs} 
Let $W$ be a compact 3-manifold whose boundary is a (possibly empty) disjoint union of tori. Let $L$ be a prime link in $W$ whose exterior is irreducible, and let $\widetilde W$ be an $n$-fold cyclic branched cover of $(W, L)$. Then
\begin{enumerate}[leftmargin=*] 
\setlength\itemsep{0.5em}
\item[{\rm (1)}] either $\widetilde W$ is irreducible or $W$ contains a non-separating 2-sphere meeting $L$ transversely in two points;
\item[{\rm (2)}] either $\partial \widetilde W$ is incompressible or $(W, L) \cong (S^1 \times D^2, S^1 \times \{0\})$.
\end{enumerate}
\end{proposition}

\begin{proof}
Let $X = W \setminus \mbox{int} N(L)$ be the exterior of $L$; so $\widetilde X = p^{-1}(X)$ is the exterior of $\tilde L$. The theorem is clearly true if $L$ is trivial, for then $W \cong S^3$ since $X$ is irreducible, and so $\widetilde W \cong S^3$.

Since $\widetilde X$ is an $n$-fold cyclic cover of $X$ and $X$ is irreducible, $\widetilde X$ is irreducible by the Equivariant Sphere Theorem; see \cite[Theorem 3]{MSY}.

(1) Suppose $\widetilde W$ is reducible. Then, by the equivariant sphere theorem there is an essential 2-sphere $S \subset \widetilde W$ such that either $h(S) \cap S = \emptyset$ or $h(S) = S$.

If $h(S) \cap S = \emptyset$ then $S$ is disjoint from $\mbox{Fix}(h) = \tilde L$ and hence can be isotoped into $\widetilde X$, contradicting the fact that $\widetilde X$ is irreducible. So suppose $h(S) = S$. Then either $n = 2$ and $S \cap \tilde L$ is a component $\tilde K$ of $\tilde L$, or $S$ meets $\tilde L$ transversely in two points.

In the first case, $p(S)$ is a disk $D$ in $W$ such that $D \cap L = \partial {D} \cap L$ is a component $K$ of $L$. Thus either $L = K$ is trivial, or the boundary of a regular neighbourhood of $D$ is a 2-sphere in $W$ separating $K$ and $L \setminus K$, contradicting the irreducibility of $X$.

In the second case, $p(S) = S'$ is a 2-sphere in $W$ meeting $L$ transversely in two points. If $S'$ is non-separating we have the second conclusion in (1). So suppose $S'$ separates $W$. Then it induces a decomposition $(W, L) \cong (W_1, L_1) \# (W_2, L_2)$, and $S$ induces a corresponding decomposition
$(\widetilde W, \tilde L) \cong (\widetilde W_1, \widetilde L_1) \# (\widetilde W_2, \tilde L_2)$. Since $L$ is prime, $L_2$ (say) is trivial in $W_2$. Then $\tilde L_2$ is trivial in $\widetilde W_2$, and hence $S$ can be isotoped off $\tilde L$. This contradicts the irreducibility of $\tilde X$.

(2) Note that $\partial X = \partial W \sqcup \partial N(L)$. We claim that $\partial W$ is incompressible in $X$. For otherwise $X \cong (S^1 \times D^2) \# Y$, where $\partial Y \supset \partial N(L)$, contradicting the irreducibility of $X$.

Suppose $\partial \widetilde W$ is compressible in $\widetilde W$. Then by the equivariant loop theorem/Dehn's lemma \cite{MY} there is a compressing disk $D$ for $\partial \widetilde W$ such that for each $k$, $1 \le k < n$, either $h^{k}(D) \cap D = \emptyset$ or $h^{k}(D) = D$. It is easy to see that this implies that either $h^{k}(D) \cap D = \emptyset$ for all $k$ such that $1 \le k < n$, or $h(D) = D$.

In the first case we can assume that $D \subset \widetilde X$, and hence $p(D)$ is a compressing disk for $\partial W$ in $X$, contradicting the claim above.

In the second case, $D$ meets $\tilde L$ transversely in a single point. Then $D' = p(D)$ is a compressing disk in $W$ for a component $T$ of $\partial W$, which meets $L$ transversely in a single point. Surgering $T$ along $D'$ gives a 2-sphere $S'$ which, when pushed into $\mbox{int} \; W$ induces a decomposition $(W, L) \cong (S^1 \times D^2, S^1 \times \{0\}) \# (W', L')$. Since $L$ is prime, $L'$ is trivial in $W'$, so $S'$ can be isotoped off $L$. Since $X$ is irreducible, $S'$ then bounds a 3-ball in $X$, so $(W, L) \cong (S^1 \times D^2, S^1 \times \{0\})$.
\end{proof}

\subsection{Cyclic branched covers of links with companion manifolds in homology $3$-spheres} 
\label{subsec: branched cover toroidal}
Consider a link $L$ in an integer homology $3$-sphere $W$ whose exterior is irreducible and contains a companion manifold $M$. Then $W = V \cup_T M$, where $T = \partial V = \partial M$ and $L \subset \mbox{int}(V)$. Note that $V$ and $M$ are integer homology solid tori.  

Let $\Sigma_n(L)$ denote the $n$-fold cyclic cover of $W$ branched over $L$.

\begin{theorem} 
\label{thm: nls detn zhs} 
Suppose that $L$ is a prime link in an integer homology $3$-sphere $W$ whose exterior is irreducible and contains a companion manifold $M$. Then $\Sigma_{n}(L)$ is not an $L$-space and has a left-orderable fundamental group for each $n \geq 2$. Furthermore, each $\Sigma_{n}(L)$ admits a co-oriented taut foliation if the longitudinal slope of $V$ is foliation-detected in $M$. In particular this holds if $M$ is fibred.
\end{theorem}

\begin{proof}
Our hypotheses imply that $\Sigma_n(L)$ is irreducible (Proposition \ref{prop: irr bccs}), so without loss of generality we assume that $\Sigma_n(L)$ is a rational homology $3$-sphere. 

Let $\lambda_V$ and $\lambda_M$ be oriented simple closed curves on $T$ representing the longitudinal classes of $V$ and $M$. 
Since $W$ is an integer homology $3$-sphere, the slopes $[\lambda_V]$ and $[\lambda_M]$ are of distance $1$ from each other. Hence $\lambda_V$ and $\lambda_M$ represent a basis of $H_1(T)$. 

The curve $\lambda_V$ bounds an orientable surface $F$ in $V$ which we can suppose is transverse to $L$ and has minimal intersection with it. By a 
mild abuse of notation we also use $\lambda_V$ to denote the class in $H_1(T)$ it carries. The intersection number $w$ of $L$ with $F$ can be thought of as either the winding number of $L$ in $V$ or the linking number of $\lambda_V$ and $K$. 

Fix an integer $n \geq 2$ and set $d = \gcd(n, w)$ and $m = n/d$. We use $V_n$ to denote the $n$-fold cyclic branched cover of $(V, L)$ 
and observe that $\partial V_n = \bigsqcup_{i=1}^d T_i$ where each $T_i$ is an incompressible torus (Proposition \ref{prop: irr bccs}). If $\widehat F$ denotes the $n$-fold cyclic branched cover of $(F, F \cap K)$ contained in $V_n$, then $\partial \widehat F \cap T_i$ is an essential simple closed curve $\hat \lambda_i$, say. If $Y$ is a manifold obtained by attaching copies of the $HF$-generalised solid torus $v2503$ to $V_n$ so that their rational longitudes are identified with the $\hat \lambda_i$, then $Y$ has a positive first Betti number, since $\widehat F$ is non-separating in $V_n$. Proposition \ref{prop: ctf implies lo detd genl} then shows that $([\hat \lambda_1], [\hat \lambda_2], \ldots , [\hat \lambda_d])$ is $\ast$-detected in $V_n$  for $\ast = LO$ or $NLS$. It is $CTF$-detected by \cite{Gab1}. 

Let $M_m$ denote the $m$-fold cyclic cover of $M$, and $\mu_m \subset \partial M_m$ the inverse image of $\lambda_V$. Then 
\begin{displaymath}
\Sigma_n(L) = V_n \cup \big(\bigsqcup_{i=1}^d M^i_m\big).
\end{displaymath}
where the $M^i_m$ are copies of $M_m$ glued to $V_n$ along $T_i$ in such a way that  $\hat \lambda_i$ is identified with the copy $\mu_m^i$ of $\mu_m$ on $\partial M_m$.  Each $M^i_m$ is a rational homology solid torus since we have assumed that $\Sigma_n(L)$ is a rational homology $3$-sphere. 
The reader will verify that $\lambda_M$ lifts to the longitude  $\lambda_{M_m^i}$ of $M_m^i$, and that $\lambda_{M_m^i}$ is integrally null-homologous in $M_m^i$. Further, as $[\lambda_V]$ and $[\lambda_M]$ are of distance $1$ from each other, $[\mu_m^i]$ is of distance $1$ from the longitude $[\lambda_{M_m^i}]$ of $M_m^i$. Then Proposition \ref{prop: mu 
+ k lambda NLS-detd} shows that $[\mu_m^i] \in \mathcal{D}_{NLS}(M_m^i)$. 
It follows that $([\hat \lambda_1], [\hat \lambda_2], \ldots, [\hat \lambda_d])$ is an $NLS$-gluing coherent family of  rational slopes and therefore Theorem \ref{thm: general * gluing} implies that $\Sigma_n(L)$ is not an $L$-space.

By Theorem \ref{thm: meridional detn}(2) $[\lambda_V]$ is order-detected in $M$, by a left-order $\mathfrak{o}$, say, on $\pi_1(M)$. Then the restriction of $\mathfrak{o}$ to $\pi_1(M_m^i)$ order-detects $[\mu_m^i]$. Then, exactly as above, Theorem \ref{thm: general * gluing} shows that $\Sigma_n(L)$ has a left-orderable fundamental group.

Finally, by hypothesis, $[\lambda_V]$ is contained in $\mathcal{D}_{CTF}(M)$. Thus there is a co-oriented taut foliation $\mathcal{F}$ on $M$ detecting $[\lambda_V]$. It is easy to see that the pull-back of $\mathcal{F}$ to $M_m^i$ foliation-detects $[\mu_m^i]$, and a similar application of Theorem \ref{thm: general * gluing} shows that $\Sigma_n(L)$ admits a co-oriented taut foliation.
\end{proof}

\begin{proof}[Proof of Theorem \ref{thm: detn and bccs}]
Theorem \ref{thm: nls detn zhs} immediately implies parts (1) and (2) of Theorem \ref{thm: detn and bccs}. Part (3) is a consequence of Theorem \ref{thm: nls detn zhs} combined with part (3) of Theorem \ref{thm: meridional detn}.
\end{proof}

\subsection{Cyclic branched covers of satellite  links with braided patterns}
\label{subsec: satellite 2}
Here we prove those parts of Theorem \ref{thm: results GL conjecture} dealing with satellite  links with braided patterns. 

Suppose $P$ is the closure of a $w$-braid in $V = S^1 \times D^2$, $w \ge 2$. Fix a meridional disk $D = D^2 \times \{*\}$ and set ${\bf p} = P \cap D^2$, a collection of $w$ points in the interior of $D$.  The $n$-fold cyclic cover $V_n$ of $(V, P)$ fibres over the circle with fibre an $n$-fold cyclic cover $\widehat F = \Sigma_n(D, {\bf p}) \to D$ branched over ${\bf p}$. Note that $\widehat F$ has $d = \gcd(w, n)$ boundary components. 

\begin{lemma} 
 \label{lem: genus branched surface}
The fibre $\widehat F = \Sigma_n(D, {\bf p})$ of $V_n$ has genus $g = \frac12 \big((n - 1)(w - 1) - (d-1) \big)$. Hence $g = 0$ if and only if 
$n = w =2$. 
\end{lemma}

\begin{proof}
Since $\chi(D \setminus {\bf p}) = 1-w$, we have 
$$2 - 2g - d = \chi(\widehat F) = n\chi(D \setminus {\bf p}) + w = n(1 - w) + w = 1 - (n  - 1)(w - 1),$$ 
which gives  
$$g = \frac12 \big((n - 1)(w - 1) - (d-1) \big).$$
Thus $g = 0$ if $n = w = 2$. On the other hand, if $g = 0$ then 
$$(n - 1)(w - 1) = d-1 \leq \min \{n-1, w - 1\},$$
and therefore $n = w = 2$. 
\end{proof}

Following Delman and Roberts \cite{DR1}, we say that a knot $K$ in $S^3$ is {\it persistently foliar} if each non-meridional rational slope on the exterior of $K$ is strongly foliation-detected. Conjecturally (see \cite[Conjecture 1.8]{DR2}), $K$ is persistently foliar if $K$ is not an $L$-space knot (in particular, if $K$ is not fibred) and has no reducible surgery. Examples of knots known to be persistently foliar are non-torus alternating knots, Montesinos knots, and composite knots with a summand that is either non-torus alternating, Montesinos, or non-prime fibred \cite{DR2}. Note that if $M_1$ and $M_2$ are exteriors of persistently foliar knots then any closed manifold of the form $M_1 \cup_{\partial} M_2$ admits a co-oriented taut foliation. 

The following proposition is the first assertion of Theorem \ref{thm: results GL conjecture}(3)(b). 

\begin{proposition} 
 \label{prop: CTF braided satellites}
Let $P(K)$ be a prime satellite  link in the $3$-sphere where the pattern $P$ is braided and the companion knot $K$ is persistently foliar. Then $\Sigma_n(P(K))$ admits a co-oriented taut foliation for all $n \ge 2$. 
\end{proposition}

\begin{proof} 
As in the proof of Theorem \ref{thm: nls detn zhs}, 
\begin{displaymath}
\Sigma_n(P(K)) = V_n \cup \big(\bigsqcup_{i=1}^d M^i_m\big),
\end{displaymath}
where $V_n$ is irreducible and the $M^i_m$ are copies of $M_m$, the $m$-fold cyclic cover of the exterior $M$ of $K$, glued to $V_n$ along the torus $T_i = \partial M_m^i$. The intersection $\partial \widehat F \cap T_i$ is an essential simple closed curve $\hat \lambda_i$. 

First assume that either $n > 2$ or $w > 2$, and therefore the fibre $\widehat F$ of $V_n$ has positive genus (Lemma \ref{lem: genus branched surface}). By \cite[Theorem 1.1]{KR}, there is a neighbourhood $U$ of  rational multislopes on $\partial V_n$ about the multislope $[\partial \widehat F] = ([\hat \lambda_1], [\hat \lambda_2], \ldots , [\hat \lambda_d])$ such that $[\beta] \in U$ implies $[\beta] \in \mathcal{D}^{str}_{CTF}(V_n)$ (cf. \S \ref{subsec: fln-detn multislopes}). If $[\mu]$ denotes the meridional slope of $K$, then as $K$ is persistently foliar, any rational slope on $\partial M$ not equal to $[\mu]$ is strongly $CTF$-detected in $M$. Hence any rational slope on $\partial M_m^i$ not equal to $[\mu_m^i]$ is strongly $CTF$-detected in $M_m^i$. Since $\hat \lambda_i$ is identified with $\mu^{i}_m \subset \partial M^{i}_m$, we see that there exists a  rational multislope $[\beta] = ([\beta_1],[\beta_2],...,[\beta_d]) \in U$ such that $[\beta_i] \in \mathcal{D}^{str}_{CTF}(M^{i}_m)$, $1 \le i \le d$. Hence $\Sigma_{n}(P(K))$ admits a co-oriented taut foliation.

 Now suppose $n = w = 2$. Then $\Sigma_n(D,{\bf p})$ is an annulus and so $V_2 = \Sigma_2(V,P)$ is homeomorphic to $T^2 \times I$. Therefore $\Sigma_2(P(K)) \cong M \cup_{\partial}  M'$, where $M'$ is a copy of $M$. Since $K$ is persistently foliar, it follows that $\Sigma_{2}(P(K))$ admits a co-oriented taut foliation by the remark immediately before Proposition \ref{prop: CTF braided satellites}.
\end{proof}

The second assertion of Theorem \ref{thm: results GL conjecture}(3)(b) follows the next proposition.

\begin{proposition} 
 \label{prop: CTF braided satellites 2}
Let $P(K)$ be a prime satellite  link where the pattern $P$ is braided of winding number $w$. Then $\Sigma_n(P(K))$ admits a co-oriented taut foliation for each $n \ge 2$ which has a prime factor not dividing $w$. 
\end{proposition}

\begin{proof}
We continue to use the notation developed above.

Without loss of generality we can suppose that $\Sigma_n(P(K))$ is a rational homology $3$-sphere. Then $V_n$ is a rational homology solid torus. 

Suppose that $p$ is a prime factor of $n$ which is coprime to $w$ and write $n = kp$. Then $\partial V_p$ is connected, so all  rational slopes distance $1$ from its fibre slope are $CTF$-detected (Theorem \ref{thm: meridional detn}(3)). In particular, this applies to the longitudinal slope of $M_p$, the $p$-fold cyclic cover of the exterior of $K$. So $\Sigma_p(P(K))$ admits a co-oriented taut foliation by Theorem \ref{thm: fln gluing}. Since the foliation constructed in Theorem \ref{thm: meridional detn}(3) is carried by a branched surface in $V_p$ consisting of a fibre surface and a product disk, it follows that the inverse image $\widehat P$ of $P$ in $V_p$ is transverse to this foliation, and so by taking the $k$-fold cyclic cover of $\Sigma_p(P(K))$ branched over $\widehat P$ we deduce that $\Sigma_{n}(P(K)) = \Sigma_{kp}(P(K))$ admits a co-oriented taut foliation, which completes the proof.
\end{proof}

\subsection{Cyclic branched covers of toroidal links in integer homology $3$-spheres}
\label{subsec: proof sigma2 toroidal link nls}
Here we prove Theorems \ref{thm: theorem y} and \ref{thm: toroidal links nls}. 

\begin{proof}[Proof of Theorem \ref{thm: theorem y}]
Equip $L$ with any orientation and fix $n \geq 2$. Our hypotheses imply that each cyclic branched cover of $L$ is irreducible (Proposition \ref{prop: irr bccs}), and so without loss of generality we may assume that $\Sigma_n(L)$ is a rational homology $3$-sphere; otherwise it has a left-orderable fundamental group by \cite[Theorem 1.1]{BRW}. 

Let $T$ be an essential torus in the exterior of $L$ and write $W = X_1 \cup_T X_2$,  where $T$ is the boundary of both $X_i$.  Note that if $L_i = \emptyset$ for some $i$, then the proof of Theorem \ref{thm: detn and bccs} shows that $\Sigma_n(L)$ has a left-orderable fundamental group, so we assume that  $L_i \ne \emptyset$ for both $i$.

Let $p: \Sigma_n(L) \rightarrow W$ denote the branched cover. Since $L_i \ne \emptyset$ for both $i$ and we assumed that $\Sigma_n(L)$ is a rational homology $3$-sphere, $p^{-1}(T)$ must be connected. Set $\widetilde{T} = p^{-1}(T)$ and $M_i = p^{-1}(X_i)$, so $\Sigma_n(L) = M_1 \cup M_2$ and $M_1\cap M_2 = \widetilde{T}$. Denote the longitudinal slope of $X_i$ by $\lambda_i$ and let $\tilde \lambda_i$ be the slope on $\partial M_i$ that covers $\lambda_i$. Then $\tilde\lambda_i$ is the longitudinal slope of $M_i$ and hence is order-detected by Example \ref{exam: long is ord-det}. 

By assumption, one of $X_i$ is not a solid torus, say $X_2$.  We know that  the slope of $\tilde \lambda_1$ is order-detected in $M_1$ (Example \ref{exam: long is ord-det}), and if we can show that it is order-detected in $M_2$ it will follow from \cite[Theorem 1.3]{BC2} (see Theorem \ref{thm: * gluing}) that $\Sigma_n(L)$ has a left-orderable fundamental group. 

The irreducibility of $W$ implies that of $X_2$. Further, $\Delta(\lambda_1, \lambda_2) = 1$ since $W$ is an integer homology sphere and therefore $\lambda_1$ is order-detected in $X_2 \ne S^1 \times D^2$ by Theorem \ref{thm: meridional detn}(2). Indeed, the proof of Theorem \ref{thm: meridional detn}(2) shows that there is a homomorphism  $\rho: \pi_1(X_2) \rightarrow  \mbox{Homeo}_\mathbb Z(\mathbb R)$ for which $\rho(\lambda_1)$ has a fixed point in $\mathbb{R}$, but $\rho(\lambda_2)$ is fixed point free. (With a mild abuse of notation we have identified $\lambda_i$ with a primitive class in $\pi_1(T) \leq \pi_1(X_2)$ it carries.)

If $\varphi: \pi_1(M_2) \rightarrow \pi_1(X_2)$ denotes the homomorphism induced by the branched cover $M_2\rightarrow X_2$, then $\varphi|_{\pi_1(\widetilde T)}$ is injective. Indeed, there are integers $k_1, k_2 \geq 1$ so that $\varphi(\tilde \lambda_i) = \lambda_i^{k_i}$. It follows that if $\tilde \rho = \rho \circ \varphi: \pi_1(M_2) \to \mbox{Homeo}_\mathbb Z(\mathbb R)$, then $\tilde \rho(\lambda_1) = \rho(\lambda_1)^{k_1}$ has a fixed point in $\mathbb{R}$, but $\tilde \rho(\lambda_2) = \rho(\lambda_2)^{k_2}$ does not. Proposition \ref{prop: rep to ord det} implies that the slope carried by $\tilde \lambda_1$ is order-detected in $M_2$, which completes the proof. 
\end{proof}

\begin{proof}[Proof of Theorem \ref{thm: toroidal links nls}]
Let $L$ be a prime oriented link in an integer homology $3$-sphere $W$ whose exterior $X$ is irreducible and contains an essential torus $T$. If we can show that neither $\Sigma_2(L)$ nor $\Sigma_3(L)$ is an $L$-space, then an inductive application of Corollary 1.2 of \cite{HLL} shows that the same is true for $\Sigma_n(L)$ when $n = 2^k$, for some $k \geq 1$, or $3 \cdot 2^k$, for some $k \geq 0$, thus completing the proof. 

Fix $n \geq 2$ and suppose that $\Sigma_n(L)$ is an $L$-space. We show that $n \ne 2$ or $3$. 

We can assume that $L$ is non-split; otherwise $\Sigma_n(L)$ is not a rational homology $3$-sphere. 

By Theorem \ref{thm: nls detn zhs} we can suppose that $L$ has at least two components. Further, if $L$ is contained in one of the complementary components of $T$ we can argue as in the proof of Theorem \ref{thm: nls detn zhs} to see that $\Sigma_n(L)$ is not an $L$-space. Thus, at least one 
component of $L$ lies to each side of $T$. 

As $T$ is separating, we can write $W = M_1 \cup_T M_2$ where $M_1$ and 
$M_2$ are irreducible, boundary incompressible integer homology solid tori. The fact that $W$ is an integer homology $3$-sphere implies that there 
are simple closed curves $\lambda_1, \lambda_2$  on $T$ such that each $\lambda_i$ is nullhomologous in $H_1(M_i)$ and $\lambda_1, \lambda_2$ represents a basis of $H_1(T)$. Then $\lambda_i$ bounds an orientable surface $F_i$ 
in $M_i$ which we can suppose is transverse to $L$. By a mild abuse of notation we use $\lambda_i$ to denote the homology class in $H_1(T)$ it carries. 

If $L_i$ denotes the (non-empty) link $L \cap M_i$, the intersection $w_i$ of $L_i$ with $F_i$ can be thought of as either the winding number of $L_i$ in $M_i$ or the linking number of $\lambda_i$ and $L_i$. 

Fix $n \geq 2$ and let $\widetilde T$ be the inverse image of $T$ in $\Sigma_n(L)$ and $\widetilde M_i$ the inverse image of $M_i$. Since each $L_i$ is non-empty, both $\widetilde M_1$ and $\widetilde M_2$ are connected, and since $\Sigma_n(L) = \widetilde M_1 \cup_{\widetilde T} \widetilde M_2$, if $\widetilde T$ is disconnected then $\Sigma_n(L)$ is not a rational homology $3$-sphere, contrary to our assumptions. Thus $\widetilde T$ is connected and each $\widetilde M_i$ is a rational homology solid torus. Since $L$ is non-split and prime, $\Sigma_n(L)$ is irreducible and $\widetilde T$ is incompressible in $\Sigma_n(L)$ (Proposition \ref{prop: irr bccs}).  

Suppose that $n$ is prime and note that under the homomorphism $H_1(X) \to \mathbb Z/n$ which sends each (oriented) meridian of $L$ to 1 (mod $n$), $\lambda_i$ is sent to $w_i$ and therefore the connectivity of $\widetilde T$ implies $\gcd(w_1, w_2, n) = 1$. As $n$ is prime we can suppose, without loss of generality, that $\gcd(w_1, n) = 1$. Then the inverse image $\widetilde \lambda_1$ of $\lambda_1$ in $\widetilde T$ is also connected. Further, $\widetilde \lambda_1 = \partial \widetilde F_1$, where $\widetilde F_i$ is the $n$-fold cyclic branched cover of $(F_i, F_i \cap L_i)$, so $\widetilde \lambda_1$ is null-homologous in $\widetilde M_1$. Similarly $\widetilde \lambda_2$ is null-homologous in $\widetilde M_2$ when 
$\gcd(w_2, n) = 1$.

When $n$ is $2$ or $3$, $w_1 \equiv \pm 1$ (mod $n$) and up to reversing the orientation of $L$, which doesn't affect the proof, we can assume that  $w_1 \equiv 1$ (mod $n$). Further, either $w_2 \equiv 0$ (mod $n$) or $w_2 \equiv \pm 1$ (mod $n$). We show that both possibilities lead to a contradiction, thus completing the proof. 

If $w_2 \equiv 0$ (mod $n$), the inverse image of $\lambda_2$ in $\widetilde T$ consists of $n$ parallel curves whose  coherently oriented union is the boundary of $\tilde F_2$. Then one of them, $\widetilde \lambda_2$ say, represents the longitude of $\widetilde M_2$ and is therefore $NLS$-detected in $\widetilde M_2$. Since $\lambda_1, \lambda_2$ is 
a basis of $H_1(T)$, $\lambda_1 \cdot \lambda_2 = \pm 1$, and therefore as $\widetilde T \to T$ is an $n$-fold cover, $\widetilde \lambda_2 \cdot \widetilde \lambda_1 = 1$. Proposition \ref{prop: mu + k lambda NLS-detd} 
then implies that $\widetilde \lambda_2$ is $NLS$-detected in $\widetilde M_1$. But then $\Sigma_n(L) = \widetilde M_1 \cup_{\widetilde T} \widetilde M_2$ is not an $L$-space by Theorem \ref{thm: HRW1}, contrary to our 
assumptions. 

If $w_2 \equiv \varepsilon \in \{\pm 1\}$ (mod $n$), the inverse image $\widetilde \lambda_2$ of $\lambda_2$ is connected and is null-homologous in $\widetilde M_2$. As the distance of $\widetilde \lambda_1$ to $\widetilde \lambda_2$ is $n \in \{2, 3\}$, there is a  rational slope $\widetilde \beta$ on $\widetilde T$ of distance $1$ from both $\widetilde \lambda_i$. Then $\widetilde \beta$ is $NLS$-detected in both $\widetilde M_1$ and $\widetilde M_2$, by Proposition \ref{prop: mu + k lambda NLS-detd}. As in the previous case, this contradicts our assumption that $\Sigma_n(L)$ is an $L$-space, which completes the proof of Theorem \ref{thm: toroidal links nls}. 
\end{proof}

\subsection{Proof of Theorem \ref{thm: satellite not qa}}
\label{subsec: pf qa satellite implies thick}

Let $L$ be a prime $\mathbb Z/2$-Khovanov thin link. By Corollary \ref{cor: satellite not H-thin} we need only consider the case that the exterior of 
$L$ is Seifert fibred.
Such links are classified (\cite{BuMu}), and each can be oriented to be the closure of a positive braid. Since the $\mathbb Z/2$-Khovanov thinness of 
$L$ is independent of its orientation, we can assume that $L$ is a positive braid link and therefore strongly quasipositive. 

Our hypotheses imply that $\Sigma_2(L)$ is an $L$-space, and so $L$ is a definite link by Theorem 1.1 of \cite{BBG1}. The main result of \cite{Baa} then shows that $L$ is either 
\vspace{-.2cm}
\begin{itemize}

\item a $(2, m)$ torus link for some $m \geq 2$, or

\vspace{.2cm} \item one of the torus knots $T(3, 4)$ or $T(3,5)$, or 

\vspace{.2cm} \item a pretzel link $P(-2, 2, m)$ for some $m \geq 2$, or

\vspace{.2cm} \item the pretzel link $P(-2,3,4)$.

\end{itemize}
\vspace{-.2cm}
Since $(2, m)$ torus links are alternating, they are $\mathbb Z/2$-Khovanov thin. On the other hand, 
the torus knots $T(3, 4)$ and $T(3,5)$ are known to be $\mathbb Z/2$-Khovanov thick (\cite{Kh}, \cite{Sto}), as are the pretzel links $P(-2, 2, m)$, $m \geq 2$, and $P(-2,3,4)$ (\cite[Theorem 2.6]{Ma}) (these pretzel links 
had previously been shown not to be quasi-alternating by Greene (\cite[Proposition 2.2]{Gr})), which completes the proof.

\end{document}